\documentclass[opre,nonblindrev]{informs3} 

\DoubleSpacedXI 

\usepackage{endnotes}
\let\footnote=\endnote

\usepackage{natbib}
 \bibpunct[, ]{(}{)}{,}{a}{}{,}%
 \def\bibfont{\small}%
\TheoremsNumberedThrough     
\ECRepeatTheorems

\EquationsNumberedThrough  

\usepackage{enumerate} 
\usepackage{tikz}
\usetikzlibrary{arrows,
                petri,
                topaths}
\definecolor{mycolor}{RGB}{215,255,244}

\usepackage{multirow}
\usepackage[normalem]{ulem}
\usepackage{makecell}
\usepackage[colorlinks,linkcolor=blue,citecolor=cyan,anchorcolor=blue]{hyperref}
\newcommand{\R}{\mathbb{R}}
\newcommand{\D}{\mathbb{D}}

\newcommand{\Z}{\mathbb{Z}}
\newcommand{\E}{\mathbb{E}}

\newcommand{\III}{\mathcal{I}}
\newcommand{\SSS}{\mathcal{S}}
\newcommand{\JJJ}{\mathcal{J}}

\newcommand{\Prob}{\mathbb{P}}

\providecommand{\abs}[1]{\left\lvert#1\right\rvert}
\providecommand{\norm}[1]{\lVert#1\rVert}

\newcommand{\blue}[1]{#1}
\newcommand{\delete}[1]{}
\newcommand{\bluetwo}[1]{#1}
\newcommand{\deletetwo}[1]{}

\begin{document}

\RUNAUTHOR{Braverman, Dai, Liu and Ying}
\RUNTITLE{Empty-car Routing in Ridesharing Systems}
\TITLE{Empty-car Routing in Ridesharing Systems}
\ARTICLEAUTHORS{%
\AUTHOR{Anton Braverman}
\AFF{School of Operations Research and Information Engineering, Cornell University, Ithaca, New York 14853, \EMAIL{ab2329@cornell.edu}} 
\AUTHOR{J.G. Dai}
\AFF{School of Operations Research and Information Engineering, Cornell University, Ithaca, New York 14853, \EMAIL{jim.dai@cornell.edu}} 
\AUTHOR{Xin Liu}
\AFF{School of Electrical, Computer and Energy Engineering Arizona State University Tempe, Arizona 85287, \EMAIL{xliu272@asu.edu}} 
\AUTHOR{Lei Ying}
\AFF{School of Electrical, Computer and Energy Engineering Arizona State University Tempe, Arizona 85287, \EMAIL{lei.ying.2@asu.edu}} 
} 

\ABSTRACT{%
This paper considers a closed queueing network model of ridesharing systems such as Didi Chuxing, Lyft, and Uber. We focus on empty-car routing, a mechanism by which we control car flow in the network to optimize  system-wide utility functions,  e.g. the availability of empty cars when a passenger arrives. 
We establish both process-level and steady-state convergence of the queueing network to a fluid limit in a large market regime where demand for rides and supply of cars tend to infinity, and use this limit to study a fluid-based optimization problem.
We prove that the optimal network utility obtained from the fluid-based optimization is an upper bound on the utility in the finite car system for any routing policy, both static and dynamic, under which the closed queueing network has a stationary distribution. This upper bound is achieved asymptotically under the fluid-based optimal routing policy. Simulation results with real-\blue{world} data released by Didi Chuxing demonstrate the benefit of using the fluid-based optimal routing policy compared to various other policies.
}%



\maketitle

%


\section{Introduction}
\label{sec:introduction}
This paper studies the modelling and control of ridesharing systems such as Didi Chuxing, Lyft and Uber. We consider a system with $r >0$ regions and $N>0$ cars. The regions can be interpreted as geographic regions in a city and cars drive around between regions transporting passengers. At time $t = 0$, all cars start off idling empty in some region, waiting for a passenger. Passengers arrive to region $i$ according to a Poisson process with rate $N\lambda_i > 0$, and arrivals to different regions are independent. When a passenger arrives to region $i$, if there is an empty car available there, then the passenger occupies that car and travels to region $j$ with probability $P_{ij}$. If no empty car is available, the passenger abandons the system and finds an alternative form of transportation to her destination. We allow $P_{ii}> 0$ to represent trips within a region. Travel times from region $i$ to $j$ have mean $1/\mu_{ij}$ and are assumed to be i.i.d.\ exponential random variables, although this assumption is not essential (see Remark~\ref{rem:general_travel} in Section~\ref{sec:main_results}).
Once the passenger arrives at region $j$, the car becomes empty. The empty car can either \deletetwo{decide to} stay in region $j$ with probability $Q_{jj}$ \bluetwo{(it becomes available to take new passengers immediately)}, or with probability $Q_{jk}$, \bluetwo{relocate without a passenger} to a different region $k$ and wait for a passenger there. \blue{The time spent driving empty from $j$ to $k$ is identical in distribution to that of driving with a customer.} In general, the routing matrix (also called the routing policy) $Q= (Q_{ij})$ is allowed to be \emph{state-dependent}, i.e.\ $Q$ may depend on the current distribution of cars across the regions. \bluetwo{In this paper, $Q$ will be a decision variable.}

We model this ridesharing network with a closed queueing network consisting of both single-server and infinite-server stations, where cars are ``jobs'' moving through the queueing network. Cars waiting in a region for a passenger are modeled with a single-server station. The buffer content of the station corresponds to the number of cars waiting, and passenger ride requests correspond to service completions at the station. Thus, the service times at the single-\blue{server} station are the interarrival times of passengers to the region, although there is no physical ``server'' at the station. The infinite-server stations are used to model car travel between regions. When the routing policy $Q$ is static, i.e.\ not state dependent, our queueing network belongs to a class of networks called BCMP networks \cite{BaskChanMuntPala1975}. The precise formulation of the model can be found in Section~\ref{sec:model}.



\begin{figure}[h!]
\begin {center}
\begin {tikzpicture}[auto, node distance=4 cm and 5cm,
semithick, state/.style ={ circle, top color=white, draw, text=black , minimum width = 1.1 cm, line width = 0.25 mm}]
\tikzset{>=latex}
\tikzstyle{dashed}=                  [dash pattern=on 3pt off 3pt]
\tikzstyle{densely dashed}=          [dash pattern=on 3pt off 2pt]
\tikzstyle{loosely dashed}=          [dash pattern=on 4pt off 5pt]
\node[state] (A) at (0,0) {1};
\node (Atext) at (0,-1.75) {$800$ passengers/unit time};
\draw[->] (Atext) -- (0, -0.75);
\node[state] (B) at (4,0) {2};
\node (Btext) at (4,1.75) {$400$ passengers/unit time};
\draw[->] (Btext) -- (4, 0.75);
\draw[thick] (0.75,0.7) -- (3.25,0.7);
\draw[thick,->] (1.4,0.35) -- (2.6,0.35);
\draw[thick, loosely dashed] (0.75,0) -- (3.25,0);
\draw[thick,->] (2.6,-0.35) -- (1.4,-0.35);
\draw[thick] (0.75,-0.7) -- (3.25,-0.7);
\end{tikzpicture}
\end{center}
\caption{A two-region example}
\label{fig:example}
\end{figure}
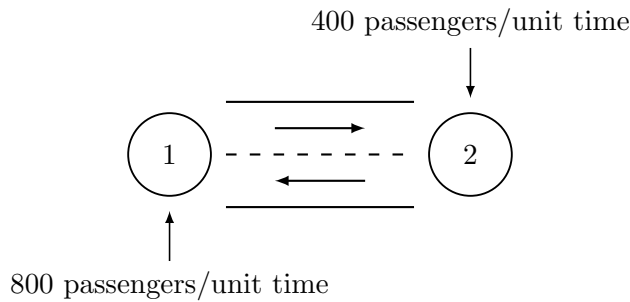

Because of the proliferation of ridesharing and bikesharing services, modelling and control of these systems have become important research topics over the last few years \cite{Adel2007, GeorXia2011, JostWase2013, JostWase2016, BaneFreuLyko2016, IglePavoRossZhan2016, FrazPavoRusSmit2012, PavoZhan2016, OzkaWard2016, FrazIyerYang2016, BimpCandDani2016}. Our paper focuses on empty-car routing as a mechanism to improve the efficiency of the system. To illustrate the effect of this mechanism, consider the two-region example in Figure \ref{fig:example}. Passengers arrive at region $1$ to go to region $2$ according to a Poisson process with rate $800$ passengers/unit time, and arrive at region $2$ to go to region $1$ with rate $400$ passenger/unit time. After dropping off a passenger at region $2$, a driver stays at region $2$ with probability $Q_{22}$, or drives empty to region $1$ with probability $Q_{21}$. The probabilities $Q_{11}$ and $Q_{12}$ are defined analogously. We define the \emph{availability} at region $i$ to be the long-run fraction of time that there is at least one empty car at the region available. Since Poisson arrivals see time averages (PASTA property), this is also the probability that a passenger's request for a ride originating from region $i$ will get fulfilled.

A number of existing models consider \emph{one-way} vehicle sharing systems, in which a vehicle can only be moved from one region to another when carrying a passenger \cite{GeorXia2011, JostWase2016, BaneFreuLyko2016}. This is a realistic assumption for bikesharing systems, where a bicycle cannot move autonomously from one region to another, and only moves when a passenger rides it. In such a case, the performance of the system is largely determined by the passengers' arrival rates and destination probabilities.  In our example, a one-way system would correspond to $Q_{12}=Q_{21}=0$, and a bike taken from $1$ to $2$ will only return to $1$ if it is brought there by a passenger. Since, on average, region $1$ sees twice as many passengers as region $2$,  the availability of bikes at region $1$ will always be at most $50\%$, regardless of the number of bikes in the system. That is, region $1$ will lose at least half of its passengers to alternative modes of transportation. This inefficiency due to passenger imbalance has been well recognized
in the literature. Proposed solutions include demand throttling via
pricing \cite{JostWase2016, BaneFreuLyko2016}, or periodic bike
rebalancing using trucks \cite{ChemMeunWolf2013,
  HendOmahShmo2016}.

Empty-car routing is an appropriate mechanism for commercial ridesharing systems, where drivers often wander around to find passengers, and it is not surprising that a good routing policy can increase the efficiency of the system. Returning to our example, we assume that there are a
total of $1200$ cars in the system, and that the mean travel time in either direction is one unit. Table~\ref{tab:2routing} compares the availability in each region under several different routing policies. In particular, we see that the policy with $Q_{21}=1/3$ is preferable to the policy with no empty-car routing ($Q_{12}=Q_{21}=0$). \bluetwo{ The question we answer in this paper is the following. Given a utility function measuring performance in the system (average availability for example), how should one choose a routing policy to maximize said utility.} 
\begin{table}[htb]
  \centering
\begin{center}
  \begin{tabular}{c|cc}
Empty-car     &     \multicolumn{2}{c}{Availability} \\
               \cline{2-3}
routing policy  &  Region $1$  & Region $2$ \\
\hline
$Q_{12}=0,\ Q_{21}=1/3$    &     73.19\%        & 97.59\% \\
$Q_{12}=0,\ Q_{21}=1/2$    &     74.64\%        & 74.64\% \\
$Q_{12}=Q_{21}=0$   &     50\%        & 100\%
  \end{tabular}
\end{center}
\caption{Availabilities under several empty-car routing policies with $1200$ cars in the system computed using the MVA algorithm.}
  \label{tab:2routing}
\end{table}

Recall that with a static routing policy $Q$, our queueing network is a BCMP network. The stationary distribution of a BCMP network has a product form, but the normalization constant is expensive to compute because the state space of the network is too large. However, algorithms such as mean value analysis (MVA) \cite{LaveReis1980} or approximate mean value analysis (AMVA) \cite{SahuSuri2007} can bypass computing the normalization constant, and directly compute performance metrics of interest, e.g.\ mean queue sizes. In other words, given $N,$ $\lambda, \mu,$  $P,$ and a static $Q$, steady-state performance analysis of the system can be done efficiently. However, the problem of optimizing some performance metric over $Q$ is difficult.  Even in the $2$-region case, \bluetwo{it can be verified numerically that}
\begin{align}
\max_{Q} \quad  \lambda_1 \times \text{(region $1$ availability) + $\lambda_2$ $\times$ (region $2$ availability)} \label{eq:example_opt}
\end{align}
is a non-convex optimization problem. \bluetwo{The reason for this is that availabilities have a non-linear dependence on $Q$.}
When we tried solving \eqref{eq:example_opt}, MATLAB's built-in solver failed to converge to a solution. We also tried using NEOS \cite{NEOS}, which is a collection of more sophisticated optimization tools. However, even the solvers in NEOS could not reliably solve the problem.


The main contribution of this paper is provide comprehensive solution to the empty-car routing problem, which is both theoretically grounded and efficient to implement. Namely, we show that as the number of cars and the passenger arrival rates tend to infinity, i.e.\ $N \to \infty$, the optimal solution of \eqref{eq:example_opt} converges to the optimal solution of a fluid-based optimization problem that can be solved by solving a linear program. Our results also hold for much broader class of utility functions; cf. Remark~\ref{rem:general_utility} in Section~\ref{sec:efficient}. We also show in Theorem~\ref{thm:lb} that the performance under the optimal static routing policy coming from the fluid-based optimization problem is an upper bound on the performance under any  \emph{state-dependent} routing policy. Furthermore, this upper bound is asymptotically tight.  \blue{  For any stochastic control problem of realistic size, the true optimal policy is rarely known. Thus, any upper bound, particularly a tight bound on the performance is valuable to develop good policies for finite-sized systems. }

The typical asymptotic regime considered for a closed BCMP network is one where the number of jobs in the network goes to infinity, but the service rates at each station remain fixed. In the context of our ridesharing network, this would correspond to a regime where the number of cars $N$ increases to infinity, while both passenger request rates and travel times remain fixed. A lot is known about the asymptotic behavior of BCMP networks in this asymptotic regime, see for instance \cite[Section 4.1]{Geor2012} and the references within. Most importantly, the limiting network always has at least one region with an infinite number of empty cars, i.e.\ a region where availability equals one. \delete{Using queueing jargon, this is equivalent to the bottleneck stations in the network having $100\%$ utilization (asymptotically).}  For this reason, we refer to this as the \emph{infinite supply} regime. The asymptotic regime considered in this paper has both the number of cars $N$, and the passenger arrival rates $N \lambda$ going to infinity together. We refer to this regime as the \emph{large market} regime. The infinite supply and large market regimes are qualitatively different. The latter is not guaranteed to have a region where availability converges to $100\%$; we discuss the implications of this further in Section~\ref{sec:related}.

In practice, the large market regime is more realistic than the infinite supply regime. For starters, it is natural that the supply of drivers in a city increases with demand for rides. Furthermore, the large market regime does not impose any restrictions on supply-demand imbalance in a city. The parameter $\lambda_i$ is the rate of arriving passengers per car to region $i$, and gives our model the flexibility to distinguish between cases when there is an oversupply, undersupply, or critical level of supply of cars with respect to passenger demand. The latter two cases are more common during morning or afternoon rush hours, and are precisely the cases where an effective choice of routing matters the most. \deletetwo{Figure~\ref{fig: didi} shows real-life ridesharing traffic data released by Didi Chuxing in January 2016. We see there the total number of passenger orders compared to the total number of accepted orders of two nearby regions during the afternoon rush hour (5PM to 6PM) over 21 consecutive days. From the figure, we can observe that there was a significant shortage of cars in region 47, pointing to a potential undersupply of cars in the system. Furthermore, there was a great disparity in car availability between the two regions, which illustrates the potential benefit of using an efficient car-routing policy.}


\blue{Our main results are stated in the steady-state setting, assuming the system parameters stay constant. However, in practice it is very common for parameters to depend on the time of day, e.g.\ passenger arrivals spike during rush hours. In Section~\ref{sec:lookahead}, we leverage our fluid-based optimization problem to suggest a time-dependent lookahead policy that  anticipates future parameter changes and routes cars accordingly.  We present several numerical examples with time-varying parameters where this anticipative behavior yields significant performance benefits over the typical approach of dividing time into smaller intervals and assuming constant parameters on each interval. We remark that the purpose of this lookahead heuristic is to demonstrate that the fluid-based optimization can guide the design of high performance empty-car routing policies in practical ridesharing systems where some of the modeling assumptions made in this paper may not hold. However, this lookahead is not meant to the ``optimal'' algorithm for these settings. }

Before moving on to the literature review, we wish to say a few things about the routing mechanism in this model. Our model uses a centralized mechanism for routing, i.e.\ the ridesharing company generates routing decisions according to the routing matrix $Q$, and cars then have to obey that decision. This mechanism is perfectly fine for systems with autonomous vehicles, which are already experimented with by Uber \cite{Chaf2016}. However, a model with centralized control is still useful even when human drivers are free to make their own decisions. For instance, a centralized control mechanism is a best-case performance benchmark against which one can compare decentralized mechanisms. 
A centralized mechanism can also be used to quantify the ``price of anarchy'', i.e.\ the difference in revenue between  centralized control and the case when drivers are free to make their own decisions. Part of this cost can then be used to incentivize drivers to obey routing instructions, e.g.\ by subsidizing fuel costs for driving empty.


\subsection{Related Literature} \label{sec:related}
%

BCMP networks are natural choice for modelling ridesharing systems. The closest papers to us are \cite{PavoZhan2016, IglePavoRossZhan2016}, where the authors also consider the supply repositioning problem. BCMP networks have also been used to study fleet-sizing \cite{GeorXia2011} and pricing \cite{JostWase2016} problems in one-way vehicle systems, where no resource repositioning takes place.



In \cite{PavoZhan2016} and \cite{IglePavoRossZhan2016}, authors use an empty-car routing mechanism to rebalance their ridesharing networks. However, each paper focuses on a single optimization problem, whereas our technical approach is more robust and allows us to consider a large class of optimization problems. Like us, \cite{IglePavoRossZhan2016} also faces a non-linear optimization problem. To deal with it, the authors pass to the \emph{infinite supply} regime where they obtain a simpler optimization problem that they solve and apply as a heuristic to the original problem from the finite system. The connection between the finite system and the limiting system in the infinite supply regime received rigorous treatment in \cite{BaneFreuLyko2016}, where the authors establish bounds on the gap between the optimal values of the optimization problems of the finite and infinite supply systems.

The aforementioned papers \cite{BaneFreuLyko2016, PavoZhan2016, IglePavoRossZhan2016} all share the same feature that the optimization problems of those papers equalize availabilities across all regions. This is either enforced via an explicit constraint \cite{IglePavoRossZhan2016}, \cite[Equation 10]{PavoZhan2016}, or arises implicitly in the optimal solution of the approximating optimization problem in \cite{BaneFreuLyko2016}. The rationale behind enforcing the equal availability constraint is that in the infinite supply regime, at least one region achieves 100\% availability. Equalizing availabilities among regions then ensures 100\% availability everywhere. Therefore, the solutions proposed in those papers rely on an abundance of vehicle supply in the system, and are unsuitable for problems where demand is comparable to, or even exceeds, supply. Indeed, in our $2$-region example in Figure~\ref{fig:example}, the only way to equalize availabilities is to use the routing policy $Q_{12}=0$, and $Q_{21} = 1/2$, and we see in Table~\ref{tab:2routing} that this hurts system performance.

\blue{After submission of this paper, a newer version of  \cite{BaneFreuLyko2016} appeared, namely \cite{BaneFreuLyko2017}. In both versions, the authors study approximations of optimization problems arising in finite sized ridesharing systems. The authors directly analyze the stationary distribution of BCMP networks to establish asymptotically tight guarantees for their approximating optimization problems. In \cite{BaneFreuLyko2016}, the focus is on pricing as a control mechanism in the infinite supply regime. Motivated by the asymptotic regime and control mechanism considered in this paper, the new version \cite{BaneFreuLyko2017} proves that in our Theorem~\ref{thm:lb}, \bluetwo{ the ratio between the optimal value of the optimization problem for the finite-sized system and the optimal value of the fluid-based optimization goes to one at a rate of $1 - 1/\sqrt{N}$;} see appendix D.1 there.}


Fluid models have been used by \cite{FrazPavoRusSmit2012, JostWase2013} to study ridesharing networks. Those fluid models are different from the fluid model in this paper. Furthermore, they are only used as heuristics, and are not shown to be connected to an underlying stochastic system. In contrast, our fluid model is proven to be the limit of the queue length process of our BCMP network, and all fluid optimization problems considered are rigorously proven to be the limits of their stochastic counterparts.

We also wish to highlight some papers that either consider problems closely related to ridesharing, or address some of the other aspects of the ridesharing problem not focused on in this paper. In \cite{Adel2007}, the author considers the problem of managing a network of shipping containers. That paper uses a combination of optimization and approximate dynamic programming techniques to determine a policy to accept/reject requests for containers. Both \cite{MaWolfZhen2013, SantXavi2015} study a carpooling problem of having multiple riders with different destinations share the same car. In \cite{OzkaWard2016}, the authors adopt a matching approach to the setting where a passenger requesting a ride from a region with no available cars is willing to wait a little bit for a driver to arrive from a nearby region. They study the problem of matching drivers to passengers when the two may have two different initial locations. In \cite{BimpCandDani2016}, the authors consider the issue of pricing rides in a ridesharing network. Also related to the ridesharing problem is \cite{FrazIyerYang2016}, where the authors study a mean field equilibrium of a system where agents explore and compete for resources that are both time-varying and location-dependent in nature. 

Outside of the ridesharing setting, fluid models are a widely used tool in the study of closed queueing networks; we refer the reader to \cite{AnseDaurWalt2013} for a recent discussion of the literature. For the type of network considered in this paper, i.e.\ a closed network with single-server and infinite-server stations, process level convergence to the fluid model was established in \cite{Kric1992}. The technical novelty of our paper lies in Theorem~\ref{thm:stability}, which characterizes the limiting behavior of the fluid model. A related paper is \cite{AnseDaurWalt2013}, which characterizes the limiting behavior of the fluid model corresponding to a multiclass, closed queueing network consisting entirely of single-server stations, with no infinite-server stations. The results in \cite{AnseDaurWalt2013} were established using a relative entropy based Lyapunov function. In this paper, we use the $L_1$ distance from the equilibrium as our Lyapunov function.

\subsection{Contributions}
The following is a summary of the main contributions of this paper.
\begin{itemize}
\item Fixing a static empty-car routing policy
    $Q=(Q_{ij})$, we consider a fluid model associated with a closed
  queueing network composed of single and infinite server stations. We
  establish process level convergence of the scaled queue length
  process in our closed queueing network to a fluid limit.  The fluid
  model's equilibrium set is explicitly characterized, and we show
  that the fluid model converges to this equilibrium set from any
  initial starting condition. We then elevate the process level
  convergence result to convergence of steady-state distributions, cf.\
  Section~\ref{sec:fluid}. \bluetwo{One consequence of our result is an answer to the following question: what is the minimum number of cars in the system needed to achieve $100\%$ availability (asymptotically as $N \to \infty$) everywhere? This discussion can be found in Appendix~\ref{sec:perfectavailability}.}

\item To find an optimal static empty-car routing policy $Q^*$, we formulate a fluid-based optimization problem that is able to accommodate a broad class of utility functions. The latter can depend on availabilities at different regions, and fractions of both empty or occupied cars on different roads. Then $Q^*$ can be solved efficiently by solving a related problem with only linear constraints, cf. Lemma~\ref{lem:relax}.

\item We prove in Theorem~\ref{thm:lb} that as the number of car grows to infinity,
  the routing policy $Q^*$  from the fluid-based optimization is asymptotically optimal among
  all state dependent routing  policies. \blue{For any MDP or stochastic control problem of realistic size, the true optimal policy is rarely known. Thus, any upper bound, particularly a tight bound is valuable to developing good policies.}


\delete{ 
\item Using real-world network and passenger order data from a dataset
  released by Didi Chuxing, we simulated the performance of the
  proposed routing policy $Q^*$.  We observe from simulation
    estimates that (a.) actual availability converges to the
    fluid-based solution at rate $1/\sqrt{N},$ (b.)  the static
  routing policy $Q^*$ outperforms heuristic state-dependent routing policies.}
\end{itemize}

The rest of the paper is structured as follows. In Section~\ref{sec:model}, we formulate the fluid-based optimization problem and state our main results, Theorems~\ref{thm:asymptotic_utility} and \ref{thm:lb}. In Section~\ref{sec:numerical}, we describe the numerical study performed using real-world data from Didi Chuxing, China's largest ridesharing company. Section~\ref{sec:fluid} is devoted to studying the fluid model of the ridesharing network, and establishing the machinery needed to prove our main results. Section~\ref{sec:conclusion} concludes.

\subsection{Notation}
For a function $f: \R \to \R^n$, we use $\dot f(t)$ to denote the derivative of $f(t)$ when the derivative exists. For any integer $n > 0$, we use $\D^n$ to denote the space of all cadlag functions $x: \R_+ \to \R^n$, i.e. functions that are right-continuous on $[0, \infty)$ with left limits on $(0, \infty)$. We define
\begin{align*}
\D^n_0 =&\ \{ x \in \D^n :\ x(0) = 0\}, \\
\D^n_1 =&\ \Big\{ x \in \D^n :\ x(0) \in [0,1]^n \text{ and } \sum_{i=1}^{n} x_i(0) = 1 \Big\}, \\
\D^n_{0+} =&\ \{ x \in \D^n :\ x(0) \geq 0\}.
\end{align*}
For any $x \in \D^n$ and any $T > 0$, we define
\begin{align}
\norm{x}_T = \max_{1 \leq i \leq n} \sup_{0 \leq t \leq T} \abs{x_i(t)} = \sup_{0 \leq t \leq T} \max_{1 \leq i \leq n} \abs{x_i(t)}. \label{eq:norm}
\end{align}
We let $C^n \subset \D^n$ be the subspace of continuous functions $x: \R_+ \to \R^n$, and define $C^n_1$ analogously to $\D^n_1$. For any $x \in \D^n$, we write $\int_{0}^{t} x(s) ds$ to denote a vector in $\R^n$ whose $i$th component is $\int_{0}^{t}x_i(s) ds$. For a vector $a \in \R^n$, we use $\abs{\cdot}$ to denote the max-norm, i.e. $\abs{a} = \max_{1 \leq i \leq n} \abs{a_i}$. For a set $A \subset \Z$, we write $\abs{A}$ to denote the number of elements contained by this set. For a collection of random variables $Y, \{X_n\}_{n=1}^{\infty}$, we write $X_n \Rightarrow Y$ to denote weak convergence of $X_n$ to $Y$ (as $n \to \infty$).

\section{The Ridesharing Optimization Problem} \label{sec:model} In
this section we formally introduce the sequence of ridesharing
networks discussed in the introduction. We then introduce the
fluid-based optimization and state our main results,
Theorems~\ref{thm:asymptotic_utility} and \ref{thm:lb}. We then show in Lemma~\ref{lem:relax} that the fluid-based optimization can be solved efficiently by solving a related optimization problem with linear constraints.

In our model, there are $N$ cars serving $r$ regions in a city. For any time $t \geq 0$, let $E^{(N)}_{ij}(t)$  be the number of empty cars en route from region $i$ to region $j \neq i$, and let $E^{(N)}_{ii}(t)$ be the number of empty cars that are waiting in region $i$ for a new passenger. Similarly, let $F^{(N)}_{ij}(t)$ be the number of full cars driving from region $i$ to $j$ (by full car we mean a car with a passenger). We allow $F^{(N)}_{ii}(t)$ to be non-zero, because a passenger's destination can be located in the same region as he was picked up. Let $E^{(N)}(t)$ and $F^{(N)}(t)$ be the $r \times r$ matrices whose $(i,j)$th elements are $E^{(N)}_{ij}(t)$ and $F^{(N)}_{ij}(t)$, respectively. 
\begin{align}
&E^{(N)} = \{E^{(N)}(t) \in \Z_+^{r \times r},\ t \geq 0\} \quad \text{ and } \quad F^{(N)} = \{F^{(N)}(t) \in \Z_+^{r \times r},\ t \geq 0 \}, \notag \\
&\bar E^{(N)} = \Big\{\frac{1}{N} E^{(N)}(t) \in \R_+^{r \times r},\ t \geq 0 \Big\} \quad \text{ and } \quad \bar F^{(N)} = \Big\{\frac{1}{N} F^{(N)}(t) \in \R_+^{r \times r},\ t \geq 0  \Big\}, \notag \\
&\mathcal{T} = \bigg\{(e, f) \in [0,1]^{r\times r} \times [0,1]^{r\times r} : \sum_{i=1}^{r} \sum_{j=1}^{r} (e_{ij} + f_{ij}) = 1\bigg\}. \label{eq:T}
\end{align}
Recall the dynamics introduced in Section~\ref{sec:introduction}.  \deletetwo{that passengers arrive to region $i$ according to a Poisson process with rate $N\lambda_i$, and that travel times between regions $i$ and $j$ are exponentially distributed with mean $1/\mu_{ij}$. \blue{We can interpret $\lambda_i$ as the number of passengers arriving per time unit, per car.}}  We allow the empty-car routing probabilities to be state-dependent, i.e.  the empty-car routing probability matrix at time $t \geq 0$ is
\begin{align*}
Q(\bar E^{(N)}(t) ,\bar F^{(N)}(t)) = \Big(Q_{ij}(\bar E^{(N)}(t) ,\bar F^{(N)}(t)) \Big).
\end{align*}
\deletetwo{where, for each $(e,f)\in \mathcal{T}$, $Q(e,f)$ is an $r\times r$ stochastic matrix.}
The process $(E^{(N)}, F^{(N)})$ is then a continuous time Markov chain (CTMC), whose transition properties are listed in Table~\ref{tab:transitions}.
\begin{table}[h!]
\caption{Markov Chain Transition Rates}
\centering
\label{tab:transitions}
\begin{tabular}{|c|c|l|}
\hline
 \multicolumn{1}{|c|}{Rate} & Transition      & \multicolumn{1}{c|}{Description} \\ \hline
\multirow{2}{*}{$N\lambda_i P_{ij} 1(E^{(N)}_{ii}(t) > 0)$} & $E^{(N)}_{ii}(t) - 1,$ &  Passenger arrives to region $i$,\\ 
& $F^{(N)}_{ij}(t) + 1$  & and starts ride from $i$ to $j$. \\ \hline
 \multirow{3}{*}{$\mu_{ij}F^{(N)}_{ij}(t) Q_{jk}(\bar E^{(N)}(t) ,\bar F^{(N)}(t))$} & $F^{(N)}_{ij}(t) - 1,$  & Passenger dropped off at region $j$. Car\\
& $ E^{(N)}_{jk}(t) + 1$ &  stays put if $k = j$, or starts driving \\
& &  empty to region $k$ if $k \neq j$\\ \hline 
 \multirow{2}{*}{$\mu_{ij}E^{(N)}_{ij}(t) 1(j \neq i)$} & $E^{(N)}_{ij}(t) - 1,$  & Empty car arrives to region  $j$. Stays\\
& $ E^{(N)}_{jj}(t) + 1$ &  there until next passenger. \\ \hline 
\end{tabular}
\label{tab:notations}
\end{table}
\deletetwo{Clearly, the fluid-scaled version $(\bar E^{(N)},\bar F^{(N)})$ is also a CTMC.}  Going forward, we focus on the fluid-scaled CTMC $(\bar E^{(N)},\bar F^{(N)})$ and assume it has a single recurrent class. Since the state space is finite, this implies that the CTMC is positive recurrent. Let $(\bar E^{(N)}(\infty) , \bar F^{(N)}(\infty)) \in \mathcal{T}$ be the random element having the stationary distribution of $(\bar E^{(N)}, \bar F^{(N)})$. Let
\begin{align*}
A_{i}^{(N)}(\infty) = \Prob( \bar E^{(N)}_{ii}(\infty) > 0)
\end{align*}
be the availability at region $i$, and let $A^{(N)}(\infty)$ be the $r$-dimensional vector with entries $A_{i}^{(N)}(\infty)$.

 \blue{ Our model assumes that a routing decision is made, the car is committed to the region until a passenger arrives to take it away (as opposed to being able to hop around empty between regions). This assumption is reasonable in a setting with a centralized router, where the decision maker has a global view of the system and will therefore `get it right the first time' when making the routing decisions. The incentive for cars to keep jumping around the network is further reduced by the assumption that passengers will not linger in the system if they cannot receive service immediately. This passenger behavior is realistic in settings when comparable transportation modes are available, e.g.\ hailing a yellow cab in Manhattan, or taking public transportation.}

\begin{remark}
The process $(E^{(N)},F^{(N)})$ can also be interpreted as the queue length process in a closed queueing network
of $r$ single server stations and $2r^2-r$ infinite server stations, where cars are the ``jobs'' in the network. For $1 \leq i \leq r$, the process $E^{(N)}_{ii} = \{E^{(N)}_{ii}(t),\ t \geq 0\}$ corresponds to a single server station with service rate $N \lambda_i$, and
\begin{align*}
E^{(N)}_{ij} =&\ \{E^{(N)}_{ij}(t),\ t \geq 0\}, \quad 1 \leq i \neq j \leq r,\\
F^{(N)}_{ij} = &\ \{F^{(N)}_{ij}(t),\ t \geq 0\}, \quad 1 \leq i,j \leq r,
\end{align*}
correspond to infinite server stations where the service rate of each server at station $E^{(N)}_{ij}$ or $F^{(N)}_{ij}$ is $\mu_{ij}$.  \blue{In the special case when $Q$ is not state-dependent}, this network \blue{belongs to the class of BCMP networks} \cite{BaskChanMuntPala1975}. \blue{Note that in our model, the service time at a single-server station can begin even before a job enters the station, e.g.\ the inter-arrival ``timer'' of the next passenger is counting down regardless of whether there is an available car in the region or not. However, the memoryless property of passenger inter-arrival times makes our model equivalent to one where service starts only when the station is non-empty.}
\end{remark}
We are now ready to introduce the fluid-based optimization problem, and state our main results.

\subsection{Main Results} \label{sec:main_results}
Recall the network primitives $\lambda, \mu, P$. We now consider the fluid-based
optimization problem to be fully specified from (\hyperref[eqlbm:objective]{4a}) to (\hyperref[eqlbm:prob]{4g}) below. In the
optimization problem, $c_{ij} > 0$ are rewards for completing a ride from $i$ to $j$. The variables in the optimization problem are $q$,
$\bar e$, $\bar f$, $\bar a$, where $q = (q_{ij})$ is an $r \times r$
matrix representing a static empty-car routing policy $Q$ and $(\bar
e, \bar f, \bar a)$ is a point in $ \mathcal{T} \times [0,1]^{r}$,
whose interpretation will be given after the equations
(\hyperref[eqlbm:objective]{4a})-(\hyperref[eqlbm:prob]{4g}):
\begin{subequations}
\begin{align}
&\max_{q,\bar e, \bar f, \bar a} \sum_{i=1}^{r} \sum_{j=1}^{r} \bar a_i\lambda_i P_{ij} c_{ij}  \label{eqlbm:objective} \\
\text{ subject to } \quad &\lambda_iP_{ij}\bar a_i = \mu_{ij}\bar f_{ij}, \quad 1\leq i,j \leq r, \tag{4b: full car i-j Little's Law} \label{eqlbm: fullcsv}\\
&\mu_{ij}\bar e_{ij} = q_{ij}\sum_{k=1}^{r}\mu_{ki}\bar f_{ki}, \quad 1\leq i,j \leq r, ~j\neq i, \tag{4c: empty car i-j Little's Law} \label{eqlbm: emptycsv1}\\
&\lambda_i\bar a_i = \sum_{k=1,k\neq i}^{r}\mu_{ki}\bar e_{ki} + q_{ii}\sum_{k=1}^{r}\mu_{ki}\bar f_{ki}, \quad 1\leq i \leq r, \tag{4d: car flow balance region i}  \label{eqlbm: emptycsv2}\\
&(1-\bar a_i)\bar e_{ii} = 0, \quad 1\leq i \leq r, \tag{4e: availability to idle car relation} \label{eqlbm: boundry}\\
&\blue{(\bar e, \bar f) \in \mathcal{T},} \quad \tag{4f: unit mass} \label{eqlbm:carnum} \\
&q_{ij}\geq 0, \quad \sum_{j=1}^{r}q_{ij} = 1, \quad 0 \leq \bar a_{i} \leq 1, \quad 1\leq i,j \leq r \tag{4g} \label{eqlbm:prob} 
\end{align}
\end{subequations}
\bluetwo{To help guide intuition, one can think of $\bar e, \bar f$, and $\bar a$ as placeholders for $\E[\bar E^{(N)}(\infty)], \E[\bar F^{(N)}(\infty)]$ and $A^{(N)}(\infty)$, respectively.} \blue{We can interpret $\bar a_i \lambda_i P_{ij}$ as the rate at which rides are initialized from $i$ to $j$, and since a ride from $i$ to $j$ has a reward of $c_{ij}$, the problem above aims to maximize revenue generation.} Our results actually hold for a much larger class of utility functions; cf.\ Remark~\ref{rem:general_utility} in Section~\ref{sec:efficient}. The constraints in \blue{(\hyperref[eqlbm: fullcsv]{4b}) are simply Little's Laws for the number of occupied cars on the road from $i$ to $j$ in equilibrium. That is, $\lambda_iP_{ij}\bar a_i$ is the rate at which rides are initialized, which equals  the mass of occupied cars on the road $f_{ij}$ divided by  the average travel time $1/\mu_{ij}$. \bluetwo{Another interpretation is that the inflow $\lambda_iP_{ij}\bar a_i$ into the infinite server station equals the outflow $\mu_{ij}\bar f_{ij}$}. Similarly, the constraints in (\hyperref[eqlbm: emptycsv1]{4c}) are also Little's Laws for the number of empty cars travelling from $i$ to $j$; the rate at which empty cars start their journey is $q_{ij}\sum_{k=1}^{r}\mu_{ki}\bar f_{ki}$. The constraints in (\hyperref[eqlbm: emptycsv2]{4d}) say that the total rate of outflow from region $i$, given by $\lambda_i \bar a_i$, must equal the total inflow into the region, given by the right hand side of (\hyperref[eqlbm: emptycsv2]{4d}).} Constraints (\hyperref[eqlbm: boundry]{4e}) state that \blue{a shortage of availability ($1 - \bar a_i > 0$) is only possible when the fraction of cars at the region equals zero ($\bar e_{ii} = 0$), because there are not enough cars to meet demand. Conversely, having a positive mass of cars at a region ($\bar e_{ii} > 0$) implies that all passenger requests are satisfied there ($1 - \bar a_i = 0$).}
 Additional intuition can be gained once the fluid model is introduced and its equilibrium behavior discussed in Section~\ref{sec:fluid}, where the reason why (\hyperref[eqlbm:objective]{4a})--(\hyperref[eqlbm:prob]{4g}) is called the fluid-based optimization problem becomes apparent. Finally, although this optimization problem is stated for static empty-car routing policies, the connection to state-dependent polices will be made in Theorem~\ref{thm:lb}.

The following are our main results. The first establishes the connection between the fluid-based optimization problem and $(\bar E^{(N)}, \bar F^{(N)})$. The second shows that asymptotically, the optimal static policy from the fluid-based optimization outperforms all state dependent policies. The ingredients to prove Theorem~\ref{thm:asymptotic_utility} are developed in Section~\ref{sec:fluid}, and the proof is left to Appendix~\ref{app:main_thm_proof}. Theorem~\ref{thm:lb} is proved in Appendix~\ref{app:lb}.
\begin{theorem}\label{thm:asymptotic_utility}
Let $q,\bar{e}, \bar{f}, \bar{a}$ be a feasible solution to the optimization problem in (\hyperref[eqlbm:objective]{4a})--(\hyperref[eqlbm:prob]{4g}). Set $Q=q.$  \begin{color}{black} Assume $P_{ij} > 0$ for all $1 \leq i,j \leq r$ and  $q_{ii}>0$ for all $1\leq i \leq r.$ Then
\begin{align}
\bar F^{(N)}(\infty) &\Rightarrow \bar f, \label{eq:main1} \\
\bar E_{ij}^{(N)}(\infty) &\Rightarrow \bar e_{ij}, \quad 1 \leq i \neq j \leq r, \\
\bar E_{ii}^{(N)}(\infty) &\Rightarrow 0, \quad \text{ for $i$ such that $\bar a_i < 1$}, \\
\sum_{i : \bar a_i = 1} \bar E_{ii}^{(N)}(\infty) &\Rightarrow \sum_{i : \bar a_{i} = 1} \bar e_{ii}, \label{eq:main4}
\end{align}
and
\begin{align}
\Prob(E_{ii}^{(N)}(\infty) > 0) \to \bar a_i, \quad 1 \leq i \leq r, \label{eq:busy_probs}
\end{align}
as $N \to \infty.$  
\end{color}
\end{theorem}
\begin{remark}
The assumptions that $P_{ij} > 0$ for all $i,j$ and $q_{ii}>0$ for all $i$ are made to facilitate exposition in the proof of Theorem~\ref{thm:stability}, which plays a central role in establishing Theorem~\ref{thm:asymptotic_utility}. We expect that Theorem~\ref{thm:stability}, and hence Theorem~\ref{thm:asymptotic_utility}, holds under the simpler assumption that $(\bar E^{(N)},\bar F^{(N)})$ has a single recurrent class. Assuming $q_{ii}>0$ for all $i$ is not very restrictive, as it simply means that a driver has a positive probability to stay in a region after dropping off a passenger there. \deletetwo{Moreover, in our numerical study in Section~\ref{sec:numerical}, we have observed that all of the optimal $q_{ii}$'s were far from zero.}
\end{remark}

\begin{remark}\label{rem:general_travel} 
The exponentially distributed travel time assumption is non-essential. We know that $\big(E^{(N)}, F^{(N)}\big)$ is a BCMP network, and in full generality, BCMP networks only require that the service time distributions in the infinite server stations have rational Laplace transforms \cite{BaskChanMuntPala1975}. This class of distributions is dense in the set of all probability distributions on $(0,\infty)$ \cite{Asmu2003}. Furthermore, stationary distribution of BCMP networks is known to depend only on the service rates at the stations, and not on the entire distribution. Hence, the results of Theorem~\ref{thm:asymptotic_utility} hold for travel time distributions with rational Laplace transforms.
\end{remark}

\begin{color}{black}
\begin{theorem} \label{thm:lb}
(a) Suppose $(\bar E^{(N)},\bar F^{(N)})$ has a single recurrent class under $P$ and $Q$, where $Q$ can be a state-dependent empty-car routing policy. Let $(q^*,\bar{e}^*, \bar{f}^*, \bar{a}^*)$ be an optimal solution of the optimization problem in (\hyperref[eqlbm:objective]{4a})--(\hyperref[eqlbm:prob]{4g}). Then
\blue{
\begin{align*}
\sum_{i=1}^{r} \sum_{j=1}^{r}  A^{(N)}_i(\infty) \lambda_i P_{ij} c_{ij} < \sum_{i=1}^{r} \sum_{j=1}^{r} \bar a_i^{*}\lambda_i P_{ij} c_{ij} , \quad N > 0.
\end{align*}
}
(b) Let $(\bar E^{(N)*}, \bar F^{(N)*})$ denote the CTMC under the static routing policy $q^*.$ If $P_{ij} > 0$ for all $1 \leq i,j \leq r$ and  $q^*_{ii}>0$ for all $1\leq i \leq r,$ then
\blue{
\begin{align*}
\lim_{N\rightarrow\infty}\sum_{i=1}^{r} \sum_{j=1}^{r}  A^{(N)*}_i(\infty) \lambda_i P_{ij} c_{ij} = \sum_{i=1}^{r} \sum_{j=1}^{r} \bar a_i^{*}\lambda_i P_{ij} c_{ij}.
\end{align*}
}
\end{theorem}

\begin{remark} Part (a) of Theorem \ref{thm:lb} states that the optimal value of the fluid-based optimization problem (\hyperref[eqlbm:objective]{4a})--(\hyperref[eqlbm:prob]{4g}) is a \blue{strict} upper bound on the expected system utility of the system with $N$ cars under any state-dependent routing policy under which the CTMC has a single recurrent class. Part (b) states that the upper bound is asymptotically achievable  under the static routing policy $q^*$ if $P_{ij} > 0$ for all $1 \leq i,j \leq r$ and  $q^*_{ii}>0$ for all $1\leq i \leq r.$
\end{remark}

\end{color}

\subsection{Efficient Solution of the Fluid-Based Optimization}
\label{sec:efficient} 
Having established the relevance of the fluid-based optimization, we now discuss how to solve it efficiently.
\blue{The main issue with solving (\hyperref[eqlbm:objective]{4a})--(\hyperref[eqlbm:prob]{4g}) is the presence of bilinear constraints, e.g.\ (\hyperref[eqlbm: emptycsv1]{4c}) and (\hyperref[eqlbm: emptycsv2]{4d}) are bilinear in $q$ and $\bar f$, and in (\hyperref[eqlbm: boundry]{4e}), the bilinearity is in $\bar a$ and $\bar e$.} In this section \blue{we show that the problem can be transformed into one with only} linear constraints. Our first step is the following lemma, \blue{which transforms the constraints in (\hyperref[eqlbm: emptycsv1]{4c}) and (\hyperref[eqlbm: emptycsv2]{4d})}. It is proved in Appendix~\ref{app:relax0}.

\begin{lemma}\label{lem:relax0}
Consider the set of constraints
\begin{subequations}
\begin{align}
&\lambda_iP_{ij}\bar a_i = \mu_{ij}\bar f_{ij}, \quad 1\leq i,j \leq r, \label{rrp: fullcsv}\\
&\mu_{ij}\bar e_{ij} \leq  \sum_{k=1}^{r}\mu_{ki}\bar f_{ki}, \quad 1\leq i,j\leq r, ~j\neq i, \label{rrp:prob1}\\
&\sum_{k=1,k\neq i}^{r}\mu_{ki}\bar e_{ki} \leq \lambda_i\bar a_i \leq \sum_{k=1,k\neq i}^{r}\mu_{ki}\bar e_{ki} + \sum_{k=1}^{r}\mu_{ki}\bar f_{ki}, \quad 1\leq i \leq r, \label{rrp:prob2}\\
&\lambda_i\bar a_i + \sum_{j=1,j\neq i}^r\mu_{ij}\bar e_{ij} = \sum_{k=1,k\neq i}^r\mu_{ki}\bar e_{ki} + \sum_{k=1}^{r}\mu_{ki}\bar f_{ki}, \quad 1\leq i \leq r, \label{rrp: carcsv}\\
&(\bar e, \bar f) \in \mathcal{T}, \quad \label{rrp: carnum} \\
& 0 \leq \bar a_{i} \leq 1, \quad 1\leq i \leq r, \label{rrp: frac}\\
&(1-\bar a_i)\bar e_{ii} = 0, \quad 1\leq i \leq r. \label{rrp: boundry}
\end{align}
\end{subequations}
\blue{
If $(\bar{e}, \bar{f}, \bar{a})$ and $q$ satisfy (\hyperref[eqlbm: fullcsv]{4b})--(\hyperref[eqlbm:prob]{4g}), then $(\bar{e}, \bar{f}, \bar{a})$ satisfy \eqref{rrp: fullcsv}--\eqref{rrp: boundry}. Conversely, suppose $(\bar{e}, \bar{f}, \bar{a})$ satisfy \eqref{rrp: fullcsv}--\eqref{rrp: boundry} and let 
\begin{align}
q_{ij}=&\ \frac{\mu_{ij}\bar e_{ij}}{\sum_{k=1}^{r}\mu_{ki}\bar f_{ki}}, \quad 1 \leq i\neq j \leq r, \quad q_{ii}= \frac{\lambda_i\bar a_i -\sum_{k=1,k\neq i}^{r}\mu_{ki}\bar e_{ki}}{\sum_{k=1}^{r}\mu_{ki}\bar f_{ki}} \quad 1 \leq i \leq r.\label{rp:intra}
 \end{align}
 Then $(\bar{e}, \bar{f}, \bar{a})$ and $q$ satisfy (\hyperref[eqlbm: fullcsv]{4b})--(\hyperref[eqlbm:prob]{4g}).}
\end{lemma}

With the help of this lemma, the fluid-based optimization problem can be rewritten as
\begin{align}
&\max_{\bar{e}, \bar{f}, \bar{a}} \sum_{i=1}^{r} \sum_{j=1}^{r} \bar a_i\lambda_i P_{ij} c_{ij}  \label{opt-o:obj}\\
\hbox{subject to: } & (\ref{rrp: fullcsv})-(\ref{rrp: boundry}).\label{opt-o:cdts}
\end{align}
Observe that \eqref{rrp: fullcsv}--\eqref{rrp: frac} are all linear constraints, and that only (\ref{rrp: boundry}) is \blue{bilinear}. The following result says that we can safely ignore \eqref{rrp: boundry}. It is proved in Appendix~\ref{app:relax}.
\begin{lemma}
Consider the relaxed optimization problem
\begin{align}
&\max_{\bar{e}, \bar{f}, \bar{a}} \sum_{i=1}^{r} \sum_{j=1}^{r} \bar a_i\lambda_i P_{ij} c_{ij} \label{opt-r:obj}\\
\hbox{subject to: } & (\ref{rrp: fullcsv})-(\ref{rrp: frac}), \label{opt-r:cdts}
\end{align}
and let $(\bar{e}^*, \bar{f}^*, \bar{a}^*)$ be an optimal solution. 
\blue{
\begin{enumerate}
\item If $\bar{a}^{*}_{i} < 1$ for all $1 \leq i \leq r$, then $\bar e^*_{ii} = 0$ for all $ 1 \leq i \leq r$, implying $(\bar{e}^*, \bar{f}^*, \bar{a}^*)$ also satisfies \eqref{rrp: boundry}.
\item \bluetwo{Otherwise, choose any} $i'$ such that $\bar{a}^*_{i'}=1$, define $\tilde e$ by letting $\tilde e_{ij} = \bar e_{ij}^*$ for all $1 \leq i \neq j \leq r$, 
\begin{align*}
 \quad \tilde{e}_{i' i'}= \sum_{i=1}^r \bar e_{ii}^*, \quad \text{ and } \quad  \tilde e_{ii} =  0, \quad i \neq i'.
\end{align*}
Then $(\tilde{e},\bar{f}^{*}, \bar{a}^{*})$ is a feasible solution that yields the same objective value as $(\bar e^*, \bar f^*, \bar a^*)$, and also satisfies \eqref{rrp: boundry}.
\end{enumerate}
}
\label{lem:relax}
\end{lemma}
\blue{
Solving the fluid-based optimization can be broken down into the following procedure. We first solve \eqref{opt-r:obj}--\eqref{opt-r:cdts}, which has ($4r^2+2r+1$) linear constraints and $2r^2 + r$ variables, and can be solved efficiently using any standard linear program solver. We then modify the solution using Lemma~\ref{lem:relax} (if needed) so that it also satisfies \eqref{rrp: boundry}, and recover the optimal routing policy using  \eqref{rp:intra}.
}

\begin{remark}
\label{rem:general_utility}
 So far, we only considered the utility function defined in \eqref{eqlbm:objective}. In fact, both Lemma~\ref{lem:relax} and Theorem~\ref{thm:lb} can be extended to hold for any function $U(\bar e, \bar f, \bar a)$ that is
\begin{enumerate}[(a)]
\item  nondecreasing in $\bar a_i$ for all $i,$

\item nondecreasing in $\bar f_{ij}$ for all $i$ and $j,$

\item nonincreasing in $\bar e_{ij}$ for all $i\not=j,$
\item  independent of $\bar e_{ii}$ for all $i,$ and

\item concave in $(\bar e, \bar f)$.
\label{assumd}
\end{enumerate}
Lemma~\ref{lem:relax} can be extended because its proof relies only on the objective value in \eqref{opt-r:obj} satisfying conditions (a)--(d). Theorem~\ref{thm:lb} can be extended to say that
\begin{align*}
\E \big[ U(A^{(N)}(\infty), E^{(N)}(\infty), F^{(N)}(\infty)\big] \leq U\big(A^{(N)}(\infty), \E \big[E^{(N)}(\infty)\big], \E \big[F^{(N)}(\infty)\big]\big)   \leq U(\bar a^{*}, \bar e^{*}, \bar f^{*}),
\end{align*}
where the first inequality follows from Jensen's inequality and (e), and the second inequality is proved by repeating the proof of Theorem~\ref{thm:lb}.
\end{remark}

\section{\blue{A Numerical Study}} \label{sec:numerical}
This section is devoted to a numerical study of our empty-car routing policy. To ground the study, we use a data set obtained from a data challenge by the Didi Research Institute \citep{Didi}. Using the data set, we extract a nine-region network and the associated realistic parameters $N, \mu, \lambda$, and $P$. \bluetwo{For more details about the dataset and how we obtained our parameters, see Appendix~\ref{app:nineregion}.}

Figure \ref{fig: traffic} shows order fulfillment levels for each of the nine regions during the 5PM-6PM evening rush hour, plotted over $21$ days. Each point represents the total number of orders received, and orders fulfilled during that one hour window. We can see that three of the nine regions, regions 13, 47, and 50, had significant supply shortages in most of the 21 days, while most orders in the remaining six regions were fulfilled. Our data did not permit us to deduce the surplus of drivers in those six regions, but these figures illustrate the significance of a good empty-car routing policy.
\begin{figure}[h]
\begin{minipage}{0.3\textwidth}\centering
  \includegraphics[width=1\textwidth]{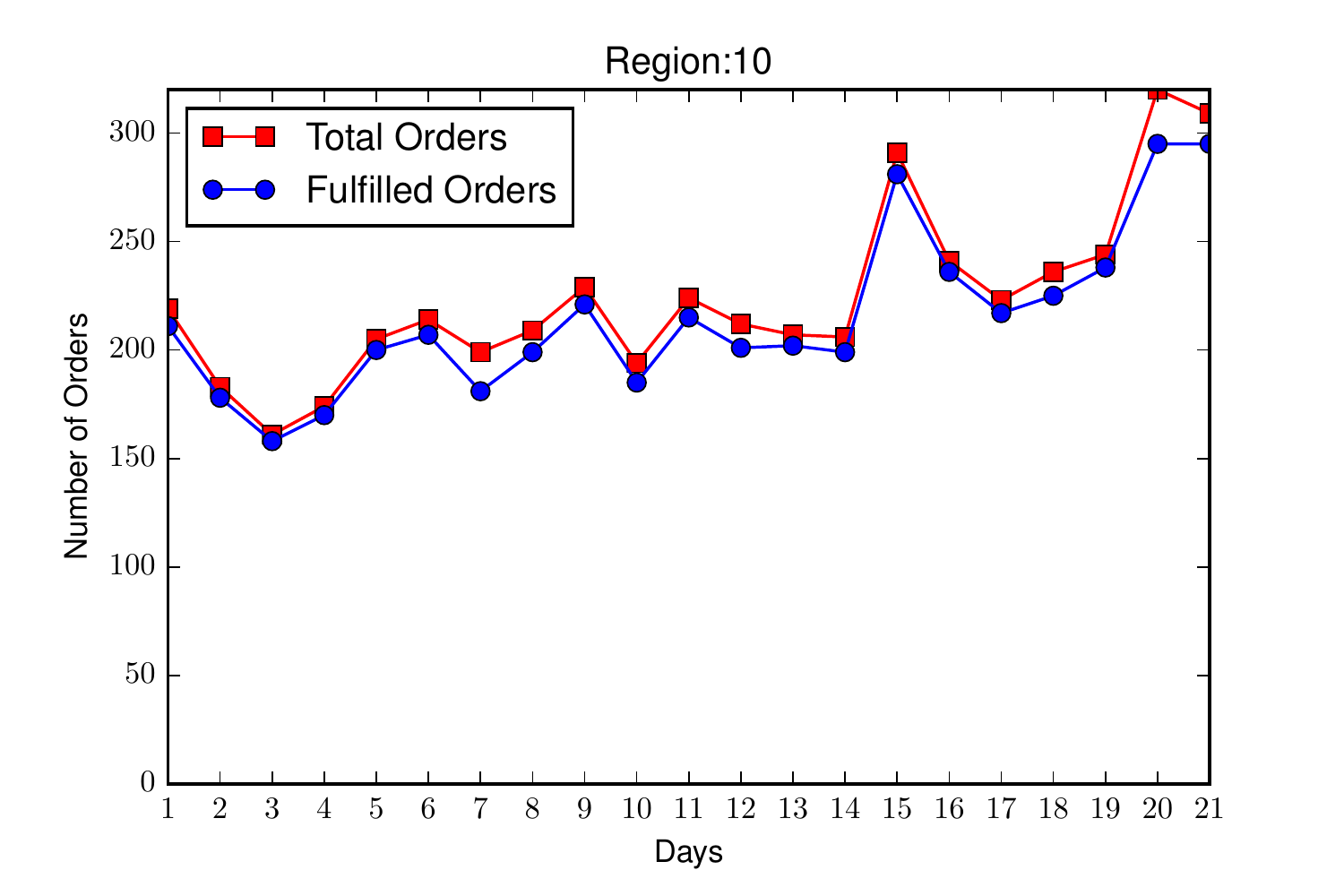}
  \end{minipage}
\begin{minipage}{0.3\textwidth}\centering
    \includegraphics[width=1\textwidth]{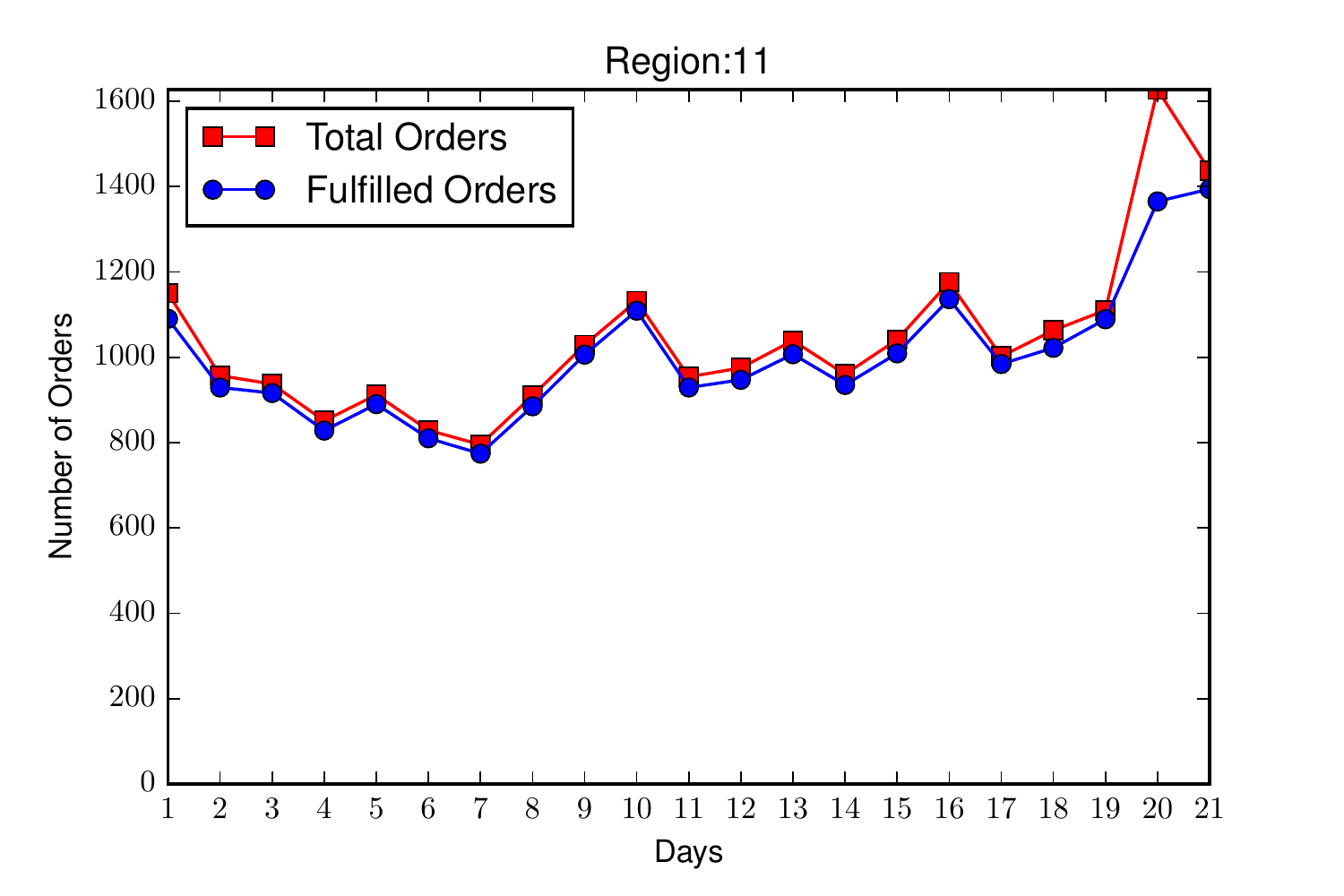}
   \end{minipage}
   \begin{minipage}{0.3\textwidth}\centering
    \includegraphics[width=1\textwidth]{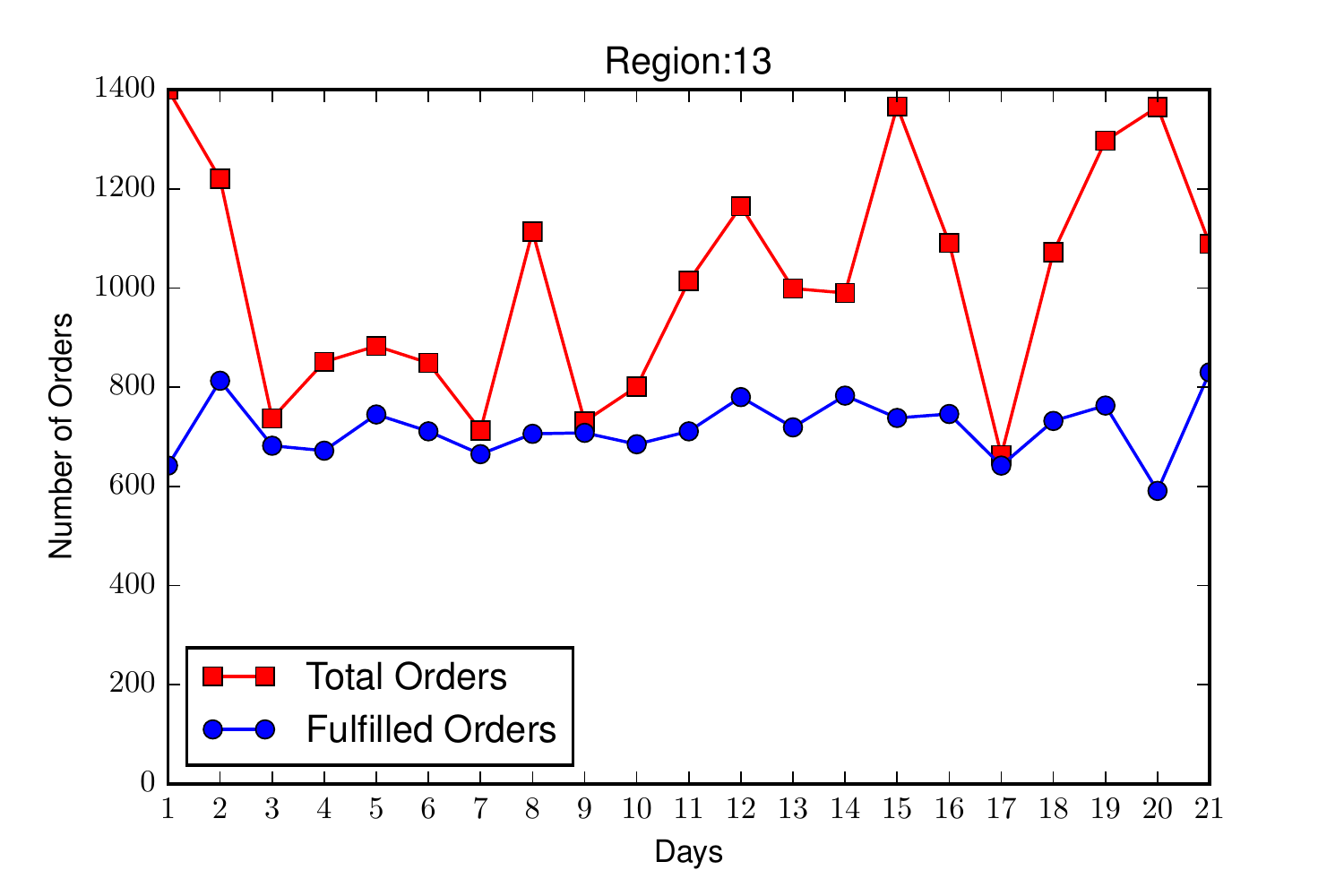}
   \end{minipage}
  \\
\begin{minipage}{0.3\textwidth}\centering
  \includegraphics[width=1\textwidth]{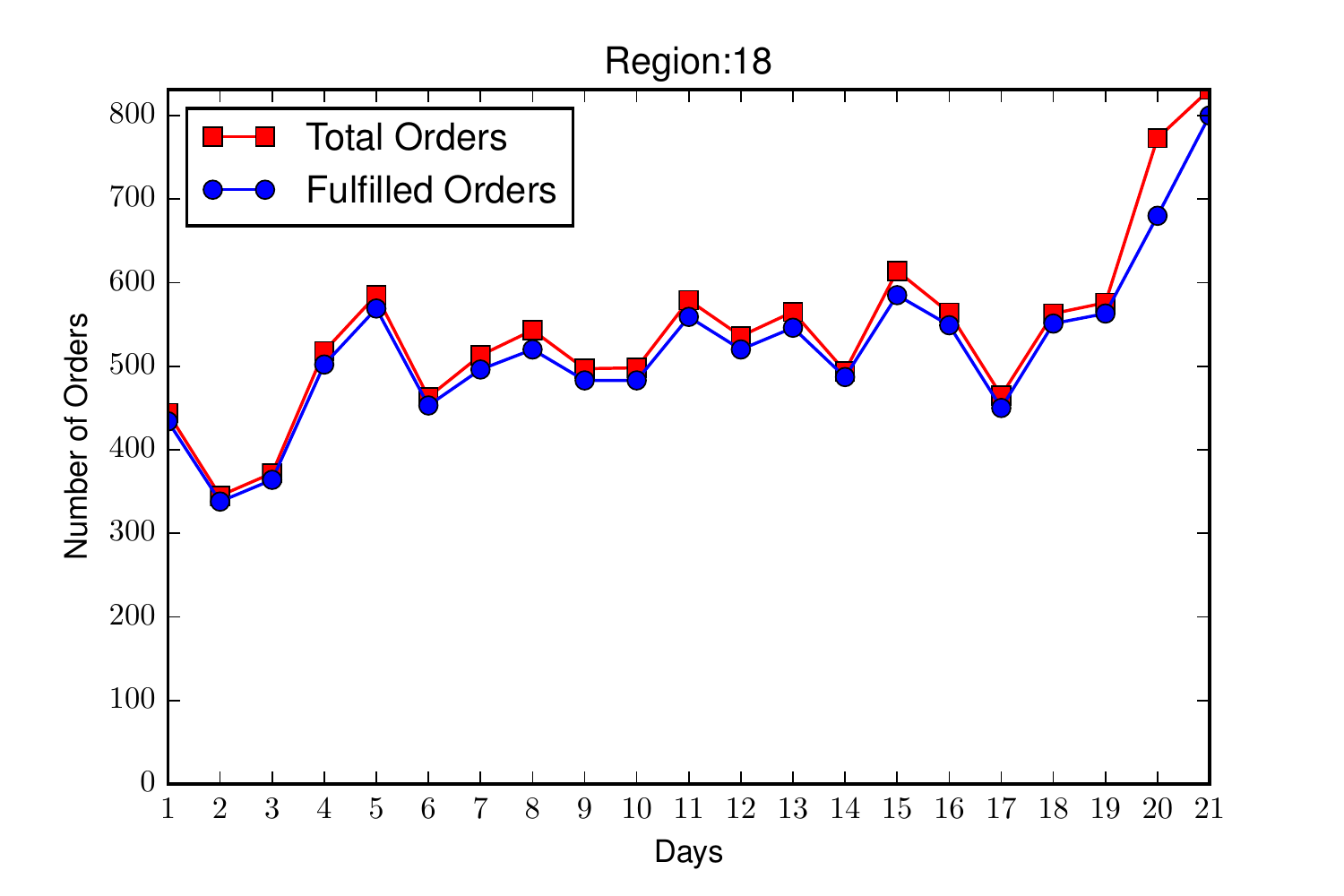}
  \end{minipage}
\begin{minipage}{0.3\textwidth}\centering
    \includegraphics[width=1\textwidth]{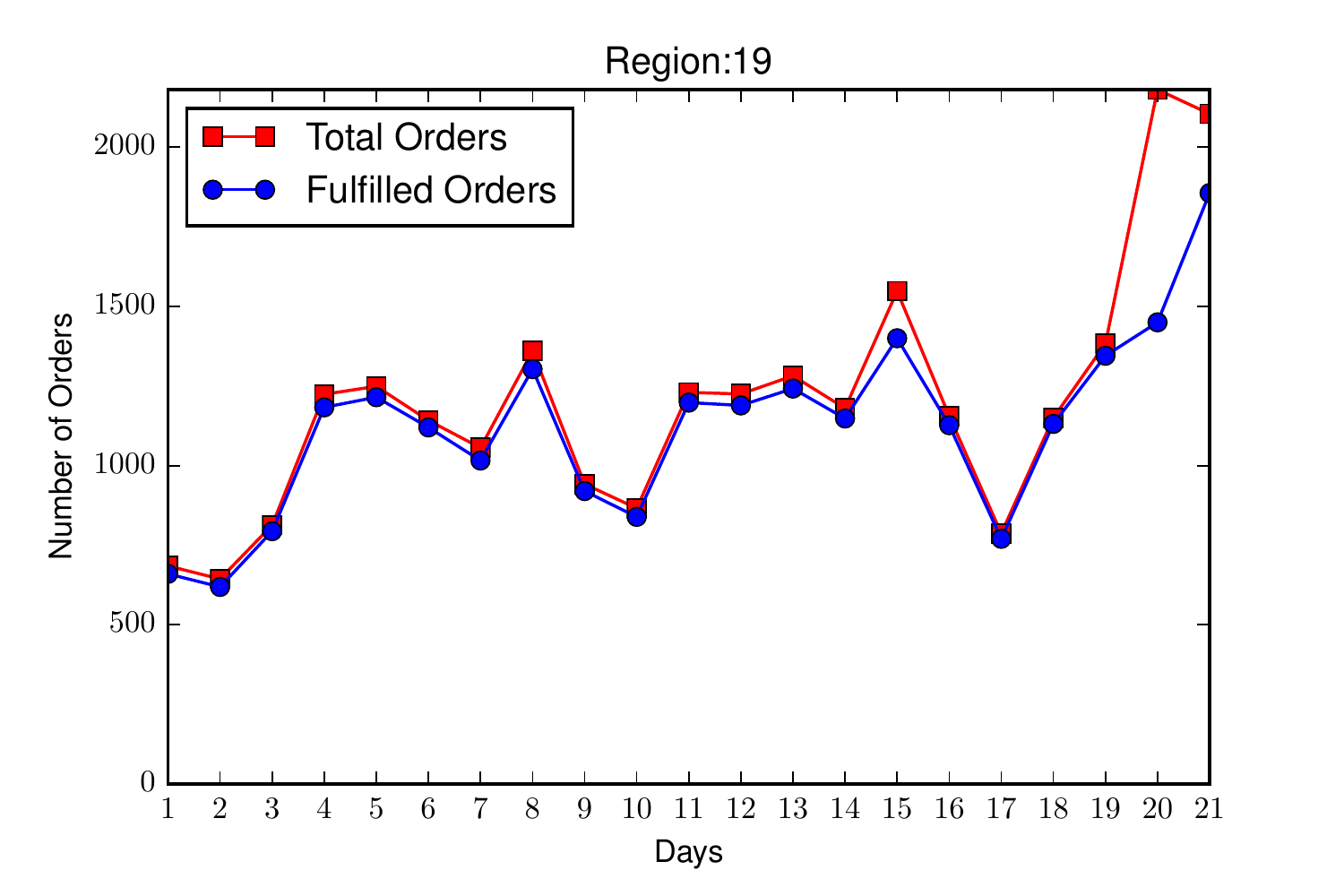}
   \end{minipage}
   \begin{minipage}{0.3\textwidth}\centering
    \includegraphics[width=1\textwidth]{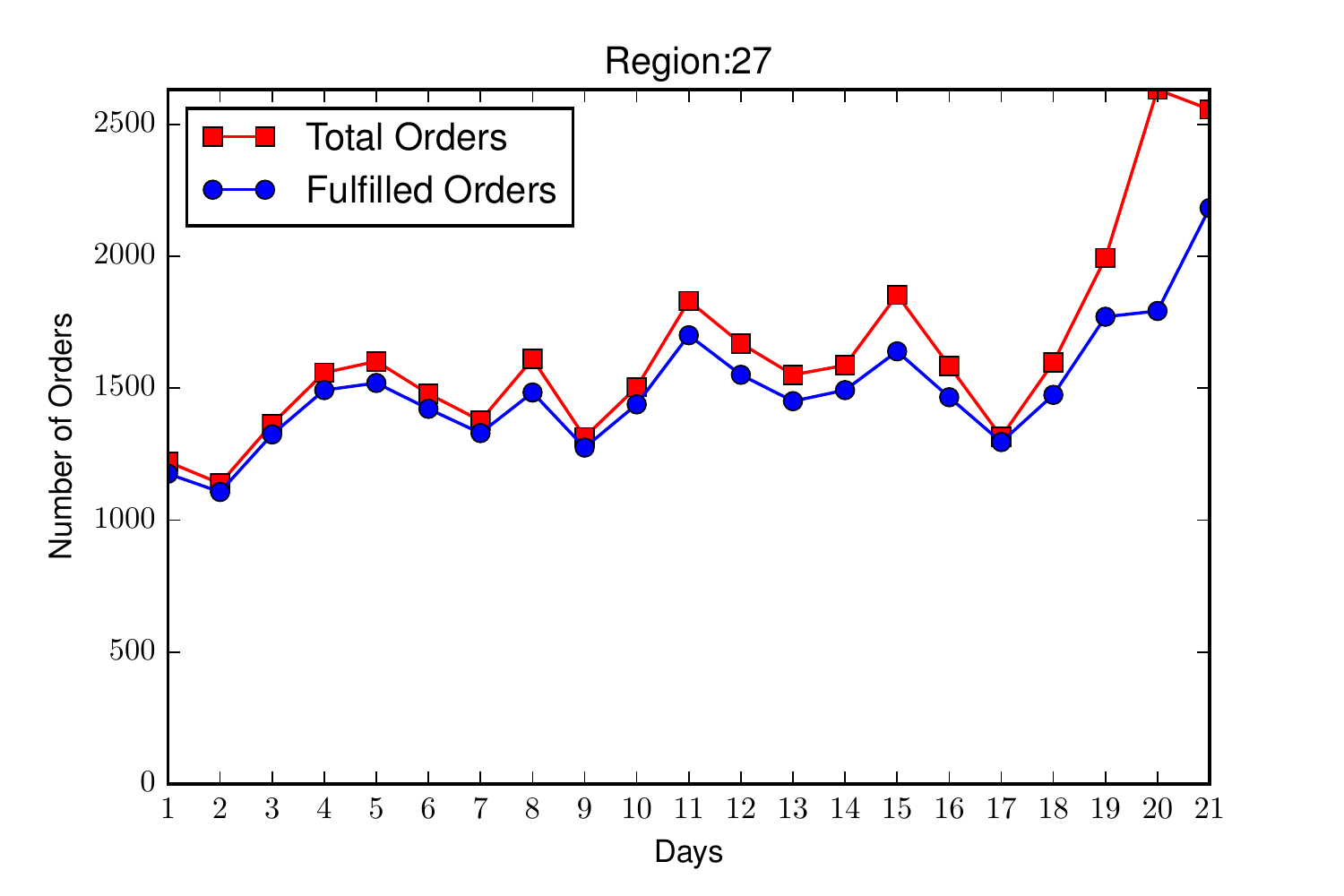}
   \end{minipage}
 \\
   \begin{minipage}{0.3\textwidth}\centering
  \includegraphics[width=1\textwidth]{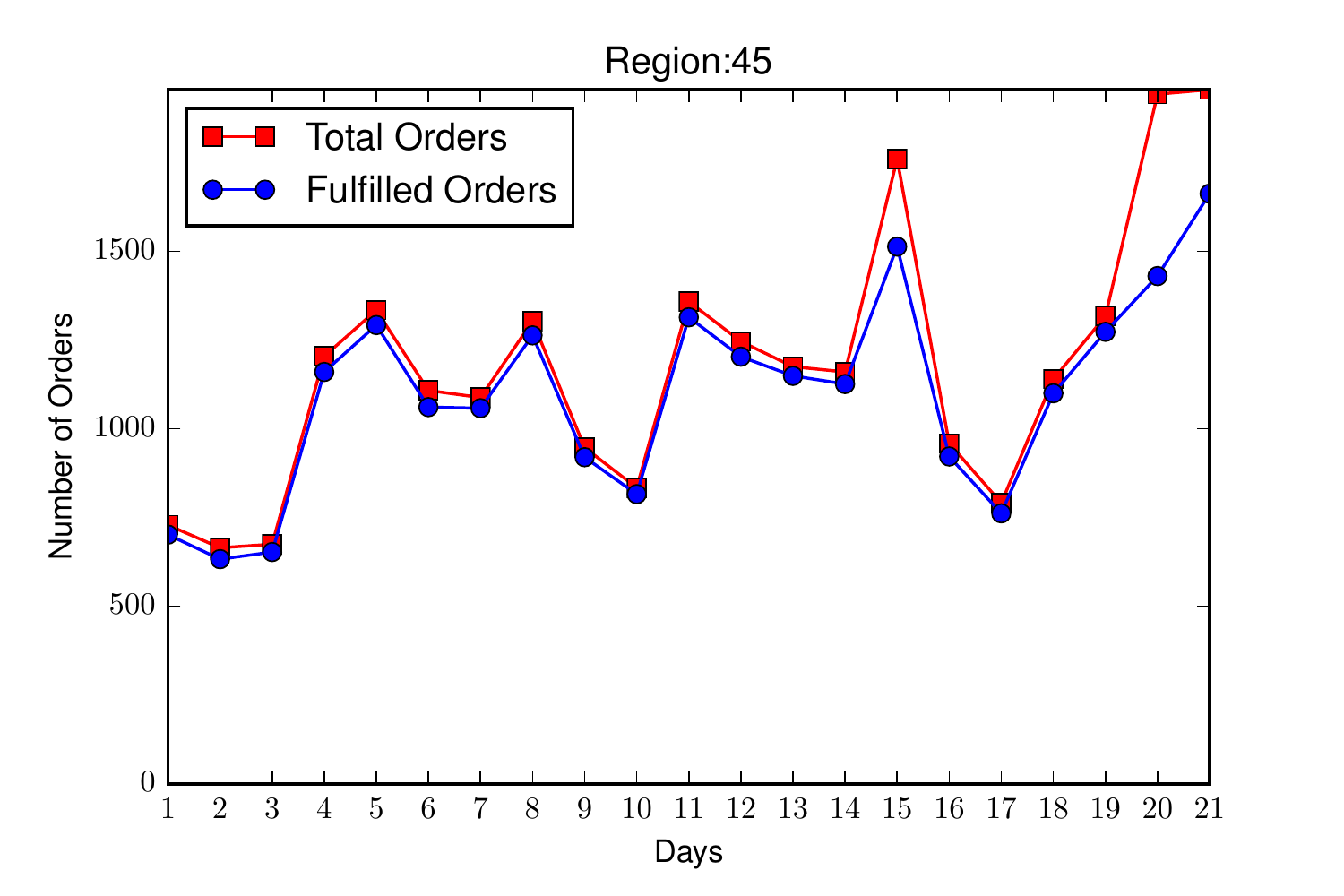}
  \end{minipage}
\begin{minipage}{0.3\textwidth}\centering
    \includegraphics[width=1\textwidth]{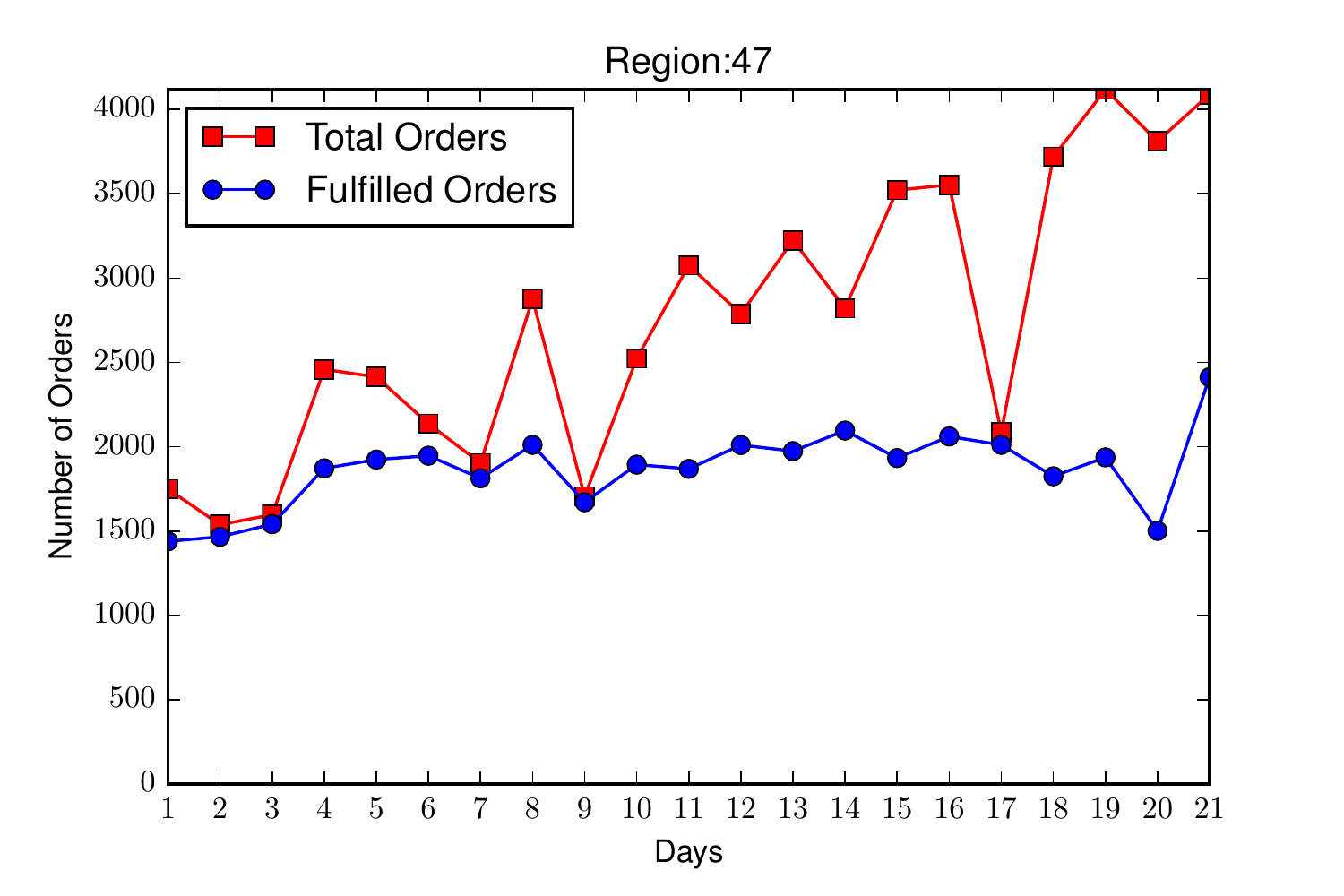}
   \end{minipage}
   \begin{minipage}{0.3\textwidth}\centering
    \includegraphics[width=1\textwidth]{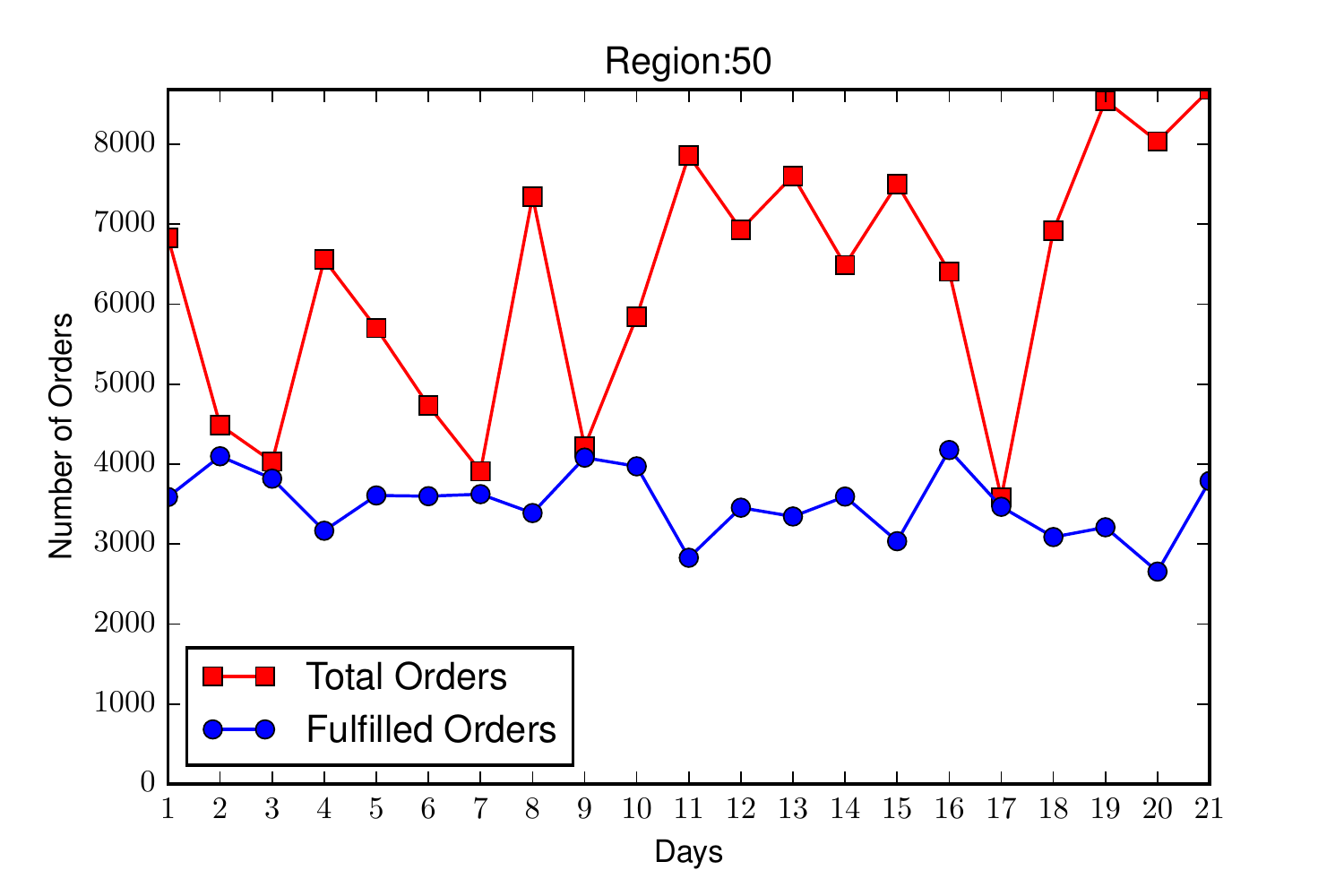}
   \end{minipage}
   \caption{The gaps between the number of passenger orders and the number of fulfilled orders of the nine regions}
   \label{fig: traffic}
\end{figure}

The rest of the section is structured as follows. We first compare the fluid-based optimal routing policy to a state-dependent routing policy in Section~\ref{sec:dynamic}. We then use the fluid-based optimization to propose a lookahead policy to deal with cases where the system parameters are not constant over time in Section~\ref{sec:lookahead}. Lastly, we perform some numerical robustness tests in Section~\ref{sec:robust} to see how the fluid-based policy performs when there is estimation noise present in the parameters.

\subsection{\blue{Performance Comparison with Dynamic Routing Polices}}
\label{sec:dynamic}

In Theorem \ref{thm:lb}, we showed that the expected utility under any state-dependent routing policy for the finite-sized system is upper bounded by the optimal utility of the fluid-based optimization. However, the question remains open whether a state-dependent routing policy can outperform the static routing policy $Q^*$ for a finite-sized system. \bluetwo{Choosing $c_{ij} = 1/\sum_{i=k}^{r} \lambda_k$ in (\hyperref[eqlbm:objective]{4a})}, we consider the utility function
\begin{align}
U(\bar e, \bar f, \bar a) = U(\bar a) = \frac{\sum_{i=1}^{r} \bar a_{i}\lambda_i}{\sum_{i=1}^{r}\lambda_i}. \label{eq:obj_data}
\end{align}
This can be thought of as the probability that a passenger requesting a ride at \emph{any} region is fulfilled. It is impossible to consider all possible dynamic routing policies, so we focus on the following two intuitive heuristics.

\noindent{\bf Join-the-Least-Congested-Region with Threshold $\eta$ (JLCR-$\eta$):} When a car drops off a passenger at region $i$ at time $t,$ the driver stays at region $i$ if
\begin{align}
(1-\eta)\frac{\sum_{k} E^{(N)}_{ki}(t)}{\lambda_i} \leq \min_{j=1,j\not=i} \frac{\sum_{k} E^{(N)}_{kj}(t)}{\lambda_j}.
\end{align}
Otherwise, the driver drives empty to region $j^*$, where
\begin{align}
j^* \in \arg \min_{j=1,j\not=i} \frac{\sum_{k} E^{(N)}_{kj}(t)}{\lambda_j}.
\end{align} Ties are broken uniformly at random. \hfill{$\square$}

To understand the JLCR-$\eta$ policy, we note that  $\sum_{k} E^{(N)}_{ki}(t) = E^{(N)}_{ii}(t) + \sum_{k \neq j} E^{(N)}_{ki}(t)$ is the number of empty cars both currently waiting and en-route to region $i$. Therefore, $\frac{\sum_{k} E^{(N)}_{ki}(t)}{\lambda_i}$ is a measure of congestion, in terms of empty cars, at region $i$. When $\eta=0$, the policy routes empty cars to the least congested region. However, such a policy can be wasteful if congestion levels among regions are similar, because it takes time for a car to go from one region to another. We therefore introduce the threshold $\eta$ such that a driver drives empty from $i$ to $j$ only if the difference in congestion levels surpasses $\eta \frac{\sum_{k} E^{(N)}_{ki}(t)}{\lambda_i}$. We test the policy on our $9$-region network with parameters as in \eqref{eq:pdidi}-\eqref{eq:lambdadidi}. \bluetwo{With $\eta$ ranging from $0$ to $1$, we find that JLCR-$\eta$ performs best when $\eta$ is around $ 0.5$.}
In addition to JLCR, we also consider the following policy where drivers aim to minimize the time until their next pickup.

\noindent{\bf Shortest-Wait (SW):} When a car drops off a passenger at region $i$ at time $t,$ the driver stays at region $i$ if
\begin{align*}
\frac{E_{ii}^{(N)}(t)}{N\lambda_i} \leq \min_{j \neq i} \frac{1}{\mu_{ij}} +   \frac{\Big(E_{jj}^{(N)}(t) + \frac{1}{\mu_{ij}} \sum_{k \neq j} \mu_{kj} E_{kj}^{(N)}(t) - \frac{ N\lambda_j}{\mu_{ij}} \Big)^+ }{ N\lambda_j},
\end{align*}
Otherwise, the driver drives empty to region  that minimizes the right hand side above. Ties are broken uniformly at random. \hfill{$\square$}

The intuition behind the SW policy is that drivers want to minimize the time to get their next passenger. Since passengers arrive to region $i$ every $1/N \lambda_i$ time units (on average), then $E_{ii}^{(N)}(t)/N\lambda_i$ is a proxy (we have not assumed any priority scheme for choosing how to allocate passengers between multiple cars in the same region) for the amount of time it will take until the driver gets a passenger if he stays in region $i$. If the driver chooses to go to region $j$, then the time taken until he gets a passenger is the sum of the travel time $1/\mu_{ij}$ and the time spent idling in region $j$. We use $\Big(E_{jj}^{(N)}(t) + \frac{1}{\mu_{ij}} \sum_{k \neq j} \mu_{kj} E_{kj}^{(N)}(t) - \frac{ N\lambda_j}{\mu_{ij}} \Big)^+ / N\lambda_j$ as a proxy for the latter quantity: $E_{jj}^{(N)}(t)$ is the number of cars idling in region $j$ at the decision point, $\frac{1}{\mu_{ij}} \sum_{k \neq j} \mu_{kj} E_{kj}^{(N)}(t)$ estimates the number of cars to arrive to $j$, and $\frac{ N\lambda_j}{\mu_{ij}}$ estimates the number of cars that will leave the region due to passenger arrivals by the time the driver makes it from $i$ to $j$.

\begin{figure}[ht]
  \centering
  \includegraphics[width=4in]{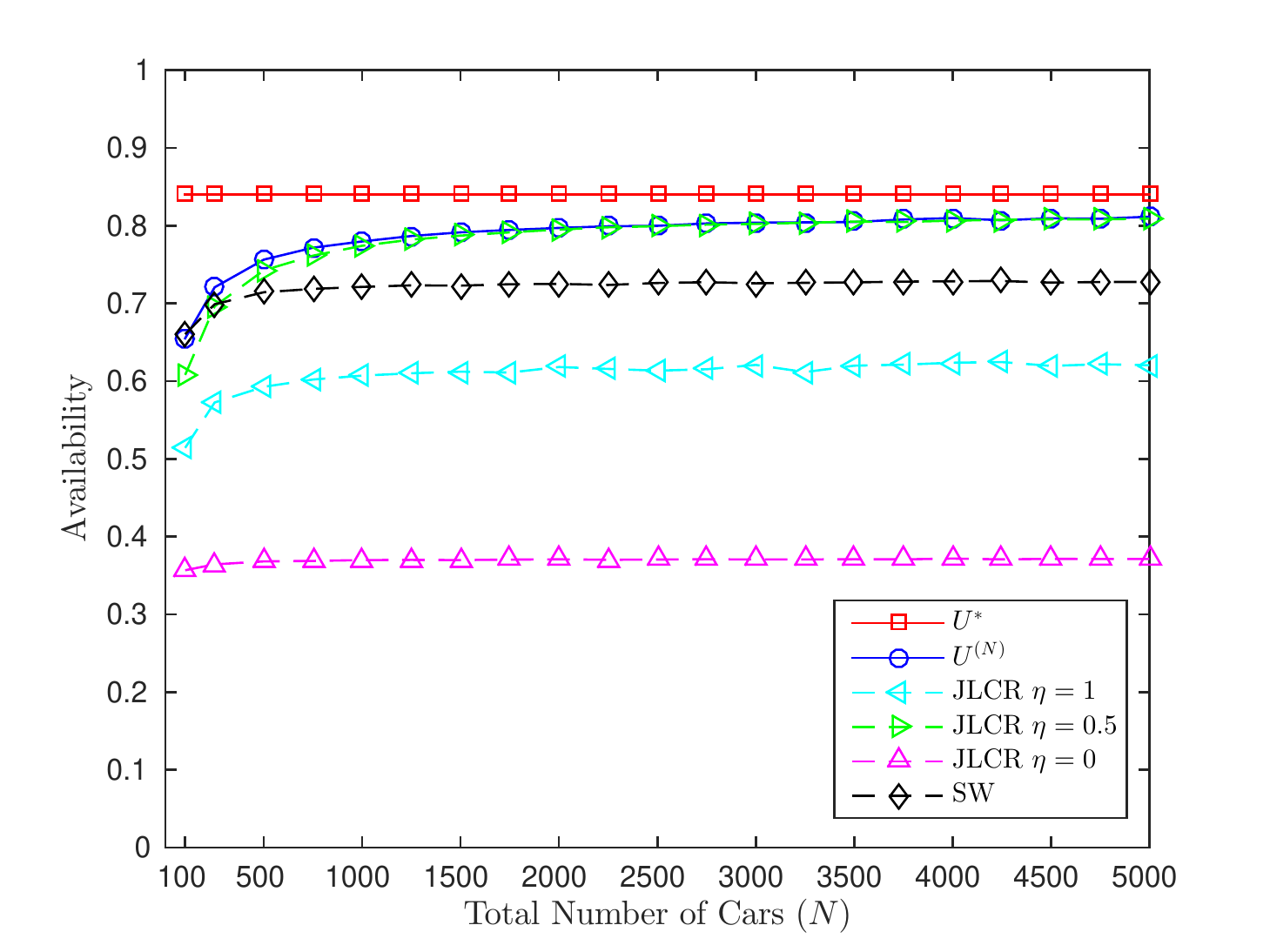}
  \caption{Performance comparison between the static routing, JLCR and SW policies in the nine-region network. In the plot, $U^*$ is the optimal utility \eqref{eq:obj_data} from the fluid-based optimization, and $U^{(N)}$ is the steady-state utility in the finite sized system under optimal static routing.}
\label{fig:comp-didi}
\end{figure}

Figure \ref{fig:comp-didi} compares the static routing policy under $Q^*$ to $SW$ and JLCR-$\eta$ with different values of $\eta$ in the 9-region network. In particular, we included JLCR policies with
\begin{itemize}
  \item $\eta=0:$ Under JLCR-0, an empty car always goes to the least congested region.

  \item $\eta=1:$ Under JLCR-1, after a car drops off a passenger, it always stays at its current region.

  \item $\eta=0.5:$ JLCR-0.5 maximizes system-wide availability among all JLCR-$\eta$ when $N=2,000$.
\end{itemize}
A few remarks are in order. The figure confirms that static routing with $Q^*$ outperforms both the $SW$ and JLCR-$\eta$ family of policies. However, a typical quality of state-dependent policies is robustness to system parameters. In our case, computing $Q^*$ requires knowledge of $\lambda, \mu$, and $P$, whereas a JLCR-$0.5$ only requires knowledge of $\lambda$. Robustness to parameters is a particularly important quality when one only has noisy observations of the true parameters, or when the true parameters change over time. The objective of this paper is not to pursue optimal state-dependent policies. Rather, it is to establish an initial, rigorous theoretical foundation for the study of ridesharing networks. Our optimal static policy can then be used as a benchmark against which one compares the performance of other routing policies. As a case in point, had Figure~\ref{fig:comp-didi} not included the performance of static-routing under $Q^*$, it would have been impossible to say whether JLCR-$\eta^*$ was a good policy or not.

\subsection{\blue{Time-varying Parameters and the Lookahead Heuristic}}
\label{sec:lookahead}
\begin{figure}[h] 
    \includegraphics[width=0.5\textwidth]{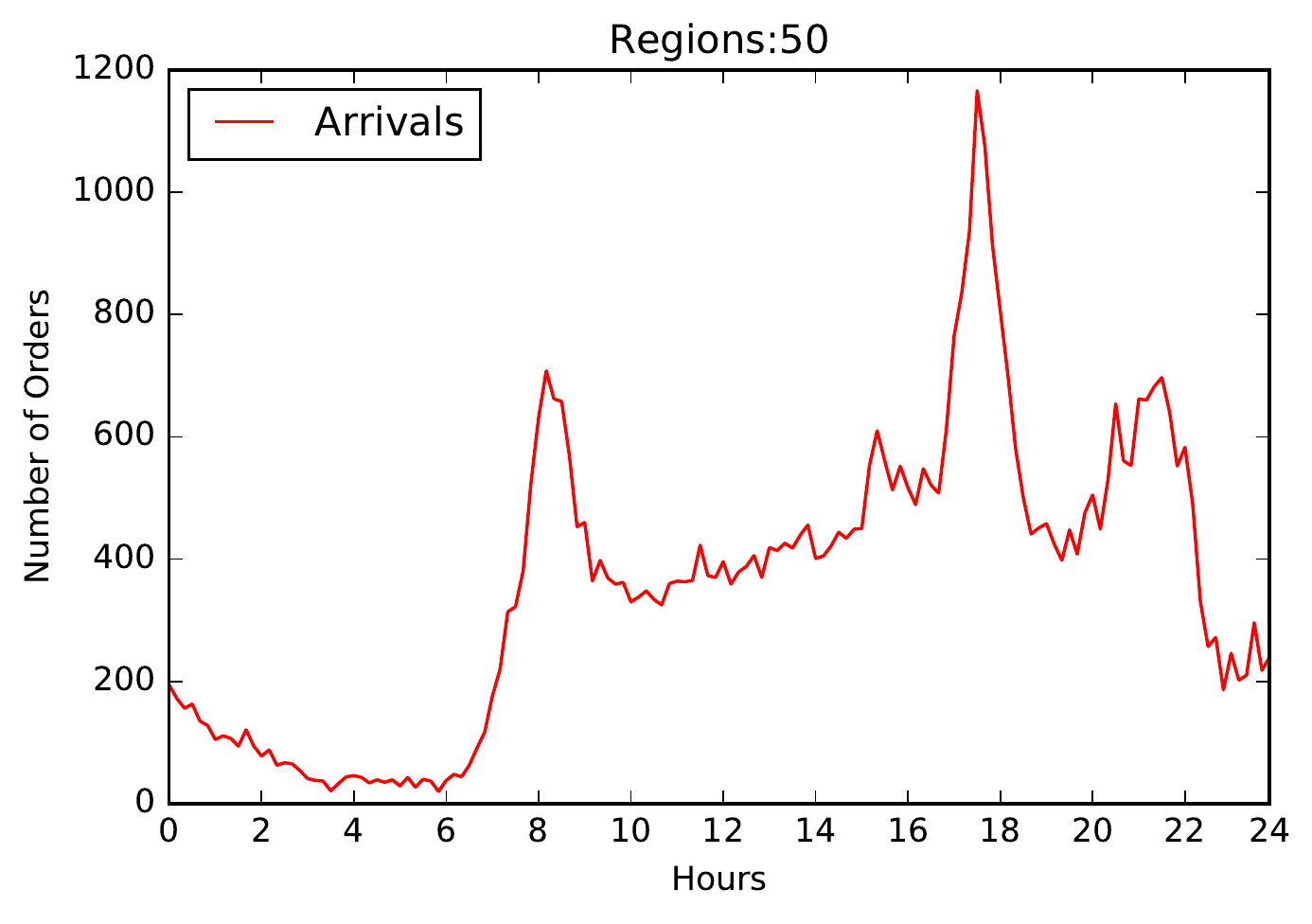}
    \centering
    \caption{Ride requests on January 5, 2016 in region 50 of the Didi dataset. The x-axis shows the time of day, with $0$ being midnight and $12$ being noon.} \label{fig:arrivals}
\end{figure}

The effectiveness of the fluid-based empty-car routing policy rests on the assumptions that a) parameters remain constant over time, and b) that the system has reached equilibrium. In practice a) is violated, e.g.\ when a rush hour starts. \deletetwo{; see Figure~\ref{fig:arrivals} where we use our dataset to present a typical snapshot of passenger ride requests throughout the day. Steady-state analysis can still be useful even in the presence of time-varying parameters.} A common solution is to divide the day into periods of time where parameters are assumed constant, and treat the system as if it were in steady-state during each of those periods \cite{GreeKoleWhit2007}. This approach can yield good results provided our time windows are long enough for the system to converge to equilibrium in each of them. 

The Didi dataset we use suggests it is reasonable to assume constant parameters over time windows 1-2 hours in length; see Figure~\ref{fig:arrivals}.
 However, numerical experiments suggest that for certain choices of parameters and initial conditions, convergence of our system to equilibrium can occur on timescales on the order of 10 hours. With the rate at which parameters vary, and the slow convergence to equilibrium, the system never really reaches steady-state. To address this, we now propose a lookahead policy that is grounded in the fluid-based optimization problem. \deletetwo{We then demonstrate via numerical examples that this lookahead heuristic is useful even in cases when both a) and b) are violated.} The purpose is to demonstrate that our fluid model framework, even though developed in a stationary environment, can potentially be used for other purposes.

\noindent{\bf The $T$-lookahead policy:} Recall the per-ride rewards $c_{ij}$ from (\hyperref[eqlbm:objective]{4a}), and suppose that they can depend on $\lambda, \mu$, and $P$, i.e.\ $c_{ij} = c_{ij}(\lambda, \mu,P)$ (like in \eqref{eq:obj_data}). Given time-varying parameter $\big\{\lambda(t), \mu(t), P(t) \big\}_{t \geq 0}$, the per-ride rewards $c_{ij}$ also depend on time; for simplicity, let us write $c_{ij}(t)$ instead of $c_{ij}(\lambda(t), \mu(t),P(t))$. Instead of using the time-independent routing matrix $q^{*}$ from the fluid-based optimization problem (\hyperref[eqlbm:objective]{4a})--(\hyperref[eqlbm:prob]{4g}), at time $t \geq 0$ we use the routing matrix $q^{*}(t)$ that solves
\begin{align*}
&\max_{q,\bar e, \bar f, \bar a} \frac{1}{T} \int_{0}^{T}\sum_{j=1}^{r} \bar a_i\lambda_i(t+u) P_{ij}(t+u) c_{ij}(t+u) du   \\
\text{ subject to } \quad & \frac{1}{T} \int_{0}^{T} \lambda_i(t+u)P_{ij}(s+u)\bar a_i du = \frac{1}{T} \int_{0}^{T}\mu_{ij}(s+u)\bar f_{ij} du, \quad 1\leq i,j \leq r,\\
&\frac{1}{T} \int_{0}^{T}\mu_{ij}(t+u) du\bar e_{ij} = q_{ij}\sum_{k=1}^{r}\frac{1}{T} \int_{0}^{T}\mu_{ki}(t+u) \bar f_{ki} du, \quad 1\leq i,j \leq r, ~j\neq i,\\
&\frac{1}{T} \int_{0}^{T}\lambda_i(t+u)\bar a_i du \\
=&\  \sum_{k=1,k\neq i}^{r}\frac{1}{T} \int_{0}^{T}\mu_{ki}(t+u)\bar e_{ki} du + q_{ii}\sum_{k=1}^{r}\frac{1}{T} \int_{0}^{T}\mu_{ki}(t+u)\bar f_{ki}du, \quad 1\leq i \leq r,\\
&(1-\bar a_i)\bar e_{ii} = 0, \quad 1\leq i \leq r, \\
&(\bar e, \bar f) \in \mathcal{T}, \quad q_{ij}\geq 0, \quad \sum_{j=1}^{r}q_{ij} = 1, \quad 0 \leq \bar a_{i} \leq 1, \quad 1\leq i,j \leq r.
\end{align*}
The above problem is solved in the exact same way as the fluid-based optimization. In other words, the $T$-lookahead policy is a time-varying routing policy that uses paramater averages over the time window $[t,t+T]$ to make a decision at time $t$. In theory, $q^{*}(t)$ can be computed in real time at any decision point. In our numerical results, we discretize time into $\Delta$-spaced intervals and use $q^{*}(k\Delta)$ as the routing decision for all times $t \in [k \Delta, (k+1)\Delta]$, where $k$ is some non-negative integer. \deletetwo{Note that the above policy can be extended to more general utility functions, just like the fluid-based optimization; cf.\ Remark~\ref{rem:general_utility} in Section~\ref{sec:efficient}.}


We now consider \bluetwo{some} example networks with time-varying parameters, and compare the $T$-lookahead policy to both the standard fluid-based routing policy, and the state-dependent JLCR and SW policies from Section~\ref{sec:dynamic}. \bluetwo{Since ridesharing trip lengths are typically on the order of 10-30 minutes, we choose $T$ to equal $30$ and $45$ minutes in our examples. This corresponds to looking ahead a few trips into the future.}

\subsubsection{The 5-Region City}  \label{sec:fiveregion}
\bluetwo{We consider a simplified model of a city that consists of $5$ regions: a downtown area, a midtown area, and three suburban areas. We take a typical evening, from 5pm-11pm, and divide it into three 2-hour slots each having different parameters $\lambda, \mu, P$ to represent different traffic patterns. The details of this 5-region network are provided in Appendix~\ref{app:fiveregion}. }
 
 Table~\ref{tab:toycity} displays  the results of a simulation that compares \bluetwo{ the system-wide availability introduced in \eqref{eq:obj_data}} under several routing policies: `standard fluid', JLCR-$0.5$, SW, 30 minute lookahead, and 45 minute lookahead. The `standard fluid' policy treats the system as if it were in equilibrium during each of the 2-hour time slots. It solves three separate fluid-based optimization problems, one for each of the three 2-hour time slots, and uses the resulting routing policy in the appropriate time slot. The JLCR-$0.5$ and SW policies are used because they adapt quickly to changes in parameters due to their state-dependent nature.

At 5pm, the system is initialized by assuming all cars are idle and distributing them across regions proportionally with the expected demand in the region, i.e.\ region $i$ gets $\lambda_i/ \sum_{j} \lambda_j$ of the cars, where $\lambda$ are the 5pm arrival rates. The results are consistent with what we would expect, and in our experiments, the performance of the lookahead policy was consistent across different choices of initial system configurations (e.g.\ starting system according to fluid equilibrium). The standard fluid policy performs very well during the first two hours because the system was initialized with a favorable initial condition, and parameters stay constant during that time. However, the performance of this policy degrades as soon as the parameter change happens at 7pm. The explanation for this is simple: operating under the 5-7pm parameters for the first two hours puts the system in a state that is very far from the 7-9pm fluid equilibrium, and performance suffers as a result. Unlike the standard fluid policy, the lookahead policies anticipate parameter changes and prepare for them by pre-positioning the cars in the system into a more favorable state. For example, by sacrificing $4\%$ in performance during 6-7pm, the $45$ minute lookahead policy increases performance during 7-8pm to $86\%$ from the $67\%$ of the standard fluid policy. \bluetwo{Factoring in the different ride request rates during the above two hours, the $4\%$ performance drop corresponds to $0.04 * N \sum_{i}\lambda_i^{6-7\text{pm}} \approx 61$ customers lost, while the $19\%$ increase during 7-8pm means an extra $433$ customers served.}
\begin{table}[h] 
\centering
\begin{tabular}{|c|c|c|c|c|c|c|c|}
\hline
          & 5-6pm & 6-7pm & 7-8pm & 8-9pm & 9-10pm & 10-11pm & Total 5-11pm \\ \hline
Standard Fluid   &0.93 &0.82 &0.67 &0.66 &0.76 &0.72 &0.75 \\ \hline
JLCR-$0.5$   &0.91 &0.73 &0.74 &0.75 &0.97 &0.66 &0.79 \\ \hline
SW   &0.92 &0.77 &0.74 &0.75 &0.94 &0.78 &0.81\\ \hline
$T$-Lookahead, $T=0.5$  &0.93 &0.82 &0.82 &0.8 &0.96 &0.71 &0.84\\ \hline
$T$-Lookahead,  $T=0.75$  &0.93 &0.78 &0.86 &0.76 &0.93 &0.76 &0.83\\ \hline
\end{tabular}
\caption{Comparing the performance of the lookahead policy ($N=1000$). The time unit for $T$ is hours and the time-discretization $\Delta$ is chosen to be one minute. Based on the utility function in \eqref{eq:obj_data}, the values displayed are the fraction of fulfilled ride requests.  We used constant travel times, but the results remain qualitatively unchanged if the travel times are exponentially distributed.  \label{tab:toycity}}
\end{table}
\subsubsection{The 9-Region Didi Network}
\label{sec:dididata}
Here we present some numerical results based on the $9$-region network detailed in Appendix~\ref{app:nineregion}. In our first experiment, we consider a 4-hour period where passenger arrivals are constant for the first two hours, and then change abruptly in the last two hours. Namely, for the first two hours we let our arrival rates equal $0.3 \lambda$, where $\lambda$ is as in \eqref{eq:lambdadidi}, and for the last two hours we use $0.85 \lambda^p$, where 
\begin{align*}
\lambda^{p} 
=
\begin{pmatrix}
    0.0131 & 0.0624 & 0.1178 & 0.0870 & 0.0652 & 0.0381 & 0.0762 &  0.2751 & 0.1438
\end{pmatrix}
\end{align*}
is a permutation of $\lambda$. With this experiment, we want to create a case where a) the passenger request rate in the last two hours is significantly higher than the first two, and b) where the sources of passenger requests change significantly. Similar to the 5-region city example, we initialize cars in our system in proportion to $\lambda_i/ \sum_{j} \lambda_j$, and we present our results in Table~\ref{tab:CounterTraffic1000}. The results again show that the performance of the standard fluid policy suffers right after the change in parameters and that the lookahead policy corrects for this. In this example, using the longer $45$ minute lookahead window does not degrade performance during hour $2$ when compared to the $30$ minute window (or even the standard fluid). This happens because there is an excess supply of cars during the first two hours, e.g.\ availabilities are almost $100\%$. This excess allows one to act suboptimally in hour 2 to preposition cars for hour 3 without hurting performance.

\begin{table}[h!]
\centering
\begin{tabular}{|c|c|c|c|c|c|c|c|}
\hline
           & Hour 1 &Hour 2 &Hour 3 &Hour 4 & 4-hour Total\\ \hline
Standard Fluid   &0.995   &0.988   &0.705   &0.759   &0.8\\ \hline
JLCR-0.5   &1   &1   &0.715   &0.717   &0.79\\ \hline
SW   &1   &1   &0.67   &0.641   &0.745\\ \hline
$T$-Lookahead, $T=0.5$   &0.995   &0.991   &0.759   &0.771   &0.824\\ \hline
$T$-Lookahead, $T=0.75$   &0.995   &0.993   &0.788   &0.779   &0.838\\ \hline
\end{tabular}
\caption{Comparing the performance of the lookahead policy ($N=2000$). The time unit for $T$ is hours and the time-discretization $\Delta$ is chosen to be one minute. Based on the utility function in \eqref{eq:obj_data}, the values displayed are the fraction of fulfilled ride requests.  We used constant travel times, but the results remain qualitatively unchanged if the travel times are exponentially distributed.} \label{tab:CounterTraffic1000}
\end{table}


\subsection{\blue{Robustness}} 
\label{sec:robust}
The optimal fluid-based routing policy depends on perfect knowledge of system parameters $\lambda, P, \mu$. In practice, one uses estimates of these parameters, and the estimates contain noise. In this section we provide several numerical examples to get a sense of the robustness of our optimal fluid-based routing policy. The setup is the following: given a set of true parameters $\lambda, \mu, P$ and the utility function in \eqref{eq:obj_data}, we compute the optimal routing policy $Q^{*}$. 
Let  $\sigma > 0$, and let $\eta_{ij}$ and $\xi_{ij}$ be i.i.d.\ random variables that equal $1$ or $-1$ with equal probability. To simulate noisy estimates of passenger request rates along a route $\lambda_{i} P_{ij}$ and mean travel times $1/\mu_{ij}$, we use
\begin{align*}
\lambda_i P_{ij} (1+\sigma\eta_{ij}), \quad \text{ and } \quad (1/\mu_{ij})(1+\sigma\xi_{ij}), \quad 1 \leq i,j \leq r,
\end{align*}
respectively. We use four sets of parameters: the first set corresponds to the $9$-region network \eqref{eq:pdidi}--\eqref{eq:lambdadidi}, and the other three sets correspond to the three parameters sets used with the $5$-region model in Section~\ref{sec:fiveregion}. Given a true parameter set and a realization of the $\eta_{ij}$ and $\xi_{ij}$'s, we   compute a routing matrix $\hat Q^*$ based on the fluid optimization with the noisy parameters. To obtain a measure of suboptimality, we then evaluate \eqref{eq:obj_data} with the true parameters $\lambda,\mu,P$ but with routing matrix $\hat Q^*$. \bluetwo{We treat the objective value under $\hat Q^*$ as a random variable (it depends on $\eta_{ij}$ and $\xi_{ij}$), and report its mean and standard deviation in Table~\ref{tab:robust}}. We see that $5\%$ and $10\%$ estimation errors lead to policies that are approximately $5\%$ and $10\%$ suboptimal, respectively. Outside these four examples, the sensitivity of the optimal solution will depend on the choice of parameters. However, the table gives us an idea of the magnitude of the suboptimality due to estimation error.



\begin{table}[h!]
\centering
\begin{tabular}{|c|c|c|c|c|c|c|c|c|c|}
\hline
\thead{ Parameter \\ Set} & \thead{Optimal \\ Fluid}  &  \thead{$\hat Q^*$ Performance \\$\sigma=0.05$\\(mean, std. dev.)} & \thead{Relative \\ Suboptimality  \\($\sigma=0.05$)}&   \thead{$\hat Q^*$ Performance \\($\sigma=0.1$)\\(mean, std. dev.)} & \thead{Relative \\ Suboptimality  \\($\sigma=0.1$)}\\ \hline
 5-region, 5-7pm  & 0.91 & (0.87,0.024) & 0.044 & (0.82,0.044)  & 0.099 \\ \hline
  5-region, 7-9pm & 0.92 & (0.88,0.015) & 0.043 & (0.84,0.028)  & 0.087 \\ \hline
 5-region, 9-11pm & 0.92 & (0.86,0.031) & 0.065 & (0.81,0.053)  & 0.120 \\ \hline
   9-region & 0.8403 & (0.818868, 0.01503) & 0.021478 & (0.800520, 0.029091)   & 0.039826 \\ \hline
\end{tabular} 
\caption{The third and fifth columns report the performance of the routing policy obtained from noisy parameters, $\hat Q^*$,  under the true parameters $\lambda,\mu,P$. We ran 1000 replications, where in each replication we generated a value for $\eta_{ij}$ and $\xi_{ij}$;  $2r^2$ random variables per replication in total.}
\label{tab:robust}
\end{table}

In the next section we introduce the fluid model together with the tools needed to prove Theorems~\ref{thm:asymptotic_utility} and \ref{thm:lb}.


%
%

\section{The Fluid Model} \label{sec:fluid}
In this section we study the fluid model of the ridesharing network in Section~\ref{sec:model} and prepare all the ingredients needed to prove Theorem~\ref{thm:asymptotic_utility}. We start by introducing the fluid model for static empty-car routing matrices $Q$, and establish process-level convergence in Theorem~\ref{thm:main_conv}. In Section~\ref{sec:equil_points}, we then characterize the fluid model's set of equilibria. This equilibrium behavior motivates the fluid-based optimization in Section~\ref{sec:model}. Theorem~\ref{thm:asymptotic_utility} then follows from a standard argument involving process-level convergence and the convergence of the fluid model to its equilibrium, and  we therefore relegate it to Appendix~\ref{app:main_thm_proof}.

Recall the primitive parameters $\lambda, \mu, P$, and assume a static empty-car routing matrix $Q$ is given.  Recall the set $\mathcal{T}$ defined in \eqref{eq:T}. Let
\begin{align*}
I_{i}^{(N)}(t) = \int_{0}^{t} 1 \big(E^{(N)}_{ii}(s) = 0\big) ds
\end{align*}
be the cumulative idle time of the single-server station corresponding to $E^{(N)}_{ii}$. The following is a process-level convergence result for  $(\bar E^{(N)}, \bar F^{(N)})$, and is proved in Appendix~\ref{app:network_description}.
\begin{theorem} \label{thm:main_conv}
Assume $\big(\bar E^{(N)}(0), \bar F^{(N)}(0) \big) \Rightarrow (e(0), f(0)) \in \mathcal{T}$ as $N \to \infty$. \blue{There exists a unique solution} $(e,f): \R_+ \to \mathcal{T}$ and $u : \R_+ \to \R_+^{r}$ to the dynamical system
\begin{align}
&f_{ij}(t) = f_{ij}(0) + \lambda_{i} P_{ij}\big( t -  u_i(t)\big) - \mu_{ij} \int_{0}^{t}f_{ij}(s)ds, \quad &1 \leq i,j \leq r, \label{eq:ufluid_fij} \\
&e_{ij}(t) = e_{ij}(0) - \mu_{ij} \int_0^t e_{ij}(s)ds + Q_{ij} \sum_{k = 1}^{r} \mu_{ki} \int_{0}^{t}f_{ki}(s)ds, \quad &1 \leq i \neq j \leq r, \label{eq:ufluid_eij} \\
&e_{ii} (t) = e_{ii}(0) -\lambda_i \big( t -  u_i(t)\big) &  \notag \\
& \hspace{2cm} + \sum_{\substack{j = 1 \\ j \neq i}}^{r} \mu_{ji} \int_{0}^{t}e_{ji}(s)ds+ Q_{ii} \sum_{j = 1}^{r} \mu_{ji} \int_{0}^{t}f_{ji}(s)ds, \quad &1 \leq i \leq r, \label{eq:ufluid_eii} \\
&u(t) \text{ is non-decreasing with } u(0) = 0, \text{ and }  \int_{0}^{\infty} e_{ii}(s) d u_i(s) = 0 &\text{ for all $1 \leq i \leq r$}. \label{eq:ufluid_regulator}
\end{align}
\blue{Furthermore,} for all $T \geq 0$,
\begin{align}
\lim_{N \to \infty} \norm{\big(\bar E^{(N)}, \bar F^{(N)}, I^{(N)} \big) - (e, f, u)}_{T} = 0 \label{eq:uber_finite_convergence}
\end{align}
almost surely.
\end{theorem}
\begin{remark}
Theorem~\ref{thm:asymptotic_utility} assumes that $P_{ij} > 0$ for all $i,j$ and $Q_{ii}>0$ for all $i$. This assumption is not used in the proof of Theorem~\ref{thm:main_conv}, but appears in Section~\ref{sec:equil_points}.
\end{remark}
We refer to equations
  (\ref{eq:ufluid_fij})-(\ref{eq:ufluid_regulator}) as the fluid model
  of the ridesharing network, and to  $\big(e,f,u\big)$ as the fluid
  analog of $\big(\bar E^{(N)}, \bar F^{(N)},
  I^{(N)}\big)$.  Note that $P_{ij} = 0$ implies $f_{ij}(t) \equiv 0$,
  and $Q_{ij} = 0$ for $i \neq j$ implies $e_{ij}(t) \equiv 0$. It
  will come in handy later on to know that, $e(t)$, $f(t)$, and $u(t)$
  are Lipschitz continuous. To see why, observe that for any $\epsilon
  > 0$, \eqref{eq:uber_finite_convergence} says that we can choose $N$
  large enough such that
\begin{align*}
\abs{u(t) - u(s)} \leq&\ \abs{u(t) - I^{(N)}(t)} + \abs{I^{(N)}(t) - I^{(N)}(s)} + \abs{I^{(N)}(s) - u(s)} \\
\leq&\ \epsilon + \abs{t-s} + \epsilon,
\end{align*}
where in the second inequality we used the definition of $I^{(N)}$ to get
\begin{align*}
\abs{I^{(N)}(t) - I^{(N)}(s)} \leq \abs{t-s}, \quad 0 \leq s,t < \infty.
\end{align*}
Hence,
\begin{align}
\abs{u(t) - u(s)} \leq \abs{t-s}, \quad 0 \leq s,t < \infty. \label{eq:ulip}
\end{align}
Combining \eqref{eq:ufluid_fij}--\eqref{eq:ufluid_eii} with \eqref{eq:ulip} and the fact that $\big(e(t),f(t)\big)$ is bounded \bluetwo{(because it is in $ \mathcal{T}$)}, we deduce that both $e(t)$ and $f(t)$ are also Lipschitz-continuous. The next section considers the equilibrium behavior of the fluid model.

\subsection{Equilibrium Points} \label{sec:equil_points}
We begin with an informal discussion on how to characterize the equilibrium of the fluid model. In equilibrium, we expect $\big(\dot e(t), \dot f(t)\big)$ to equal zero. Taking derivatives in \eqref{eq:ufluid_fij}--\eqref{eq:ufluid_eii} and setting the left hand sides to zero gives us
\begin{align}
&0 =\lambda_{i} P_{ij}\big( 1 -  \dot u_i(t)\big) - \mu_{ij} f_{ij}(t), \quad &1 \leq i,j \leq r, \label{eq:fij_dot}\\
&0 =  - \mu_{ij} e_{ij}(t) + Q_{ij} \sum_{k = 1}^{r} \mu_{ki} f_{ki}(t), \quad &1 \leq i \neq j \leq r, \label{eq:eij_dot}\\
&0 = -\lambda_i \big( 1 -  \dot u_i(t)\big)  + \sum_{\substack{j = 1 \\ j \neq i}}^{r} \mu_{ji} e_{ji}(t)+ Q_{ii} \sum_{j = 1}^{r} \mu_{ji} f_{ji}(t), \quad &1 \leq i \leq r. \label{eq:eii_dot}
\end{align}
The derivative $\dot u(t)$ above exists for almost every $t$  because $u(t)$ is Lipschitz continuous by \eqref{eq:ulip}. Next, by adding up all the terms in \eqref{eq:ufluid_fij}--\eqref{eq:ufluid_eii} we can see that the total amount of fluid in the system always equals one, or \blue{
\begin{align}
\big( e(t), f(t) \big) \in \mathcal{T}. \label{eq:fluidmass}
\end{align}
}
Lastly, we combine Lipschitz continuity of $u(t)$ together with \eqref{eq:ufluid_regulator} to see that
\begin{align*}
\int_{0}^{\infty} e_{ii}(s) du_i(s) = \int_{0}^{\infty} e_{ii}(s) \dot u_i(s) ds = 0.
\end{align*}
Since $e_{ii}(t) \geq 0$ and $\dot u_i(t) \geq 0$ (because $u(t)$ is an increasing function), for all $t \geq 0$ where $\dot u_{i}(t)$ exists,
\begin{align}
e_{ii}(t) \dot u_{i}(t) = 0. \label{eq:ucomplement}
\end{align}
Therefore, we expect the equilibrium of the fluid model (if it exists) to satisfy \eqref{eq:fij_dot}--\eqref{eq:ucomplement}. The following lemma addresses the issue of existence of an equilibrium. In the lemma, we use $\bar a_i \in [0,1]$ as a placeholder for $1- \dot u_i(t)$ to represent the equilibrium server utilization at the station corresponding to $e_{ii}$. We let $\bar a$ be the $r$-dimensional vector whose components are $\bar a_i$. \blue{We can now see that \eqref{eq:fij_dot}--\eqref{eq:ucomplement} are exactly the constraints in the fluid-based optimization problem (\hyperref[eqlbm: fullcsv]{4b})--(\hyperref[eqlbm:prob]{4g}).}


\begin{lemma} \label{lem:equil}
Assume $P_{ij} > 0$ and $Q_{ii} > 0$ for all $i,j= 1, \ldots, r$. The system of equations
\begin{subequations}
\begin{align}
&\ \lambda_{i} P_{ij} \bar a_i  = \mu_{ij} \bar f_{ij}, \quad &1 \leq i, j, \leq r, \label{eq:one} \\
&\ \mu_{ij}\bar e_{ij} = Q_{ij} \sum_{k = 1}^{r} \mu_{ki}\bar f_{ki}, \quad &1 \leq i \leq r,\ j \neq i, \label{eq:two} \\
&\ \lambda_i \bar a_i  = \sum_{\substack{k = 1 \\ k \neq i}}^{r} \mu_{ki} \bar e_{ki} + Q_{ii} \sum_{k = 1}^{r} \mu_{ki} \bar f_{ki}, \quad &1 \leq i \leq r, \label{eq:three} \\
&(1 - \bar a_{i}) \bar e_{ii} = 0, \quad 1 \leq i \leq r, \label{eq:ycomp}\\
& (\bar e, \bar f) \in \mathcal{T}, \quad  \bar a \in [0,1]^{r} \label{eq:balance}
\end{align}
\end{subequations}
has at least one solution $(\bar e, \bar f, \bar a)$. Multiple solutions may exist, but the only difference between solutions are in the components $\bar e_{ii}$ for $i$ such that $\bar a_{i} = 1$. In other words, the components $\bar a, \bar f,$ $\bar e_{ij}$ for $i \neq j$, and $\bar e_{ii}$ for $i$ such that $\bar a_i < 1$, are identical across all solutions. 
\end{lemma}
Define
\begin{align}
 \mathcal{E} = \big\{(\bar e, \bar f) \in \mathcal{T} : \exists\ \bar a \in [0,1]^{r} \text{ such that $(\bar e, \bar f, \bar a)$ solves \eqref{eq:one}--\eqref{eq:ycomp}}\big\}. \label{eq:E}
\end{align}
Lemma~\ref{lem:equil} implies that we can associate a unique $\bar a$ to the set $\mathcal{E}$. Furthermore, the lemma implies that the quantity $\bar m$, defined as 
\begin{align}
\bar m = \sum_{i: \bar a_{i} = 1} \bar e_{ii} = 1 - \sum_{i=1}^{r} \sum_{j=1}^{r} \bar f_{ij} - \sum_{i=1}^{r} \sum_{\substack{j=1\\ j \neq i}}^{r} \bar e_{ij} \label{eq:mdef}
\end{align}
is unique. In light of the following theorem, we refer to $\mathcal{E}$ as the equilibrium set. The theorem is proved in Appendix~\ref{app:stability} and is the final ingredient to prove Theorem~\ref{thm:asymptotic_utility}.
\begin{theorem} \label{thm:stability}
Let $\big(e(t),f(t),u(t)\big)$ be the unique solution to \eqref{eq:ufluid_fij}--\eqref{eq:ufluid_regulator} with initial condition $(e(0), f(0)) \in \mathcal{T}$. Assume $P_{ij} > 0$ and $Q_{ii} > 0$ for all $i,j= 1, \ldots, r$. Then for any $\epsilon > 0$, there exists a $T > 0$ such that
\begin{align*}
\inf_{x \in \mathcal{E}} \abs{\big(e(t),f(t)\big) - x} < \epsilon, \quad t \geq T.
\end{align*}
\end{theorem}


\begin{remark}In Lemma~\ref{lem:equil}  and Theorem~\ref{thm:stability}, we make use of the assumption that $P_{ij} > 0$ for all $i,j$ and $Q_{ii} > 0$ for all $i$. We expect that the result holds if the CTMC has a single recurrent class, but our extra assumption greatly facilitates exposition. This is especially true in the proof of Theorem~\ref{thm:stability}, which is already rather cumbersome.
\end{remark}

\section{Conclusions} \label{sec:conclusion}
This paper considered empty-car routing in a ridesharing network under a regime where both supply and demand for cars tend to infinity, and provided a comprehensive analysis of the design of an optimal empty-car routing policy based on asymptotic fluid analysis. A numerical study using real-world ridesharing data confirmed the effectiveness of our solution. Our paper is only a first step to understand ridesharing networks, and poses some interesting problems \delete{are} for future research directions:

\textbf{Decentralized routing:}
Our routing policy is a centralized solution that assumes the ridesharing platform has full control over its empty cars, which provides a best-case benchmark. One future research topic is the design of decentralized incentive mechanisms for achieving the routing probabilities of the centralized solution. One can use the centralized routing policy as a benchmark to quantify the efficiency of a decentralized mechanism.

\textbf{Time-varying parameters:} \blue{Our $T$-lookahead policy is just one heuristic for dealing with time varying parameters, but it would be very interesting to be able to say something rigorous about time-varying policies. Along this line, studying the transient control problem would also be of interest (given the long time it takes the fluid to converge to equilibrium from certain initial conditions). Even studying the fluid transient control problem is non-trivial, because the fluid model is a non-linear dynamical system.}


\textbf{Robust routing policies:} We witnessed in Figure~\ref{fig:comp-didi} that certain state-dependent policies can attain performance levels close to those of the optimal static policy, yet have the benefit of not requiring explicit knowledge of system parameters. Now that a benchmark has been established, can we find state-dependent policies that are provably asymptotically optimal, yet rely as little as possible on primitive system parameters?

\ECSwitch


\ECHead{Appendix}

\section{Process Level Convergence for Closed Networks of Single and Infinite Server Stations}
\label{app:network_description}
In this section we consider closed queueing networks with exponential service times that consist of both single-server and infinite-server stations. We prove that as both the number of jobs in the network, and the service rates at the single server stations increase, the appropriately scaled queue length process converges to a fluid limit. The ridesharing model in Section~\ref{sec:model} is a special type of such networks, meaning that Theorem~\ref{thm:main_conv} will be a special case of the results here. We remark that process-level convergence for the class of networks considered here is a consequence of the results in \cite{Kric1992}. However, in that paper the limiting process is defined as the solution to a differential inclusion, and our proof technique is sufficiently different to merit a separate write-up.

We consider a closed queueing network with $N >0$ jobs and $J$ stations, consisting of both single-server and infinite-server stations.  We let the stations be indexed by the set $\mathcal{J} = \{1, \ldots, J\}$. Let $\mathcal{I} \subset \mathcal{J}$ and $\mathcal{S} \subset \mathcal{J}$ be the \textit{non-empty} index sets corresponding to the infinite server stations, and single server stations, respectively. To describe the network dynamics, we introduce the following primitives. Let $Q^{(N)}(0) \in \Z_+^{J}$ with $\sum_{i \in \JJJ} Q^{(N)}_i(0) = N$ be the vector representing the initial job distribution in the network. To keep track of service completions at each stations, we let
\begin{align*}
S_i = \{S_i(t), t \geq 0\}, \quad i \in \JJJ
\end{align*}
be a collection of unit-rate Poisson processes with $S_i$ independent of $S_j$ for $i \neq j$. Let $\lambda, \mu \in \R_+^J$ be two vectors with $\lambda_i = 0$ for $i \in \III$ and $\mu_i = 0$ for $i \in \SSS$, which we will use to represent the service rates at different stations. We assume that the service rate of each server at station $i \in \JJJ$ is
\begin{align*}
&\ N \lambda_i, \quad i \in \SSS \\
&\ \mu_i, \quad i \in \III.
\end{align*}
Let
\begin{align*}
\Big\{ \big(\Phi_{i1}(n), \ldots, \Phi_{iJ}(n)\big) \in \Z_+^J,\ n \in \Z_+\Big\}, \quad i \in \JJJ
\end{align*}
be a collection of routing processes defined as follows. For each $n \in \Z_+$ and each $i \in \JJJ$, the vector
\begin{align*}
\big(\Phi_{i1}(n), \ldots, \Phi_{iJ}(n)\big) = \sum_{m = 1}^{n} \phi_{i}(m),
\end{align*}
where $\big\{\phi_{i}(m) \in \{0, 1\}^J\big\}_{m=1}^{\infty}$ is a sequence of i.i.d. random variables with
\begin{align*}
\Prob( \phi_{i}(1) = e^{(j)} ) = R_{ij}, \quad i, j \in \JJJ.
\end{align*}
Furthermore, the sequences $\big\{\phi_{i}(m) \big\}_{m=1}^{\infty}$ and $\big\{\phi_{j}(m) \big\}_{m=1}^{\infty}$ are assumed to be independent for $i \neq j$.  Let $R$ be the routing probability matrix whose $i,j$th entry is $R_{ij}$, and observe that it is a column stochastic matrix.

Using an inductive argument similar to the proof of \cite[Theorem 2.1]{ChenMand1991}, one can show that there exists a unique process
\begin{align*}
Q^{(N)} = \{Q^{(N)}(t) = (Q_1^{(N)}(t), \ldots, Q_J^{(N)}(t)), t \geq 0\}
\end{align*}
satisfying
\begin{align}
Q_i^{(N)}(t) =&\ Q_{i}^{(N)}(0) - S_{i}\big(N\lambda_i T_{i}^{(N)}(t) \big) + \sum_{j \in \SSS} \Phi_{j i}\Big( S_{j}\big(N\lambda_j T_{j}^{(N)}(t)\big) \Big)  \notag \\
&\ + \sum_{k \in \III} \Phi_{k i}\bigg( S_{k}\Big(\mu_k\int_{0}^{t} Q_{k}^{(N)}(s) ds \Big) \bigg), \quad i \in \SSS, \label{eq:qeq1} \\
Q_i^{(N)}(t) =&\ Q_{i}^{(N)}(0) - S_{i}\Big(\mu_i\int_{0}^{t} Q_{i}^{(N)}(s) ds \Big) + \sum_{j \in \SSS} \Phi_{j i}\Big(S_{j}\big(N\lambda_j T_{j}^{(N)}(t)\big) \Big)  \notag \\
&\ + \sum_{k \in \III} \Phi_{k i}\bigg( S_{k}\Big(\mu_k \int_{0}^{t} Q_{k}^{(N)}(s) ds \Big) \bigg), \quad i \in \III, \label{eq:qeq2}
\end{align}
where
\begin{align*}
T_{i}^{(N)} = \bigg\{ T_{i}^{(N)}(t) = \int_{0}^{t} 1( Q_i^{(N)}(s) > 0) ds \bigg\}, \quad i \in \SSS,
\end{align*}
is the cumulative busy time process of the server at each single server station. At any time $t \geq 0$, $Q^{(N)}_i(t)$ is the job count at station $i \in \JJJ$. It is a straightforward exercise to verify that $Q^{(N)}$ satisfies the Markov property and is therefore a CTMC. Furthermore, the CTMC $(E^{(N)}, F^{(N)})$ introduced in Section~\ref{sec:model} is a special case of $Q^{(N)}$.

To write \eqref{eq:qeq1}--\eqref{eq:qeq2} in a form that is more convenient for analysis, for any $t \geq 0$ let us define
\begin{align*}
\widehat S_i(t) =&\ S_i(t) - t, \quad i \in \JJJ,\\
\widehat \Phi_{ij}(n) =&\ \Phi_{ij}(n) - R_{ij} n, \quad i,j \in \JJJ,\ n \in \Z_+,
\end{align*}
and
\begin{align}
\widehat M_{i}^{(N)}(t) =&\ - \widehat S_{i}\big(N \lambda_i T_{i}^{(N)}(t)\big) + \sum_{j \in \SSS} \bigg[\widehat \Phi_{ji}\Big( S_j\big(N \lambda_j T^{(N)}_j(t)\big) \Big) + R_{ji}\widehat S_{j}\big(N\lambda_j T^{(N)}_j(t)\big) \bigg] \notag \\
&+ \sum_{k \in \III} \bigg[ \widehat \Phi_{ki}\bigg( S_{k}\Big(\mu_k\int_{0}^{t} Q_{k}^{(N)}(s) ds \Big) \bigg) + R_{ki} \widehat S_{k}\Big(\mu_k\int_{0}^{t} Q_{k}^{(N)}(s) ds \Big) \bigg], \quad i \in \SSS, \label{eq:mart1}\\
\widehat M_{i}^{(N)}(t) =&\ -\widehat S_{i}\Big(\mu_i\int_{0}^{t} Q_{i}^{(N)}(s) ds \Big) + \sum_{j \in \SSS} \bigg[\widehat \Phi_{j i}\Big(S_{j}\big(N\lambda_j T_{j}^{(N)}(t)\big) \Big) + R_{ji}\widehat S_{j}\big(N\lambda_j T_{j}^{(N)}(t)\big) \bigg] \notag  \\
&+ \sum_{k \in \III} \bigg[\widehat \Phi_{k i}\bigg( S_{k}\Big(\mu_k \int_{0}^{t} Q_{k}^{(N)}(s) ds \Big) \bigg) + R_{ki} \widehat S_{k}\Big(\mu_k \int_{0}^{t} Q_{k}^{(N)}(s) ds \Big) \bigg], \quad i \in \III, \label{eq:mart2}
\end{align}
and let $\widehat{M}^{(N)}(t)$ be the vector whose components are $\widehat{M}_i^{(N)}(t)$. For $t \geq 0$, we also define
\begin{align}
I_i^{(N)}(t) = 0 \text{ for $i \in \III$} \quad \text{ and } \quad I_{i}^{(N)}(t) = t - T_{i}^{(N)}(t) \text{ for $i \in \SSS$}, \label{eq:idle}
\end{align}
and let $I^{(N)} = \{ I^{(N)}(t) \in \R_+^{J}, t \geq 0 \}$. Then for $i \in \SSS$, $I_{i}^{(N)}(t)$ represents the cumulative idle time up to time $t$. Setting
\begin{align*}
\bar Q^{(N)}(t) = \frac{1}{N} Q^{(N)}(t) \quad \text{ and } \quad \bar M^{(N)}(t) = \frac{1}{N} \widehat M^{(N)}(t),
\end{align*}
we from \eqref{eq:qeq1}--\eqref{eq:mart2} that
\begin{align}
\bar Q_{i}^{(N)}(t) =&\ \bar Q_{i}^{(N)}(0) + \bar M_{i}(t) + \Big(\sum_{j \in \SSS}R_{ji} \lambda_j - \lambda_i \Big)t + \sum_{k \in \III}  R_{ki} \mu_k \int_{0}^{t} \bar Q_{k}^{(N)}(s) ds \notag \\
&+ \lambda_i I_i^{(N)}(t) - \sum_{j \in \SSS}R_{ji} \lambda_j I_j^{(N)}(t), \quad i \in \SSS, \label{eq:qfinal1} \\
\bar Q_{i}^{(N)}(t)=&\ \bar Q_{i}^{(N)}(0)+ \bar M_{i}(t) + \Big(\sum_{j \in \SSS} R_{ji} \lambda_j\Big)t - \mu_i \int_{0}^{t} \bar Q_{i}^{(N)}(s) ds \notag \\
&+  \sum_{k \in \III} R_{ki} \mu_k \int_{0}^{t} \bar Q_{k}^{(N)}(s) ds - \sum_{j \in \SSS} R_{ji} \lambda_j I_j^{(N)}(t), \quad i \in \III.\label{eq:qfinal2}
\end{align}
In the next section, we describe the fluid model to which the process $\bar Q^{(N)}$ will converge to as $N \to \infty$.

\subsection{The Fluid Model}
\label{sec:fluid_model}
Recalling that $\mu_i = 0$ for $i \in \SSS$ and $\lambda_i = 0$ for $i \in \III$, we set
\begin{align}
M = \text{diag}(\mu) \quad \text{ and } \quad \Lambda = \text{diag}(\lambda). \label{eq:diags}
\end{align}
We also define the $J \times J$ matrix $\tilde R$ by setting
\begin{align*}
\tilde R_{ij} =&\ R_{ij}, \quad i \in \SSS,\\
\tilde R_{ij} =&\ 0, \quad i \in \III.
\end{align*}
That is $\tilde R$ is the matrix $R$ with all rows corresponding to infinite server stations being set to zero. Since $\III \neq \emptyset$, the matrix $\tilde R^T$ is sub-stochastic. The following lemma is proved in Section~\ref{sec:lemmaproof}.
\begin{lemma} \label{lem:continuity}
For each $x \in \D_1^J$, there exists a unique $(q,v) \in \D^{2J}$, with $q(t) \in \R_+^J$ and $v(t) \in \R_+^J$ for all $t \geq 0$, such that
\begin{align}
&q(t) = x(t) - (I - R^{T}) M\int_{0}^{t}q(s) ds + (I-\tilde R^{T}) v(t) \label{eq:lem1} \\
&q(t) \geq 0, \quad t \geq 0,  \label{eq:lem2}\\
&v(\cdot) \text{ is  non-decreasing with } v(0)=0, \label{eq:lem3}\\
&\int_{0}^{\infty} q_i(s) dv_i(s) = 0, \quad i \in \JJJ. \label{eq:lem4}
\end{align}
Furthermore, the map $\Upsilon: \D_1^J \to \D^{2J}$ given by $\Upsilon(x) = (q, v)$ is well-defined and is Lipschitz-continuous, in the sense that for any $x, \tilde x \in \D^J_1$, and any $T > 0$, there exists a constant $c_{\Upsilon}^{T}$ such that
\begin{align}
\norm{\Upsilon(x) - \Upsilon(\tilde x)}_T \leq c_{\Upsilon}^{T}\norm{x - \tilde x}_T, \label{eq:ups_lipschitz}
\end{align}

\end{lemma}

\begin{theorem} \label{thm:process_convergence}
Assume $\tilde Q^{(N)}(0) \to a$ as $N \to \infty$ for some $a \in [0,1]^J$ with $\sum_{i=1}^{J} a_i = 1$. Let $e \in \R^J$ be the vector of ones, and let $\gamma : \R_+ \to \R_+$ be the identity map defined by $\gamma(t) = t$. Set
\begin{align}
(q,v) = \Upsilon\Big(a + \big((R^T - I) \Lambda e\big)\gamma\Big). \label{eq:thm_statement}
\end{align} Then for any $T > 0$,
\begin{align*}
&\lim_{n \to \infty} \norm{\tilde Q^{(N)}(t) - q(t) }_T = 0, \\
&\lim_{n \to \infty} \norm{I^{(N)}(t) - v(t) }_T = 0.
\end{align*}
\end{theorem}
\proof{Proof of Theorem~\ref{thm:process_convergence}.}
From \eqref{eq:qfinal1}--\eqref{eq:qfinal2} we see that
\begin{align*}
\big(\bar Q^{(N)}, I^{(N)}\big)= \Upsilon\Big(\bar Q^{(N)}(0) + \bar M^{(N)} + \big((R^T - I) \Lambda e\big)\gamma\Big),
\end{align*}
where $\bar M^{(N)} = \{ \bar M^{(N)}(t) \in \R^J, t \geq 0 \}$. Suppose we knew that for every $T \geq 0$,
\begin{align}
\lim_{N \to \infty} \norm{\bar M^{(N)}}_{T} = 0  \quad \text{ almost surely.} \label{eq:mart_vanish}
\end{align}
Then the continuous mapping theorem \cite{Bill1999}, together with \eqref{eq:ups_lipschitz} would imply Theorem~\ref{thm:process_convergence}. The proof of \eqref{eq:mart_vanish} involves a standard argument using the functional strong law of large numbers (FSLLN), and is therefore omitted. For an example of such an argument, see the proof of (5.6) in \cite{DaiHeTezc2010}.
\endproof
\begin{remark}
Since $I^{(N)}_i(t) = 0$ for $i \in \III$ and $t \geq 0$, Theorem~\ref{thm:process_convergence} implies
\begin{align}
v_i(t) = 0, \quad i \in \III. \label{eq:single_zero}
\end{align}
Establishing \eqref{eq:single_zero} by relying on convergence of $\big(\bar Q^{(N)}, I^{(N)}\big)$ to $(q,v)$ may seem strange, because \eqref{eq:single_zero} should be a standalone property of $\Upsilon\Big(a + \big((R^T - I) \Lambda e\big)\gamma\Big)$. Indeed, it is possible to establish \eqref{eq:single_zero} using a direct argument that relies on Proposition 1 of \cite{Reim1984}. However, we avoid using said argument because it is significantly longer.
\end{remark}
We immediately see that Theorem~\ref{thm:main_conv} is a special case of Theorem~\ref{thm:process_convergence}. The rest of this section is devoted to proving Lemma~\ref{lem:continuity}.


\subsection{Proof of Lemma~\ref{lem:continuity}}
\label{sec:lemmaproof}
In order to prove Lemma~\ref{lem:continuity}, we first need to introduce the Skorohod problem. Let $\tilde Q$ be a $J \times J$ column sub-stochastic matrix with non-negative entries. For any $x \in \D^J_{0+}$,  let $(z, y) \in \D^{2J}$ with $z(t) \in \R^J_+$ and $y(t) \in \R^J_+$ for all $ t \geq 0$ be the solution to
\begin{align}
&\ z = x + (I - \tilde Q^{T}) y, \label{eq:SP1}\\
&\ z \geq 0, \label{eq:SP2}\\
&\ y(\cdot) \text{ is non-decreasing and } y(0) = 0,  \label{eq:SP3}\\
&\ \int_{0}^{\infty} z_i(s) d y_i(s) = 0, \quad  1 \leq i \leq J. \label{eq:SP4}
\end{align}
Existence and uniqueness of $(z,y)$ was proved in \cite{HarrReim1981} when $x$ is continuous, but the arguments there hold for $x \in \D^J_{0+}$ as well. We refer to \eqref{eq:SP1}--\eqref{eq:SP4} as the Skorohod problem associated with $(x, \tilde Q)$, and write SP$(x,\tilde Q)$ for short. We refer to $(z,y) \in \D^{2J}$ as the solution to SP$(x,\tilde Q)$ if it satisfies \eqref{eq:SP1}--\eqref{eq:SP4}. Furthermore, for any $x \in \D_{0+}^J$  we define the Skorohod map $\Psi^{\tilde Q}: \D_{0+}^J \to \D^{2J}$ by
\begin{align*}
\Psi^{\tilde Q}(x) = \big( \Psi_{z}^{\tilde Q}(x), \Psi_{y}^{\tilde Q}(x)\big) = (z,y),
\end{align*}
where $(z,y)$ is the solution of SP$(x,\tilde Q)$.
From \eqref{eq:SP1} it is clear that
\begin{align}
\Psi^{\tilde Q}_{z}(x) = x + (I-\tilde Q^{T})\Psi^{\tilde Q}_{y}(x), \quad x \in \D^J,\ x(0) \geq 0. \label{eq:smap}
\end{align}
Both $\Psi^{\tilde Q}_{z}$ and $\Psi^{\tilde Q}_{y}$ are Lipschitz-continuous, in the sense that for any $T > 0$, there exists constants $c_{z}^{T}, c_{y}^{T} > 0$, which depend on $\tilde Q$, such that for any $x, \tilde x \in \D_{0+}^J$,
\begin{align*}
\norm{\Psi^{\tilde Q}_{y}(x) - \Psi^{\tilde Q}_{y}(\tilde x)}_{T} \leq&\ c_{y}^{T} \norm{x - \tilde x}_T, \\
\norm{\Psi^{\tilde Q}_{z}(x) - \Psi^{\tilde Q}_{z}(\tilde x)}_{T} \leq&\ c_{z}^{T} \norm{x - \tilde x}_T,
\end{align*}
where $\norm{\cdot}_T$ is defined in \eqref{eq:norm}. This was established in \cite[p. 305]{HarrReim1981} when both $x$ and $\tilde x$ are continuous, but the argument used there holds for $x, \tilde x \in \D_{0+}^J$ as well. We are now ready to prove Lemma~\ref{lem:continuity}.

\proof{Proof of Lemma~\ref{lem:continuity}.}
Fix $x \in \D_1^J$ and omitting the superscript $\tilde R$, let $(\Psi_{z}, \Psi_{y})$ be the Skorohod map associated with SP$(x,\tilde R)$. Define the integral operator $\alpha:\D_{1}^J \times \D^J \to \D^J$ by
\begin{align}
\alpha(x, w)(t) = x (t) - (I-R^T) M \int_{0}^{t} w(s) ds, \label{eq:alpha}
\end{align}
and consider the integral equation
\begin{align}
&q = \Psi_{z}(\alpha(x, q)) = \alpha(x,q) + (I-\tilde R^{T}) \Psi_{y}(\alpha(x,q)) \label{eq:INT1} \\
&q \in \D^J, \label{eq:INT2}
\end{align}
where the second equality in \eqref{eq:INT1} follows from \eqref{eq:smap}.
Provided \eqref{eq:INT1}--\eqref{eq:INT2} has a unique solution $q$, we can set $v = \Psi_{y}(\alpha(x,q))$ and observe that by definition of $\alpha$ and $(\Psi_{z}, \Psi_{y})$,
\begin{align*}
(q,v) = \Upsilon(x).
\end{align*}
Hence, we now establish the existence and uniqueness of a solution $q$ to \eqref{eq:INT1}--\eqref{eq:INT2}. Construct a sequence $\{ q^{n} \in \D^J \}_{n=0}^{\infty}$ by letting
\begin{align*}
q^{0}(t) \equiv&\ x(0), \quad t \geq 0,\\
q^{n+1}(t) =&\ \Psi_{z}\big(\alpha(x, q^{n})\big)(t), \quad t \geq 0.
\end{align*}
Observe that $q^{n}(0) = x(0)$ for all $n \geq 0$. We first show that for any $T > 0$, this sequence is a Cauchy sequence in the Hilbert space $\big(\D^J[0,T], \norm{\cdot}_T\big)$. Let $\bar{\mu} = \max_{i \in \III} \{\mu_j\}$ (remembering that $\mu_i = 0$ for $i \in \SSS$), and observe that
\begin{align}
\norm{q^{n+1}- q^{n}}_T \leq&\ c_{z}^{T}\norm{\alpha(x,q^{n}) - \alpha(x,q^{n-1})}_T \notag \\
=&\ c_{z}^{T}\max_{i \in \JJJ}\sup_{0 \leq t \leq T}\abs{\big(\alpha(x,q^{n})\big)_i(t) - \big(\alpha(x,q^{n-1})\big)_i(t)} \notag \\
\leq&\ c_{z}^{T} \bar{\mu} J  \max_{i \in \III}\bigg\{\int_{0}^{T} \abs{q_i^{n}(s) - q_i^{n-1}(s)} ds\bigg\}\notag  \\
\leq&\ c_{z}^{T} \bar{\mu} J  \int_{0}^{T} \max_{i \in \III}\abs{q_i^{n}(s) - q_i^{n-1}(s)} ds\notag  \\
\leq&\ \frac{(c_{z}^{T} \bar{\mu}J  T)^{n} }{n!}\norm{q^{1} - q^{0}}_T. \label{eq:lip_argument}
\end{align}
The first inequality follows from the Lipschitz property of $\Psi_{z}$, the second inequality is from the form of $\alpha$, and the last inequality follows by recursion. From this point it is not hard to conclude (see for instance (11.22) of \cite{MandMassReim1998}) that $\{ q^{n}\}_{n=0}^{\infty}$ is a Cauchy sequence in $\big(|D^J[0,T], \norm{\cdot}_T\big)$ for each $T > 0$. Therefore, $q^n$ converges to some limit $q \in \D^J[0,T]$ that satisfies \eqref{eq:INT1}. Since the choice of $T> 0$ was arbitrary, we have proved existence of a solution to \eqref{eq:INT1}--\eqref{eq:INT2}. Uniqueness can be argued by taking two potential solutions $q$ and $\tilde q$, and applying the chain of arguments in \eqref{eq:lip_argument} with $\norm{q - \tilde q}_{T}$ on the left hand side there.

It remains to prove the Lipschitz-continuity of $\Upsilon$. Fix any $x, \tilde x \in \D^J_1$, and set
\begin{align*}
(q,v) = \Upsilon(x) \quad \text{ and } \quad (\tilde q, \tilde v)  = \Upsilon( \tilde x).
\end{align*}
Repeating the logic used to obtain \eqref{eq:lip_argument}, we see that for any $n \geq 1$,
\begin{align*}
\norm{q - \tilde q}_T =&\ \norm{\Psi_{z}\big( \alpha(x,q) \big) - \Psi_{z}\big( \alpha(\tilde x, \tilde q) \big)}_T \\
\leq&\ c_{z}^{T}\norm{x - \tilde x}_T  +
c_{z}^{T} \bar{\mu} J  \int_{0}^{T} \max_{i \in \III}\abs{q_i(s) - \tilde q_i(s)} ds\notag  \\
\leq&\ c_{z}^{T}\norm{x - \tilde x}_T \sum_{k = 0}^{n-1} \frac{\big(c_z^{T} \bar{\mu} JT\big)^k}{k!} + \frac{\big(c_z^{T} \bar{\mu} JT\big)^{n}}{n!} \norm{q - \tilde q}_T,
\end{align*}
where the last inequality follows by recursion. Choosing $n$ large enough so that $\frac{\big(c_z^{T} \bar{\mu} JT\big)^{n}}{n!} < 1$, we conclude the existence of a constant $c_{\Upsilon,q}^{T} > 0$ such that
\begin{align}
\norm{q - \tilde q}_T \leq c_{\Upsilon,q}^{T} \norm{x - \tilde x}_T . \label{eq:qlip}
\end{align}
Similarly, we see that
\begin{align*}
\norm{v - \tilde v}_T =&\ \norm{\Psi_{y}\big( \alpha(x,q) \big) - \Psi_{y}\big( \alpha(\tilde x, \tilde q) \big)}_T \\
\leq&\ c_{y}^{T}\norm{x - \tilde x}_T +
c_{y}^{T} \bar{\mu} J  \int_{0}^{T} \max_{i \in \III}\abs{q_i(s) - \tilde q_i(s)} ds\notag  \\
\leq&\ c_{y}^{T}\norm{x - \tilde x}_T +
c_{y}^{T} \bar{\mu} J T  \norm{q - \tilde q}_T,
\end{align*}
and by \eqref{eq:qlip}, there exists a constant $c_{\Upsilon, v}^{T} > 0$ satisfying
\begin{align*}
\norm{v - \tilde v}_T \leq c_{\Upsilon,v}^{T} \norm{x - \tilde x}_T .
\end{align*}
This establishes \eqref{eq:ups_lipschitz} and concludes the proof the lemma.
\endproof

\section{Proof of Theorem~\ref{thm:stability}}
\label{app:stability}
This section is devoted to proving Theorem~\ref{thm:stability}. For the remainder of this section, we fix an initial condition $\big(e(0), f(0)\big) \in \mathcal{T}$ and let $\big(e(t),f(t),u(t)\big)$ be the unique solution to the fluid model with this initial condition. Furthermore,  we fix $(\bar e, \bar f) \in \mathcal{E}$ and let $\bar a$ and $\bar m$ be defined by \eqref{eq:one} and \eqref{eq:mdef}, respectively.  For $(x,y) \in \mathcal{T}$, define the function $V: \mathcal{T} \to \R_+$ by
\begin{align}
V(x,y) =&\ \sum_{i= 1}^{r}\sum_{j= 1}^{r} \abs{ y_{ij} - \bar f_{ij}} + \sum_{i= 1}^{r}\sum_{\substack{j = 1\\ i \neq j}}^{r} \abs{x_{ij} - \bar e_{ij}}+ \sum_{i : \bar a_{i} < 1}^{r} x_{ii}+ \Big| \bar m - \sum_{i : \bar a_{i} = 1}^{r} x_{ii}\Big|. \label{eq:lyap}
\end{align}
To prove Theorem~\ref{thm:stability}, we will show that $V(e(t), f(t))$ is a Lyapunov function. We know that $V\big(e(t), f(t)\big)$ is a Lipschitz-continuous function from $\R_+ \to \R_+$ because $V(\cdot), e(\cdot)$, and $f(\cdot)$ are all Lipschitz-continuous.

We say $t > 0$ is a \textit{regular point} of $V\big(e(t), f(t)\big)$ if for all $i,j = 1, \ldots, r$, the functions $e_{ii}(t)$, $\abs{f_{ij}(t) - \bar f_{ij}}$, $\abs{e_{ij}(t) - \bar e_{ij}}$ for $i \neq j$, and $\big| \bar m - \sum_{i : \bar a_{i} = 1}^{r} e_{ii}(t)\big|$ are differentiable at $t$. Since these functions are all Lipschitz-continuous, then almost every point is a regular point. Furthermore, if $t$ is a regular point, we claim that
\begin{align}
f_{ij}(t) = \bar f_{ij} \quad &\Rightarrow \quad \dot f_{ij}(t) = 0, \quad 1 \leq i,j \leq r \label{eq:reg_f}\\
e_{ii}(t) = 0 \quad &\Rightarrow \quad \dot e_{ii}(t) = 0, \quad 1 \leq i \leq r, \label{eq:reg_eii}\\
e_{ij}(t) = \bar e_{ij} \quad &\Rightarrow \quad \dot e_{ij}(t) = 0, \quad 1 \leq i \neq j \leq r \label{eq:reg_eij}\\
\sum_{i : \bar a_{i} = 1}^{r} e_{ii} = \bar m  \quad &\Rightarrow \quad \sum_{i : \bar a_{i} = 1}^{r} \dot e_{ii} = 0, \label{eq:reg_m} \\
\dot u_i(t) \text{ exists }&\text{for all $1 \leq i \leq r$}. \label{eq:reg_u}
\end{align}
Most of these properties are follow directly from the definition of a regular point. For instance, the existence of $\frac{d\abs{f_{ij}(t) - \bar f_{ij}}}{dt}$ together with $f_{ij}(t) = \bar f_{ij}$ implies $\dot f_{ij}(t) = 0$; the same argument applies to \eqref{eq:reg_eii}--\eqref{eq:reg_m}. We know \eqref{eq:reg_u} is true because $\dot u_i(t)$ must exist in order for $\dot e_{ii}(t)$ to exist; cf. \eqref{eq:ufluid_eii}.
\begin{lemma} \label{lem:strict_ineq}
Assume $P_{ij} > 0$ and $Q_{ii} > 0$ for all $i,j= 1, \ldots, r$. If $t > 0$ is a regular point of $V\big(e(t), f(t)\big)$, then $\dot V\big(e(t),f(t)\big) \leq 0$. Furthermore, if $\big(e(t), f(t)\big) \notin \mathcal{E}$, then $\dot V\big(e(t),f(t)\big) < 0$.
\end{lemma}
We postpone the proof of Lemma~\ref{lem:strict_ineq} to the end of this section, and first demonstrate how it is used to prove Theorem~\ref{thm:stability} via LaSalle's Invariance Principle \cite[Theorem 4.4]{Khal2002}. To do so, we need to introduce a few definitions. A point $p \in \mathcal{T}$ is said to be a positive limit point of $\big(e(t), f(t)\big)$ if there exists a sequence $\{t_n\}_{n=1}^{\infty}$ with $t_n \to \infty$ such that $\big(e(t_n), f(t_n)\big) \to p$ as $n \to \infty$. A set $B\subset \mathcal{T}$ is said to be positively invariant if for any $s \geq 0$,
\begin{align*}
\big(e(s), f(s)\big) \in B \quad \Rightarrow \quad \big(e(t), f(t)\big) \in B, \quad \forall\ t \geq s.
\end{align*}
The following lemma is a version of \cite[Lemma 4.1]{Khal2002} adapted to the setting of this paper. We prove it in Section~\ref{sec:khallemma}.
\begin{lemma} \label{lem:khallemma}
 Let $L^{+}$ be the set of all positive limit points of $\big(e(t), f(t)\big)$. Then $L^{+}$ is a nonempty, compact, and positively invariant set. Moreover, 
 \begin{align*}
 \lim_{t \to \infty} \inf_{x \in L^{+}} \norm{\big(e(t), f(t)\big)-x } = 0
 \end{align*}
\end{lemma}
We are now ready to prove Theorem~\ref{thm:stability}. The following argument is repeated from the proof of \cite[Theorem 4.4]{Khal2002}.
\proof{Proof of Theorem~\ref{thm:stability}.}
We know that $V\big(e(t), f(t)\big)$ is a bounded function from $\R_+ \to \R_+$, because $V(\cdot)$ is continuous and $\big(e(t), f(t)\big)$  belong to the compact set $\mathcal{T}$ for all $t \geq 0$. Furthermore, Lemma~\ref{lem:strict_ineq} implies that it is a non-increasing function. Therefore, $\lim_{t \to \infty}  V\big(e(t), f(t)\big)$ exists, and we denote it by $\ell$. Recall by Lemma~\ref{lem:khallemma} that $L^+$, the set of positive limit points of $\big(e(t), f(t)\big)$, is not empty. For any point $p \in L^{+}$, there exists a sequence $\{t_n\}_{n=1}^{\infty}$ such that $\big(e(t_n), f(t_n)\big) \to p$. By continuity of $V(\cdot)$, 
\begin{align*}
V(p) = \lim_{n\to \infty} V\big(e(t_n), f(t_n)\big) = \ell, \quad \forall p \in L^{+}.
\end{align*}
Now suppose $\big(\tilde e(t), \tilde f(t)\big)$ is a solution to the fluid model with initial condition $\big(\tilde e(0), \tilde f(0)\big) \in L^{+}$. Since $L^{+}$ is positively invariant, $\tilde e(t), \tilde f(t) \in L^{+}$ for all $t \geq 0$, which implies that $V\big(\tilde e(t), \tilde f(t)\big) = \ell$, or 
\begin{align*}
\dot V\big(\tilde e(t), \tilde f(t)\big) = 0.
\end{align*}
Using Lemma~\ref{lem:strict_ineq}, we conclude that $L^{+} \subset \mathcal{E}$, which proves Theorem~\ref{thm:stability}.
\endproof
The rest of this section is devoted to proving Lemma~\ref{lem:strict_ineq}.  Fix a regular point $t > 0$. For notational simplicity, we omit the time index $t$ when referring to $V\big(e(t), f(t)\big)$, $e_{ij}(t)$, $f_{ij}(t)$, $u_i(t)$, or their derivatives. From \eqref{eq:ufluid_fij}--\eqref{eq:ufluid_eii} we can see that
\begin{align}
\sum_{i=1}^{r} \sum_{j=1}^{r} \dot e_{ij} + \sum_{i=1}^{r} \sum_{j=1}^{r} \dot f_{ij} = 0. \label{eq:flow_cons}
\end{align}
Set
\begin{align*}
\hat f_{ij} =&\ f_{ij} - \bar f_{ij}, \quad 1 \leq i,j \leq r,\\
\hat e_{ij} =&\ e_{ij} - \bar e_{ij}, \quad  1 \leq i \neq j \leq r.
\end{align*}
Recall that $t$ is a regular point, meaning \eqref{eq:reg_f}--\eqref{eq:reg_u} hold. Therefore,
\begin{align}
\frac{1}{2}\dot V(e,f) =&\ \frac{1}{2}\sum_{i : \bar a_i < 1}^{r} \dot e_{ii}+ \frac{1}{2}\sum_{i=1}^{r} \sum_{\substack{j = 1\\ j \neq i}}^{r} \dot e_{ij} 1\big( e_{ij} > \bar e_{ij}\big) - \frac{1}{2}\sum_{i=1}^{r} \sum_{\substack{j = 1\\ j \neq i}}^{r} \dot e_{ij} 1\big( e_{ij} < \bar e_{ij}\big) \notag \\
&+\frac{1}{2}\sum_{i=1}^{r} \sum_{j=1}^{r} \dot f_{ij} 1\big(f_{ij} > \bar f_{ij}\big)- \frac{1}{2}\sum_{i=1}^{r} \sum_{j=1}^{r} \dot f_{ij} 1\big(f_{ij} < \bar f_{ij}\big) \notag  \\
&+ 1\Big( \sum_{i : \bar a_i = 1} e_{ii} > \bar m\Big)\frac{1}{2}\sum_{i : \bar a_i = 1} \dot e_{ii} - 1\Big( \sum_{i : \bar a_i = 1} e_{ii} < \bar m\Big)\frac{1}{2}\sum_{i : \bar a_i = 1} \dot e_{ii} \notag \\
=&\ \sum_{i : \bar a_i < 1}^{r} \dot e_{ii}+ \sum_{i=1}^{r} \sum_{\substack{j = 1\\ j \neq i}}^{r} \dot e_{ij} 1\big( e_{ij} > \bar e_{ij}\big) +\sum_{i=1}^{r} \sum_{j=1}^{r} \dot f_{ij} 1\big(f_{ij} > \bar f_{ij}\big) \notag  \\
&+ 1\Big( \sum_{i : \bar a_i = 1} e_{ii} > \bar m\Big)\sum_{i : \bar a_i = 1} \dot e_{ii}, \label{eq:halfV}
\end{align}
where the second equality follows from \eqref{eq:flow_cons}. The expression in \eqref{eq:halfV} will soon become very bulky to work with, so before moving forward we first present an illustrative example to help the reader gain some intuition.
\subsection{An Illustrative Example}
Suppose $r = 3$ and that $(\lambda, \mu, P, Q)$ are such that $\bar a_1 < 1$, $\bar a_2 = \bar a_3 = 1$, and $\bar m > 0$. For simplicity, we also assume that all entries of $P$ and $Q$ are strictly positive. We now compute $\frac{1}{2} \dot V(e,f)$ for several choices of $(e, f)\notin \mathcal{E} $, and show that it is always strictly negative. The list of cases we consider is by no means exhaustive, but is nevertheless helpful to develop intuition.

\textbf{Case 1:} Suppose $e_{ii} > 0$ for $i = 1,2,3$ and $e_{22} + e_{33} > \bar m$. Furthermore, suppose $\hat f_{ij} < 0$ for all $i,j,$ and $\hat e_{ij} < 0$ for all $i \neq j$. Using \eqref{eq:halfV} and \eqref{eq:eii_dot}, we see that
\begin{align*}
\frac{1}{2} \dot V(e,f) = \dot e_{11} + \dot e_{22} + \dot e_{33}=\sum_{i=1}^{3} \bigg[-\lambda_i \big( 1 - \dot u_i\big) + \sum_{\substack{j = 1 \\ j \neq i}}^{3} \mu_{ji} e_{ji} + Q_{ii} \sum_{j = 1}^{3} \mu_{ji} f_{ji}\bigg].
\end{align*}
By \eqref{eq:ucomplement}, $\dot u_i = 0$ for $i=1,2,3$ because $e_{ii} > 0$. Furthermore, using \eqref{eq:three} we see that
\begin{align*}
\sum_{\substack{j = 1 \\ j \neq i}}^{3} \mu_{ji} e_{ji} + Q_{ii} \sum_{j = 1}^{3} \mu_{ji} f_{ji} =&\ \sum_{\substack{j = 1 \\ j \neq i}}^{3} \mu_{ji} \bar e_{ji} + Q_{ii} \sum_{j = 1}^{3} \mu_{ji} \bar f_{ji} + \sum_{\substack{j = 1 \\ j \neq i}}^{3} \mu_{ji} \hat e_{ji} + Q_{ii} \sum_{j = 1}^{3} \mu_{ji} \hat f_{ji} \\
=&\ \lambda_i \bar a_i + \sum_{\substack{j = 1 \\ j \neq i}}^{3} \mu_{ji} \hat e_{ji} + Q_{ii} \sum_{j = 1}^{3} \mu_{ji} \hat f_{ji}, \quad 1 \leq i \leq 3.
\end{align*}
Recalling that $\bar a_2 = \bar a_3 = 1$, we arrive at
\begin{align*}
\frac{1}{2} \dot V(e,f) =&\ \lambda_1 (1- \bar a_1)
+ \sum_{i=1}^{3} \bigg[\sum_{\substack{j = 1 \\ j \neq i}}^{3} \mu_{ji} \hat e_{ji} + Q_{ii} \sum_{j = 1}^{3} \mu_{ji} \hat f_{ji}\bigg] < 0.
\end{align*}

\textbf{Case 2:} Suppose $e_{11} = e_{22} = 0$, $e_{33} >  \bar m$, $\hat f_{ij} = 0$ for all $i,j$, $\hat e_{12} < 0$ and $\hat e_{ij} = 0$ for all other $i,j$ with $i \neq j$. In such a case, \eqref{eq:halfV} and \eqref{eq:eii_dot} tell us that
\begin{align*}
\frac{1}{2} \dot V(e,f) = \dot e_{33}= -\lambda_3(1-\dot u_3) + \sum_{\substack{j = 1 \\ j \neq 3}}^{3} \mu_{j3} e_{j3} + Q_{33} \sum_{j = 1}^{3} \mu_{j3} f_{j3}.
\end{align*}
Now $e_{33} > 0$ and \eqref{eq:ucomplement} implies that $\dot u_3 = 0$. Furthermore, we know $\hat e_{j3} = 0$ for $j=1,2$ and $\hat f_{j3} = 0$ for $j = 1,2,3$. Therefore, we can use \eqref{eq:three}
to see that
\begin{align*}
\frac{1}{2} \dot V(e,f)= -\lambda_3 + \sum_{\substack{j = 1 \\ j \neq 3}}^{3} \mu_{j3} \bar e_{j3} + Q_{33} \sum_{j = 1}^{3} \mu_{j3} \bar f_{j3} = 0,
\end{align*}
which appears to contradict what we set out to prove. However, it turns out that the time $t$ corresponding to this configuration of $(e,f)$ is not a regular point. If $t$ were a regular point, then by \eqref{eq:reg_f} we would have $\dot f_{2k} = 0$ for all $k = 1,2,3$, because $\hat f_{2k} = 0$. However, by \eqref{eq:fij_dot}, the fact that $\hat f_{2k} = 0$, and \eqref{eq:one} we see that
\begin{align*}
\dot f_{2k} = \lambda_2 P_{2k} (1-\dot u_2) - \mu_{2k} f_{2k} = \lambda_2 (1-\dot u_2) - \lambda_2 P_{2k} \bar a_2.
\end{align*}
Since $e_{22} = 0$, \eqref{eq:reg_eii} forces $\dot e_{22} = 0$, which together with \eqref{eq:eii_dot} implies that
\begin{align*}
\lambda_2 P_{2k}(1-\dot u_2) = \sum_{\substack{j = 1 \\ j \neq 2}}^{3} \mu_{j2} e_{j2} + Q_{22} \sum_{j = 1}^{3} \mu_{j2} f_{j2} < \sum_{\substack{j = 1 \\ j \neq 2}}^{3} \mu_{j2} \bar e_{j2} + Q_{22} \sum_{j = 1}^{3} \mu_{j2} \bar f_{j2} = \lambda_2 P_{2k} \bar a_2,
\end{align*}
where in the inequality we used that $\hat f_{j2} = \hat e_{32} = 0$, and $\hat e_{12} < 0$, and in the last equality we used \eqref{eq:three}. Therefore, we just showed that $\dot f_{2k} < 0$, which is a contradiction to the assumption that $t$ is a regular point.

Having worked through the example, we now present a general algebraic expansion of $\frac{1}{2} \dot V(e,f)$ in Lemma~\ref{lem:algebra}, which is proved in Section~\ref{sec:algebra}. A line by line inspection of \eqref{eq:vdot_case1} and \eqref{eq:vdot_case2} confirms that $\dot V(e,f) \leq 0$, but proving $\dot V(e,f)<0$ requires more careful arguments presented in Section~\ref{sec:cases}.
\begin{lemma} \label{lem:algebra}
 If $\{i :\ \bar a_{i} = 1\} \neq \emptyset$ and $\sum_{i : \bar a_{i} = 1}^{r} e_{ii} \leq \bar m$, then
\begin{align}
\frac{1}{2}\dot V(e,f) =&\ \sum_{i : \bar a_i < 1}^{r} \lambda_i(\bar a_i - 1) \Big(1-\sum_{j=1}^{r} P_{ij} 1(\hat f_{ij} > 0) \Big)1(\dot u_i = 0) \notag \\
&- \sum_{i= 1}^{r} \Big(1-\sum_{j=1}^{r} P_{ij} 1(\hat f_{ij} > 0) \Big) \Big(\sum_{\substack{j = 1 \\ j \neq i}}^{r} \mu_{ji} \hat e_{ji} 1(\hat e_{ji} \geq 0) + Q_{ii} \sum_{j = 1}^{r} \mu_{ji} \hat f_{ji} 1(\hat f_{ji} \geq 0)\Big) 1(\dot u_i > 0) \notag \\
&+ \sum_{i=1}^{r} \Big(\sum_{j=1}^{r} P_{ij} 1(\hat f_{ij} > 0) \Big) \Big(\sum_{\substack{j = 1 \\ j \neq i}}^{r} \mu_{ji} \hat e_{ji} 1(\hat e_{ji} \leq 0) + Q_{ii} \sum_{j = 1}^{r} \mu_{ji} \hat f_{ji} 1(\hat f_{ji} \leq 0) \Big) 1(\dot u_i > 0)\notag \\
&+ \sum_{i : \bar a_i < 1}^{r} 1( \dot u_i = 0) \bigg[ \sum_{\substack{j = 1 \\ j \neq i}}^{r} \mu_{ji} \hat e_{ji} 1(\hat e_{ji} \leq 0) + Q_{ii} \sum_{j = 1}^{r} \mu_{ji} \hat f_{ji} 1(\hat f_{ji} \leq 0)\bigg]\notag \\
&+\sum_{i= 1}^{r} \sum_{\substack{j = 1\\ j \neq i}}^{r} Q_{ij}1\big( \hat e_{ij} > 0\big) \sum_{k = 1}^{r} \mu_{ki} \hat f_{ki}1 (\hat f_{ki} \leq 0)\notag  \\
&- \sum_{i : \bar a_i < 1}^{r}\Big(1 - Q_{ii} - \sum_{\substack{j = 1\\ j \neq i}}^{r} Q_{ij}1\big( \hat e_{ij} > 0\big)\Big) \sum_{j=1}^{r} \mu_{ji} \hat f_{ji} 1(\hat f_{ji} \geq 0)\notag \\
&-\sum_{i: \bar a_i =1}^{r}  1(\dot u_i = 0)\sum_{\substack{j = 1\\ j \neq i}}^{r}  \mu_{ji} \hat e_{ji} 1\big( \hat e_{ji} \geq 0\big) \notag \\
&- \sum_{i: \bar a_i =1}^{r}\Big(1 - Q_{ii} 1(\dot u_i > 0)- \sum_{\substack{j = 1\\ j \neq i}}^{r} Q_{ij}1( \hat e_{ij} > 0)\Big) \sum_{j=1}^{r} \mu_{ji} \hat f_{ji} 1(\hat f_{ji} \geq  0).
 \label{eq:vdot_case1}
\end{align}
and if $\{i :\ \bar a_{i} = 1\} = \emptyset$, or $\{i :\ \bar a_{i} = 1\} \neq \emptyset$ and $\sum_{i : \bar a_{i} = 1}^{r} e_{ii} > \bar m$, then
\begin{align}
\frac{1}{2}\dot V(e,f) =&\ \sum_{i : \bar a_i < 1}^{r} \lambda_i(\bar a_i - 1) \Big(1-\sum_{j=1}^{r} P_{ij} 1(\hat f_{ij} > 0) \Big)1(\dot u_i = 0) \notag \\
&- \sum_{i= 1}^{r} \Big(1-\sum_{j=1}^{r} P_{ij} 1(\hat f_{ij} > 0) \Big) \Big(\sum_{\substack{j = 1 \\ j \neq i}}^{r} \mu_{ji} \hat e_{ji} 1(\hat e_{ji} \geq 0) + Q_{ii} \sum_{j = 1}^{r} \mu_{ji} \hat f_{ji} 1(\hat f_{ji} \geq 0) \Big) 1(\dot u_i > 0) \notag \\
&+ \sum_{i=1}^{r} \Big(\sum_{j=1}^{r} P_{ij} 1(\hat f_{ij} > 0) \Big) \Big(\sum_{\substack{j = 1 \\ j \neq i}}^{r} \mu_{ji} \hat e_{ji} 1(\hat e_{ji} \leq 0) + Q_{ii} \sum_{j = 1}^{r} \mu_{ji} \hat f_{ji} 1(\hat f_{ji} \leq 0) \Big) 1(\dot u_i > 0)\notag \\
&+ \sum_{i = 1}^{r} 1( \dot u_i = 0) \bigg[ \sum_{\substack{j = 1 \\ j \neq i}}^{r} \mu_{ji} \hat e_{ji} 1(\hat e_{ji} \leq 0) + Q_{ii} \sum_{j = 1}^{r} \mu_{ji} \hat f_{ji} 1(\hat f_{ji} \leq 0)\bigg]\notag \\
&+\sum_{i= 1}^{r} \sum_{\substack{j = 1\\ j \neq i}}^{r} Q_{ij}1\big( \hat e_{ij} > 0\big) \sum_{k = 1}^{r} \mu_{ki} \hat f_{ki}1 (\hat f_{ki} \leq 0)\notag  \\
&- \sum_{i = 1}^{r}\Big(1 - Q_{ii} - \sum_{\substack{j = 1\\ j \neq i}}^{r} Q_{ij}1\big( \hat e_{ij} > 0\big)\Big) \sum_{j=1}^{r} \mu_{ji} \hat f_{ji} 1(\hat f_{ji} \geq 0). \label{eq:vdot_case2}
\end{align}
\end{lemma}

\subsection{Proof of Lemma~\ref{lem:strict_ineq}}
\label{sec:cases}
\proof{Proof of Lemma~\ref{lem:strict_ineq}.}
For notational simplicity, we omit the argument $t$ when referring to $V\big(e(\cdot), f(\cdot)\big), e(\cdot), f(\cdot), u(\cdot)$ and their derivatives.\ We first work in the case when $\{i : \bar a_{i} = 1\} \neq \emptyset$ and $\sum_{i : \bar a_{i} = 1}^{r} e_{ii} \leq \bar m$, meaning that $\frac{1}{2}\dot V(e,f)$ is given by \eqref{eq:vdot_case1}. The following is a list of conditions that are necessary for $\dot V(e,f) = 0$, and are obtained by equating each line in \eqref{eq:vdot_case1} to zero. These conditions are necessary because each line in \eqref{eq:vdot_case1} is non-positive. Any condition on $\hat e_{ij}$ for some $i \neq j$ assumes $Q_{ij} > 0$, because if $Q_{ij} = 0$ then $e_{ij}$ always equals zero.
\begin{enumerate}
\item \label{cond:1} Consider $i$ such that $\bar a_{i} < 1$ and $\dot u_i = 0$. Setting line 1 of \eqref{eq:vdot_case1} to zero requires $\hat f_{ij} > 0$ for all $i,j$. Line 4 requires $\hat f_{ji} \geq 0$ for all $i,j$ and $\hat e_{ji} \geq 0$ for all $i \neq j$. Line 6 requires $\hat e_{ij} > 0$ for all $i \neq j$. 
\item \label{cond:2} Consider $i$ such that $\bar a_i < 1$ and $\dot u_i > 0$. One of the following three mutually exclusive sets of conditions must hold:
\begin{enumerate}[a.]
\item \label{cond:2a} To set line 2 to zero, we choose to enforce $\hat f_{ij} > 0$ for all $i,j$. Line 3 then requires $\hat f_{ji} \geq 0$ for all $i,j$ and $\hat e_{ji} \geq 0$ for all $i \neq j$. Line 6 requires $\hat e_{ij} > 0$ for all $i \neq j$. 
\item \label{cond:2b} This time, to set line 2 to zero, we choose to enforce $\hat f_{ji} \leq 0$ for all $i,j$, and $\hat e_{ji} \leq 0$ for $i \neq j$. Line 2 then requires $\hat f_{ij} \leq 0$ for all $i,j$. Line 5 requires that if $\hat f_{ki} < 0$ for some $k$, then $\hat e_{ij} \leq 0$ for all $j \neq i$.
\item \label{cond:2c} Suppose $f_{ij} > 0$ and $f_{ik} \leq 0$ for some $j, k = 1, \ldots, r$. The only way to make both lines 2 and 3 equal zero is to enforce $\hat f_{ji} = 0$ for all $i,j$ and $\hat e_{ji} = 0$ for $i \neq j$. 
\end{enumerate}
\item \label{cond:3} Consider $i$ such that $\bar a_i = 1$ and $\dot u_i = 0$. Since we assumed $Q_{ii} > 0$, setting line 8 to zero requires $\hat f_{ji} \leq 0$ for all $i,j$. Line 7 requires $\hat e_{ji} \leq 0$ for $i \neq j$, and line 5 requires that if $\hat f_{ki} < 0$ for some $k$, then $\hat e_{ij} \leq 0$ for all $j \neq i$.
\item \label{cond:4} Consider $i$ such that $\bar a_i  = 1$ and $\dot u_i > 0$. One of the following three mutually exclusive sets of conditions must hold:
\begin{enumerate}[a.]
\item \label{cond:4a} To set line 2 to zero, we choose to enforce $\hat f_{ij} > 0$ for all $i,j$. Line 3 then requires $\hat f_{ji} \geq 0$ for all $i,j$ and $\hat e_{ji} \geq 0$ for all $i \neq j$. Line 8 requires that $\hat e_{ij} > 0$ for all $i \neq j$. 
\item \label{cond:4b} This time, to set line 2 to zero, we choose to enforce $\hat f_{ji} \leq 0$ for all $i,j$, and $\hat e_{ji} \leq 0$ for $i \neq j$. Line 2 then requires $\hat f_{ij} \leq 0$ for all $i,j$. Line 5 requires that if $\hat f_{ki} < 0$ for some $k$, then $\hat e_{ij} \leq 0$ for all $j \neq i$.
\item  \label{cond:4c} Suppose $f_{ij} > 0$ and $f_{ik} \leq 0$ for some $j, k = 1, \ldots, r$. The only way to make both lines 2 and 3 equal zero is to enforce $\hat f_{ji} = 0$ for all $i,j$ and $\hat e_{ji} = 0$ for $i \neq j$. 
\end{enumerate}
\end{enumerate}
We now argue that there does not exist a configuration of $e_{ij}$'s and $f_{ij}$'s that satisfies conditions \ref{cond:1}--\ref{cond:4}. If  region $i$ satisfies condition \ref{cond:1}, we refer to it as a type \ref{cond:1} region. This convention is adopted for all other conditions as well.  Recall our assumption that $P_{ij} > 0$ and $Q_{ii} > 0$ for all $i,j$. While we conjecture this lemma to hold even without these conditions, we impose them to prevent the following arguments from becoming even more involved. 

  First observe that a region $i$ can never be of type \ref{cond:4a} or \ref{cond:4c}. For such a region, the fact that $\dot u_i > 0$, together with  \eqref{eq:ucomplement} and \eqref{eq:reg_eii} will imply that $\dot e_{ii}  = 0$, or
\begin{align*}
\lambda_i (1-\dot u_i) = \sum_{\substack{j=1\\j \neq i}}^{r}  \mu_{ji} e_{ji} + Q_{ii} \sum_{j=1}^{r} \mu_{ji} f_{ji}.
\end{align*}
However, the conditions in \ref{cond:4a} and \ref{cond:4c} also imply that
\begin{align*}
\sum_{\substack{j=1\\j \neq i}}^{r}  \mu_{ji} e_{ji} + Q_{ii} \sum_{j=1}^{r} \mu_{ji} f_{ji} \geq \sum_{\substack{j=1\\j \neq i}}^{r}  \mu_{ji} \bar e_{ji} + Q_{ii} \sum_{j=1}^{r} \mu_{ji} \bar f_{ji} = \lambda_i,
\end{align*}
where the equality above follows from \eqref{eq:three} and the fact that $\bar a_i = 1$. This leads to a contradiction because conditions \ref{cond:4a} and \ref{cond:4c} require that $\dot u_i > 0$.

Second, by our assumption that $\{i : \bar a_{i} = 1\} \neq \emptyset$, there must always be a region of either type \ref{cond:3} or type \ref{cond:4b}. This implies that there are no type \ref{cond:1} or type \ref{cond:2a} regions. If $i$ were a type \ref{cond:1} or \ref{cond:2a} region and $j$ were a type \ref{cond:3} or \ref{cond:4b} region, then the conditions in \ref{cond:1} and \ref{cond:2a} would require $\hat f_{ij} > 0$, but the conditions in \ref{cond:3} and \ref{cond:4b} would require that $\hat f_{ij} \leq 0$, which is a contradiction.

Third, we argue that there cannot be a type \ref{cond:2c} region. Suppose $i$ is a type \ref{cond:2c} region. Then for some region $j$, we would have $\hat f_{ij} > 0$. However, this region $j$ could not belong to any of types \ref{cond:2b}, \ref{cond:2c}, \ref{cond:3}, or \ref{cond:4b}, causing a contradiction.

Lastly, we show why it cannot be that $\dot V(e,f) = 0$. We have shown that all regions must be of types \ref{cond:2b}, \ref{cond:3}, or \ref{cond:4b}, which means that $\hat f_{ji} \leq 0$ for all $i,j$ and $\hat e_{ji} \leq 0$ for all $j\neq i$. We also assumed that $\sum_{i : \bar a_{i} = 1}^{r} e_{ii} \leq \bar m$. Since we assumed $(e,f) \notin \mathcal{E}$, at least one of the previous inequalities must be strict. This implies that
\begin{align*}
\sum_{i=1}^{r}\sum_{j=1}^{r} f_{ji} + \sum_{i=1}^{r}\sum_{\substack{j = 1 \\ j \neq i}}^{r} e_{ji} + \sum_{i=1}^{r} e_{ii} =&\ \sum_{i=1}^{r}\sum_{j=1}^{r} f_{ji} + \sum_{i=1}^{r}\sum_{\substack{j = 1 \\ j \neq i}}^{r} e_{ji} + \sum_{i: \bar a_i = 1}^{r} e_{ii} \\
<& \sum_{i=1}^{r}\sum_{j=1}^{r} \bar f_{ji} + \sum_{i=1}^{r}\sum_{\substack{j = 1 \\ j \neq i}}^{r} \bar e_{ji} + \bar m = 1,
\end{align*}
where in the first equality we used the fact that any region $i$ with $\bar a_i < 1$ is of type \ref{cond:2b} and must satisfy $\dot u_i > 0$, which together with by \eqref{eq:ucomplement} implies $e_{ii} = 0$.  The result above implies that the total mass in the system is strictly less than one, which is impossible because the total mass in the system must always equal one. Hence, we have just shown that at any regular point, $\dot V(e,f) < 0$ in the case when $\{i : \bar a_{i} = 1\} \neq \emptyset$ and $\sum_{i : \bar a_{i} = 1}^{r} e_{ii} \leq \bar m$.

We now assume that $\{i :\ \bar a_{i} = 1\} = \emptyset$, or $\{i :\ \bar a_{i} = 1\} \neq \emptyset$ and $\sum_{i : \bar a_{i} = 1}^{r} e_{ii} > \bar m$. Just as before, we assume that $\dot V(e,f) = 0$ and use \eqref{eq:vdot_case2} to list the necessary conditions required for that to happen. Observe that lines 1--3, and 5--6 of both \eqref{eq:vdot_case1} and \eqref{eq:vdot_case2} are identical. There is a slight difference in line 4 of both equations. In the former, the outside summation is over $i : \bar a_{i} < 1$, whereas in the latter the sum is over all $i$. Therefore, conditions \ref{cond:1}, \ref{cond:2}, and \ref{cond:4} must hold as before, and a region cannot be of type \ref{cond:4a} or \ref{cond:4c}. Only condition \ref{cond:3} changes slightly:
\begin{enumerate}
\item [$3'$] \phantomsection \label{cond:b3} Consider $i$ such that $\bar a_i = 1$ and $ \dot u_i = 0$. Setting line 4 of \eqref{eq:vdot_case2} to zero requires $\hat f_{ji} \geq 0$ for all $i,j$, and $\hat e_{ji} \geq 0$ for all $i \neq j$. Line 6 requires that if $\hat f_{ki} > 0$ for some $k$, then $\hat e_{ij} > 0$ for all $j \neq i$.
\end{enumerate}
Again, any condition on $\hat e_{ij}$ assumes $Q_{ij} > 0$, because otherwise $e_{ij}$ always equals zero.

First, we show that there must be a region of type \ref{cond:2b}, \ref{cond:2c} or \ref{cond:4b}. If that were not the case, then all regions would be of type \ref{cond:1}, \ref{cond:2a} or \hyperref[cond:b3]{$3'$}. Recall our assumption that $(e,f) \notin \mathcal{E}$, $\{i :\ \bar a_{i} = 1\} = \emptyset$, or $\{i :\ \bar a_{i} = 1\} \neq \emptyset$ and $\sum_{i : \bar a_{i} = 1}^{r} e_{ii} > \bar m$. If $\{i :\ \bar a_{i} = 1\} = \emptyset$, then all regions are of type \ref{cond:1} or \ref{cond:2a} and
\begin{align*}
\sum_{i=1}^{r}\sum_{j=1}^{r} f_{ij} + \sum_{i=1}^{r}\sum_{\substack{j = 1 \\ j \neq i}}^{r} e_{ij} + \sum_{i=1}^{r} e_{ii}  =&\ \sum_{i=1}^{r}\sum_{j=1}^{r} f_{ij} + \sum_{i=1}^{r}\sum_{\substack{j = 1 \\ j \neq i}}^{r} e_{ji} + \sum_{i : \bar a_i < 1}^{r} e_{ii} \\
>&\ \sum_{i=1}^{r}\sum_{j=1}^{r} \bar f_{ij} + \sum_{i=1}^{r}\sum_{\substack{j = 1 \\ j \neq i}}^{r} \bar e_{ji} \\
=&\  \sum_{i=1}^{r}\sum_{j=1}^{r} \bar f_{ij} + \sum_{i=1}^{r}\sum_{\substack{j = 1 \\ j \neq i}}^{r} \bar e_{ji} + \sum_{i : \bar a_i < 1}^{r} \bar e_{ii} =  1.
\end{align*}
The inequality holds because if region $i$ is of type \ref{cond:1} or \ref{cond:2a}, then $\hat f_{ij} > 0$ for all $j$ and $\hat e_{ji} \geq 0$ for all $j \neq i$, and the second last equality holds because $\bar e_{ii} = 0$ for $i$ such that $\bar a_i < 1$. The inequality above is a contradiction because the total fluid mass in the system must always equal one. By similar reasoning, if $\{i :\ \bar a_{i} = 1\} \neq \emptyset$ and $\sum_{i : \bar a_{i} = 1}^{r} e_{ii} > \bar m$, then 
\begin{align*}
\sum_{i=1}^{r}\sum_{j=1}^{r} f_{ij} + \sum_{i=1}^{r}\sum_{\substack{j = 1 \\ j \neq i}}^{r} e_{ij} + \sum_{i=1}^{r} e_{ii} > \sum_{i=1}^{r}\sum_{j=1}^{r} \bar f_{ij} + \sum_{i=1}^{r}\sum_{\substack{j = 1 \\ j \neq i}}^{r} \bar e_{ij} + \bar m = 1,
\end{align*}
which is again a contradiction because fluid mass is conserved. Hence, it cannot be that all regions are of type \ref{cond:1}, \ref{cond:2a} or \hyperref[cond:b3]{$3'$}. The fact that there is always a region of type \ref{cond:2b}, \ref{cond:2c} or \ref{cond:4b} implies that there cannot be any type \ref{cond:1} or \ref{cond:2a} regions. If $j$ belonged to the former group and $i$ to the latter, then the definitions of type \ref{cond:1} or \ref{cond:2a} would imply that $\hat f_{ij} > 0$, which would contradict the requirements in types \ref{cond:2b}, \ref{cond:2c}  and \ref{cond:4b}.

Second, we show that there must always be a region of type \hyperref[cond:b3]{$3'$}. If there is a region of type \ref{cond:4b} then it must be that $\{i :\ \bar a_{i} = 1\} \neq \emptyset$ and $\sum_{i : \bar a_{i} = 1}^{r} e_{ii} > \bar m$. However, if $i$ is a type \ref{cond:4b} region, then $\bar e_{ii} = 0$ because $u_i > 0$. Hence, for $\sum_{i : \bar a_{i} = 1}^{r} e_{ii} > \bar m$ to be true, the set $\{i :\ \bar a_{i} = 1\}$ must contain at least one type \hyperref[cond:b3]{$3'$} region. Now suppose there is no region of type \ref{cond:4b}. We argue that it cannot be the case that all regions are exclusively of type \ref{cond:2b} or exclusively of type \ref{cond:2c}. The former case would imply that the total mass in the system is strictly less than one, and the latter case cannot happen by definition of \ref{cond:2c} (i.e.\ if $i$ were of type \ref{cond:2c} then there would be some $j$ with $\hat f_{ij} > 0$, but this $j$ could not be of type \ref{cond:2b} or \ref{cond:2c}). The same reasoning implies that the regions cannot be a mixture of exclusively types \ref{cond:2b} and \ref{cond:2c}. Therefore, there must always exist a region of type \hyperref[cond:b3]{$3'$}.

Lastly, we show why it cannot be that $\dot V(e,f) = 0$. Let $i$ be a type \hyperref[cond:b3]{$3'$} region with $e_{ii} > 0$ (such a region must always exist because $\sum_{i : \bar a_{i} = 1}^{r} e_{ii} > \bar m$). Since fluid mass is conserved, there must exist some $k, \ell$, such that either $\hat f_{k \ell} < 0$, or $k \neq \ell$ and $\hat e_{k \ell} < 0$. Observe that $\ell$ cannot be a type \ref{cond:2c} or \hyperref[cond:b3]{$3'$} region, so it must be a type \ref{cond:2b} or \ref{cond:4b} region. Furthermore, since $i$ is of type \hyperref[cond:b3]{$3'$} and $\ell$ is of type \ref{cond:2b} or \ref{cond:4b}, it must be true that $\hat f_{\ell i} = 0$, and by \eqref{eq:reg_f} this would imply that $\dot f_{\ell i} = 0$. We will now show that $\dot f_{\ell i}$ must also be strictly less than zero, leading to a contradiction. Using \eqref{eq:ufluid_fij}, \eqref{eq:one}, and the fact that $ f_{\ell i} = \bar  f_{\ell i}$, it follows that
\begin{align*}
0 = \dot f_{\ell i} = \lambda_{\ell} P_{\ell i} (1 - \dot u_{\ell}) - \mu_{\ell i} f_{\ell i} = \lambda_{\ell} P_{\ell i} (1 - \dot u_{\ell})- \lambda_{\ell} P_{\ell i} \bar a_{\ell}.
\end{align*}
We know that $\dot u_{\ell} > 0$ by definition of a type \ref{cond:2b} and \ref{cond:4b} region. Hence, $e_{\ell \ell} = 0$ by \eqref{eq:ucomplement}, which in turn means that $\dot e_{\ell\ell} = 0$ by \eqref{eq:reg_eii}, and therefore
\begin{align}
\lambda_{\ell} \big( 1 - \dot u_{\ell}\big) =&\ \sum_{\substack{j = 1 \\ j \neq {\ell}}}^{r} \mu_{j{\ell}} e_{j\ell} + Q_{\ell \ell} \sum_{j = 1}^{r} \mu_{j\ell} f_{j\ell} \notag \\
=&\ \lambda_{\ell} \bar a_{\ell} + \sum_{\substack{j = 1 \\ j \neq \ell}}^{r} \mu_{j\ell} \hat e_{j\ell} + Q_{\ell\ell} \sum_{j = 1}^{r} \mu_{j\ell} \hat f_{j\ell}, \label{eq:dotu}
\end{align}
where the first equality follows from differentiating \eqref{eq:ufluid_eii} and setting the left hand side to zero, and the second equality follows from \eqref{eq:three}. Combining \eqref{eq:dotu} with the form of $\dot f_{\ell i}$, we see that
\begin{align*}
\dot f_{\ell i} =&\ P_{\ell i}\Big(\lambda_{\ell} \bar a_{\ell} +  \sum_{\substack{j = 1\\ j \neq \ell}}^{r} \mu_{j \ell} \hat e_{j \ell}  + Q_{\ell \ell} \sum_{j=1}^{r} \mu_{j \ell} \hat f_{j \ell}\Big) - \lambda_{\ell} P_{\ell i} \bar a_{\ell} = P_{\ell i}\Big(\sum_{\substack{j = 1\\ j \neq \ell}}^{r} \mu_{j \ell} \hat e_{j \ell}  + Q_{\ell \ell} \sum_{j=1}^{r} \mu_{j \ell} \hat f_{j \ell}\Big).
\end{align*}
Since $\ell$ is of type \ref{cond:2b} or \ref{cond:4b}, it follows that $\hat f_{j \ell} \leq 0$ for all $j$ and $\hat e_{j \ell} \leq 0$ for all $j \neq \ell$. Furthermore, we know that either $\hat f_{k \ell} < 0$, or $k \neq \ell$ and $\hat e_{k \ell} < 0$, which implies that $\dot f_{\ell i}  < 0$. However, this is a contradiction because $\dot f_{\ell i} = 0$ at regular points. Therefore, the necessary conditions required for $\dot V(e,f) = 0$ to hold cause a contradiction, and so $\dot V(e,f) < 0$. This concludes the proof of Lemma~\ref{lem:strict_ineq}.
\endproof

\subsection{Proof of Lemma~\ref{lem:algebra}} \label{sec:algebra}
\proof{Proof of Lemma~\ref{lem:algebra}.}
We begin with the case that $\{i :\ \bar a_{i} = 1\} \neq \emptyset$  and $\sum_{i : \bar a_{i} = 1}^{r} e_{ii} \leq \bar m$. Starting with \eqref{eq:halfV}, and taking derivatives in \eqref{eq:ufluid_fij}--\eqref{eq:ufluid_eii}, we know that
\begin{align}
 \frac{1}{2}\dot V(e,f) =&\ \sum_{i : \bar a_i < 1}^{r} \bigg(-\lambda_i \big( 1 - \dot u_i\big) + \sum_{\substack{j = 1 \\ j \neq i}}^{r} \mu_{ji} e_{ji} + Q_{ii} \sum_{j = 1}^{r} \mu_{ji} f_{ji}\bigg) \label{eq:l1} \\
&+ \sum_{i=1}^{r} \sum_{\substack{j = 1\\ j \neq i}}^{r} \bigg( - \mu_{ij} e_{ij} + Q_{ij} \sum_{k = 1}^{r} \mu_{ki} f_{ki} \bigg) 1\big( \hat e_{ij} > 0\big) \label{eq:l2}\\
&+\sum_{i=1}^{r} \sum_{j=1}^{r} \Big(\lambda_{i} P_{ij}\big( 1 - \dot u_i\big) - \mu_{ij} f_{ij}\Big) 1\big(\hat f_{ij} >  0\big), \label{eq:l3}
\end{align}
The first line, \eqref{eq:l1}, equals
\begin{align}
&\sum_{i : \bar a_i < 1}^{r} \bigg(-\lambda_i \big( 1 - \dot u_i\big) + \sum_{\substack{j = 1 \\ j \neq i}}^{r} \mu_{ji} \bar e_{ji} + Q_{ii} \sum_{j = 1}^{r} \mu_{ji} \bar f_{ji}\bigg) + \sum_{i : \bar a_i < 1}^{r} \bigg( \sum_{\substack{j = 1 \\ j \neq i}}^{r} \mu_{ji} \hat e_{ji} + Q_{ii} \sum_{j = 1}^{r} \mu_{ji} \hat f_{ji}\bigg) \notag \\
=&\ \sum_{i : \bar a_i < 1}^{r} \Big(-\lambda_i \big( 1 - \dot u_i\big) +\lambda_i \bar a_i\Big) + \sum_{i : \bar a_i < 1}^{r} \bigg( \sum_{\substack{j = 1 \\ j \neq i}}^{r} \mu_{ji} \hat e_{ji} + Q_{ii} \sum_{j = 1}^{r} \mu_{ji} \hat f_{ji}\bigg) \notag \\
=&\ \sum_{i : \bar a_i < 1}^{r} \lambda_i(\bar a_i - 1) 1(\dot u_i = 0) + \sum_{i : \bar a_i < 1}^{r} \Big(-\lambda_i \big( 1 - \dot u_i\big) +\lambda_i \bar a_i\Big)1(\dot u_i > 0)\notag \\
&+ \sum_{i : \bar a_i < 1}^{r} \bigg( \sum_{\substack{j = 1 \\ j \neq i}}^{r} \mu_{ji} \hat e_{ji} + Q_{ii} \sum_{j = 1}^{r} \mu_{ji} \hat f_{ji}\bigg) \notag \\
=&\ \sum_{i : \bar a_i < 1}^{r} \lambda_i(\bar a_i - 1) 1(\dot u_i = 0) \Big(1-\sum_{j=1}^{r} P_{ij} 1(\hat f_{ij} > 0) \Big) \notag \\
&+ \sum_{i : \bar a_i < 1}^{r} \lambda_i(\bar a_i - 1) 1(\dot u_i = 0)\sum_{j=1}^{r} P_{ij} 1(\hat f_{ij} > 0) \notag \\
&+ \sum_{i : \bar a_i < 1}^{r} \Big(-\lambda_i \big( 1 - \dot u_i\big) +\lambda_i \bar a_i\Big)1(\dot u_i > 0)\Big(1-\sum_{j=1}^{r} P_{ij} 1(\hat f_{ij} > 0) \Big) \notag \\
&+ \sum_{i : \bar a_i < 1}^{r} \Big(-\lambda_i \big( 1 - \dot u_i\big) +\lambda_i \bar a_i\Big)1(\dot u_i > 0)\sum_{j=1}^{r} P_{ij} 1(\hat f_{ij} > 0) \notag \\
&+ \sum_{i : \bar a_i < 1}^{r} \bigg( \sum_{\substack{j = 1 \\ j \neq i}}^{r} \mu_{ji} \hat e_{ji} + Q_{ii} \sum_{j = 1}^{r} \mu_{ji} \hat f_{ji}\bigg), \label{eq:l1.1}
\end{align}
where the first equality comes from \eqref{eq:three}. The second line, \eqref{eq:l2}, equals
\begin{align}
&\sum_{i=1}^{r} \sum_{\substack{j = 1\\ j \neq i}}^{r} \bigg( - \mu_{ij} \bar e_{ij} + Q_{ij} \sum_{k = 1}^{r} \mu_{ki} \bar f_{ki} \bigg) 1\big( \hat e_{ij} > 0\big) + \sum_{i=1}^{r} \sum_{\substack{j = 1\\ j \neq i}}^{r} \bigg( - \mu_{ij} \hat e_{ij} + Q_{ij} \sum_{k = 1}^{r} \mu_{ki} \hat f_{ki} \bigg) 1\big( \hat e_{ij} > 0\big) \notag \\
=&\ -\sum_{i=1}^{r} \sum_{\substack{j = 1\\ j \neq i}}^{r}  \mu_{ji} \hat e_{ji} 1\big( \hat e_{ji} > 0\big) + \sum_{i=1}^{r} \sum_{\substack{j = 1\\ j \neq i}}^{r} Q_{ij}1\big( \hat e_{ij} > 0\big) \sum_{k = 1}^{r} \mu_{ki} \hat f_{ki}  \label{eq:l2.1},
\end{align}
where in the equality we used \eqref{eq:two}. Similarly, we use \eqref{eq:one} to see that the third line, \eqref{eq:l3}, equals
\begin{align}
&\sum_{i=1}^{r} \sum_{j=1}^{r} \Big(\lambda_{i} P_{ij}\big( 1 - \dot u_i\big) - \mu_{ij} \bar f_{ij}\Big) 1\big(\hat f_{ij} >  0\big) - \sum_{i=1}^{r} \sum_{j=1}^{r} \mu_{ij} \hat f_{ij} 1\big(\hat f_{ij} >  0\big) \notag \\
=&\ \sum_{i=1}^{r} 1(\dot u_i = 0)\lambda_i (1 - \bar a_i) \sum_{j=1}^{r} P_{ij} 1\big(\hat f_{ij} >  0\big) \notag \\
&+\sum_{i=1}^{r} 1(\dot u_i > 0) \sum_{j=1}^{r} P_{ij} \Big(\lambda_{i} \big( 1 - \dot u_i\big) - \lambda_i \bar a_i\Big) 1\big(\hat f_{ij} >  0\big)  - \sum_{i=1}^{r} \sum_{j=1}^{r} \mu_{ji} \hat f_{ji} 1\big(\hat f_{ji} >  0\big) \notag \\ 
=&\ \sum_{i=1}^{r} 1(\dot u_i = 0)\lambda_i (1 - \bar a_i) \sum_{j=1}^{r} P_{ij} 1\big(\hat f_{ij} >  0\big)  \notag \\
&+\sum_{i : \bar a_i < 1}^{r} \Big(\lambda_{i} \big( 1 - \dot u_i\big) - \lambda_i \bar a_i\Big)1(\dot u_i > 0) \sum_{j=1}^{r} P_{ij}  1\big(\hat f_{ij} >  0\big) \notag \\
&+\sum_{i : \bar a_i = 1}^{r} \Big(\lambda_{i} \big( 1 - \dot u_i\big) - \lambda_i \bar a_i\Big)1(\dot u_i > 0) \sum_{j=1}^{r} P_{ij}  1\big(\hat f_{ij} >  0\big)  - \sum_{i=1}^{r} \sum_{j=1}^{r} \mu_{ji} \hat f_{ji} 1\big(\hat f_{ji} >  0\big) \label{eq:l3.1}
\end{align}
In addition to \eqref{eq:l1.1}--\eqref{eq:l3.1}, we will require the equation
\begin{align}
\lambda_i(1-\dot u_i) - \lambda_i a_i = \sum_{\substack{j = 1 \\ j \neq i}}^{r} \mu_{ji} \hat e_{ji} + Q_{i i} \sum_{j = 1}^{r} \mu_{ji} \hat f_{ji}, \quad \text{ for } i  \text{ such that } \dot u_i > 0, \label{eq:handy}
\end{align}
which is argued in the same way as \eqref{eq:dotu}. We now combine \eqref{eq:l1.1}--\eqref{eq:handy} to see that
\begin{align}
\frac{1}{2}\dot V(e,f) =&\ \sum_{i : \bar a_i < 1}^{r} \lambda_i(\bar a_i - 1) \Big(1-\sum_{j=1}^{r} P_{ij} 1(\hat f_{ij} > 0) \Big)1(\dot u_i = 0) \label{eq:m1}\\
&- \sum_{i : \bar a_i < 1}^{r} \Big(1-\sum_{j=1}^{r} P_{ij} 1(\hat f_{ij} > 0) \Big) \Big(\sum_{\substack{j = 1 \\ j \neq i}}^{r} \mu_{ji} \hat e_{ji} + Q_{ii} \sum_{j = 1}^{r} \mu_{ji} \hat f_{ji} \Big) 1(\dot u_i > 0) \label{eq:m2} \\
&+ \sum_{i : \bar a_{i} = 1}^{r} \Big(\sum_{k=1}^{r} P_{ik} 1(\hat f_{ik} > 0) \Big)\Big(\sum_{\substack{j = 1 \\ j \neq i}}^{r} \mu_{ji} \hat e_{ji} + Q_{ii} \sum_{j = 1}^{r} \mu_{ji} \hat f_{ji} \Big)1(\dot u_i > 0) \label{eq:m3} \\
&+ \sum_{i : \bar a_i < 1}^{r} \bigg( \sum_{\substack{j = 1 \\ j \neq i}}^{r} \mu_{ji} \hat e_{ji} + Q_{ii} \sum_{j = 1}^{r} \mu_{ji} \hat f_{ji}\bigg) \label{eq:m4} \\
&-\sum_{i=1}^{r} \sum_{\substack{j = 1\\ j \neq i}}^{r}  \mu_{ji} \hat e_{ji} 1\big( \hat e_{ji} > 0\big) + \sum_{i=1}^{r} \sum_{\substack{j = 1\\ j \neq i}}^{r} Q_{ij}1\big( \hat e_{ij} > 0\big) \sum_{k = 1}^{r} \mu_{ki} \hat f_{ki} \label{eq:m5} \\
&- \sum_{i=1}^{r} \sum_{j=1}^{r} \mu_{ji} \hat f_{ji} 1\big(\hat f_{ji} >  0\big). \label{eq:m6}
\end{align}
Since the term above is very bulky, we manipulate one line at a time to help exposition. We leave \eqref{eq:m1} as is, and decompose  \eqref{eq:m2} into
\begin{align*}
&- \sum_{i : \bar a_i < 1}^{r} \Big(1-\sum_{j=1}^{r} P_{ij} 1(\hat f_{ij} > 0) \Big) \Big(\sum_{\substack{j = 1 \\ j \neq i}}^{r} \mu_{ji} \hat e_{ji} 1(\hat e_{ji} \geq 0) + Q_{ii} \sum_{j = 1}^{r} \mu_{ji} \hat f_{ji} 1(\hat f_{ji} \geq 0 \Big) 1(\dot u_i > 0) \\
&- \sum_{i : \bar a_i < 1}^{r} \Big(1-\sum_{j=1}^{r} P_{ij} 1(\hat f_{ij} > 0) \Big) \Big(\sum_{\substack{j = 1 \\ j \neq i}}^{r} \mu_{ji} \hat e_{ji} 1(\hat e_{ji} \leq 0) + Q_{ii} \sum_{j = 1}^{r} \mu_{ji} \hat f_{ji} 1(\hat f_{ji} \leq 0 \Big) 1(\dot u_i > 0).
\end{align*}
We leave \eqref{eq:m3} as is. The term in \eqref{eq:m4} equals
\begin{align*}
&\sum_{i : \bar a_i < 1}^{r} \bigg( \sum_{\substack{j = 1 \\ j \neq i}}^{r} \mu_{ji} \hat e_{ji} 1(\hat e_{ji} \geq 0) + Q_{ii} \sum_{j = 1}^{r} \mu_{ji} \hat f_{ji} 1(\hat f_{ji} \geq 0)\bigg) \\
&+ \sum_{i : \bar a_i < 1}^{r} \bigg( \sum_{\substack{j = 1 \\ j \neq i}}^{r} \mu_{ji} \hat e_{ji} 1(\hat e_{ji} \leq 0) + Q_{ii} \sum_{j = 1}^{r} \mu_{ji} \hat f_{ji} 1(\hat f_{ji} \leq 0)\bigg).
\end{align*}
The term in \eqref{eq:m5} equals
\begin{align*}
&-\sum_{i : \bar a_i < 1}^{r} \sum_{\substack{j = 1\\ j \neq i}}^{r}  \mu_{ji} \hat e_{ji} 1\big( \hat e_{ji} \geq 0\big) + \sum_{i: \bar a_i < 1}^{r} \sum_{\substack{j = 1\\ j \neq i}}^{r} Q_{ij}1\big( \hat e_{ij} > 0\big) \sum_{k = 1}^{r} \mu_{ki} \hat f_{ki}1 (\hat f_{ki} \leq 0) \\
&+ \sum_{i: \bar a_i < 1}^{r} \sum_{\substack{j = 1\\ j \neq i}}^{r} Q_{ij}1\big( \hat e_{ij} > 0\big) \sum_{k = 1}^{r} \mu_{ki} \hat f_{ki}1 (\hat f_{ki} \geq 0)\\
&-\sum_{i: \bar a_i = 1}^{r} \sum_{\substack{j = 1\\ j \neq i}}^{r}  \mu_{ji} \hat e_{ji} 1\big( \hat e_{ji} \geq 0\big) + \sum_{i: \bar a_i = 1}^{r} \sum_{\substack{j = 1\\ j \neq i}}^{r} Q_{ij}1\big( \hat e_{ij} > 0\big) \sum_{k = 1}^{r} \mu_{ki} \hat f_{ki},
\end{align*}
and the term in \eqref{eq:m6} equals
\begin{align*}
- \sum_{i : \bar a_i < 1}^{r} \sum_{j=1}^{r} \mu_{ji} \hat f_{ji} 1\big(\hat f_{ji} \geq  0\big)- \sum_{i : \bar a_i = 1}^{r} \sum_{j=1}^{r} \mu_{ji} \hat f_{ji} 1\big(\hat f_{ji} \geq  0\big).
\end{align*}
Putting all of these expansions back into \eqref{eq:m1}--\eqref{eq:m6}, we see that
\begin{align}
\frac{1}{2}\dot V(e,f) =&\ \sum_{i : \bar a_i < 1}^{r} \lambda_i(\bar a_i - 1) \Big(1-\sum_{j=1}^{r} P_{ij} 1(\hat f_{ij} > 0) \Big)1(\dot u_i = 0) \label{eq:nbegin}\\
&- \sum_{i : \bar a_i < 1}^{r} \Big(1-\sum_{j=1}^{r} P_{ij} 1(\hat f_{ij} > 0) \Big) \Big(\sum_{\substack{j = 1 \\ j \neq i}}^{r} \mu_{ji} \hat e_{ji} 1(\hat e_{ji} \geq 0) + Q_{ii} \sum_{j = 1}^{r} \mu_{ji} \hat f_{ji} 1(\hat f_{ji} \geq 0 \Big) 1(\dot u_i > 0) \\
&+ \sum_{i : \bar a_i < 1}^{r} \Big(\sum_{j=1}^{r} P_{ij} 1(\hat f_{ij} > 0) \Big) \Big(\sum_{\substack{j = 1 \\ j \neq i}}^{r} \mu_{ji} \hat e_{ji} 1(\hat e_{ji} \leq 0) + Q_{ii} \sum_{j = 1}^{r} \mu_{ji} \hat f_{ji} 1(\hat f_{ji} \leq 0) \Big) 1(\dot u_i > 0) \\
&+ \sum_{i : \bar a_i < 1}^{r} 1( \dot u_i = 0) \bigg( \sum_{\substack{j = 1 \\ j \neq i}}^{r} \mu_{ji} \hat e_{ji} 1(\hat e_{ji} \leq 0) + Q_{ii} \sum_{j = 1}^{r} \mu_{ji} \hat f_{ji} 1(\hat f_{ji} \leq 0)\bigg)\\
&+\sum_{i: \bar a_i < 1}^{r} \sum_{\substack{j = 1\\ j \neq i}}^{r} Q_{ij}1\big( \hat e_{ij} > 0\big) \sum_{k = 1}^{r} \mu_{ki} \hat f_{ki}1 (\hat f_{ki} \leq 0)  \\
&- \sum_{i : \bar a_i < 1}^{r}\Big(1 - Q_{ii} - \sum_{\substack{j = 1\\ j \neq i}}^{r} Q_{ij}1\big( \hat e_{ij} > 0\big)\Big) \sum_{j=1}^{r} \mu_{ji} \hat f_{ji} 1(\hat f_{ji} \geq 0) \label{eq:nend} \\
&+ \sum_{i : \bar a_{i} = 1}^{r} \Big(\sum_{k=1}^{r} P_{ik} 1(\hat f_{ik} > 0) \Big)\Big(\sum_{\substack{j = 1 \\ j \neq i}}^{r} \mu_{ji} \hat e_{ji} + Q_{ii} \sum_{j = 1}^{r} \mu_{ji} \hat f_{ji} \Big)1(\dot u_i > 0)  \label{eq:n1}\\
&-\sum_{i: \bar a_i =1}^{r} \sum_{\substack{j = 1\\ j \neq i}}^{r}  \mu_{ji} \hat e_{ji} 1\big( \hat e_{ji} \geq 0\big) + \sum_{i: \bar a_i =1}^{r} \sum_{\substack{j = 1\\ j \neq i}}^{r} Q_{ij}1\big( \hat e_{ij} > 0\big) \sum_{k = 1}^{r} \mu_{ki} \hat f_{ki}\label{eq:n2}\\
&- \sum_{i: \bar a_i =1}^{r} \sum_{j=1}^{r} \mu_{ji} \hat f_{ji} 1\big(\hat f_{ji} \geq  0\big). \label{eq:n3}
\end{align}
It remains to manipulate the terms in \eqref{eq:n1}--\eqref{eq:n3} to get them into the form we need. We begin with \eqref{eq:n1}, which equals
\begin{align*}
&\sum_{i : \bar a_{i} = 1}^{r} \Big(\sum_{k=1}^{r} P_{ik} 1(\hat f_{ik} > 0) \Big)\Big(\sum_{\substack{j = 1 \\ j \neq i}}^{r} \mu_{ji} \hat e_{ji} 1(\hat e_{ji} \leq 0) + Q_{ii} \sum_{j = 1}^{r} \mu_{ji} \hat f_{ji} 1(\hat f_{ji} \leq 0) \Big)1(\dot u_i > 0)\\
&+\sum_{i : \bar a_{i} = 1}^{r} \Big(\sum_{k=1}^{r} P_{ik} 1(\hat f_{ik} > 0) \Big)\Big(\sum_{\substack{j = 1 \\ j \neq i}}^{r} \mu_{ji} \hat e_{ji} 1(\hat e_{ji} \geq 0) + Q_{ii} \sum_{j = 1}^{r} \mu_{ji} \hat f_{ji} 1(\hat f_{ji} \geq 0) \Big)1(\dot u_i > 0).
\end{align*}
The term in \eqref{eq:n2} equals
\begin{align*}
&-\sum_{i: \bar a_i =1}^{r} \sum_{\substack{j = 1\\ j \neq i}}^{r}  \mu_{ji} \hat e_{ji} 1\big( \hat e_{ji} \geq 0\big) + \sum_{i: \bar a_i =1}^{r} \sum_{\substack{j = 1\\ j \neq i}}^{r} Q_{ij}1\big( \hat e_{ij} > 0\big) \sum_{k = 1}^{r} \mu_{ki} \hat f_{ki} 1(\hat f_{ki} \leq 0) \\
&+ \sum_{i: \bar a_i =1}^{r} \sum_{\substack{j = 1\\ j \neq i}}^{r} Q_{ij}1\big( \hat e_{ij} > 0\big) \sum_{k = 1}^{r} \mu_{ki} \hat f_{ki} 1(\hat f_{ki} \geq 0),
\end{align*}
and the term in \eqref{eq:n3} equals
\begin{align*}
&- \sum_{i: \bar a_i =1}^{r}\Big(1 - Q_{ii} 1(\dot u_i > 0)\Big) \sum_{j=1}^{r} \mu_{ji} \hat f_{ji} 1\big(\hat f_{ji} \geq  0\big) - \sum_{i: \bar a_i =1}^{r} Q_{ii} 1(\dot u_i > 0) \sum_{j=1}^{r} \mu_{ji} \hat f_{ji} 1\big(\hat f_{ji} \geq  0\big).
\end{align*}
Inserting these expansions back into \eqref{eq:n1}--\eqref{eq:n3}, we see that 
\begin{align*}
& \sum_{i : \bar a_{i} = 1}^{r} \Big(\sum_{k=1}^{r} P_{ik} 1(\hat f_{ik} > 0) \Big)\Big(\sum_{\substack{j = 1 \\ j \neq i}}^{r} \mu_{ji} \hat e_{ji} + Q_{ii} \sum_{j = 1}^{r} \mu_{ji} \hat f_{ji} \Big)1(\dot u_i > 0) \\
&-\sum_{i: \bar a_i =1}^{r} \sum_{\substack{j = 1\\ j \neq i}}^{r}  \mu_{ji} \hat e_{ji} 1\big( \hat e_{ji} \geq 0\big) + \sum_{i: \bar a_i =1}^{r} \sum_{\substack{j = 1\\ j \neq i}}^{r} Q_{ij}1\big( \hat e_{ij} > 0\big) \sum_{k = 1}^{r} \mu_{ki} \hat f_{ki} \\
&- \sum_{i: \bar a_i =1}^{r} \sum_{j=1}^{r} \mu_{ji} \hat f_{ji} 1\big(\hat f_{ji} \geq  0\big) \\
=&\ \sum_{i : \bar a_{i} = 1}^{r} \Big(\sum_{k=1}^{r} P_{ik} 1(\hat f_{ik} > 0) \Big)\Big(\sum_{\substack{j = 1 \\ j \neq i}}^{r} \mu_{ji} \hat e_{ji} 1(\hat e_{ji} \leq 0) + Q_{ii} \sum_{j = 1}^{r} \mu_{ji} \hat f_{ji} 1(\hat f_{ji} \leq 0) \Big)1(\dot u_i > 0)\\
&-\sum_{i : \bar a_{i} = 1}^{r} \Big(1-\sum_{k=1}^{r} P_{ik} 1(\hat f_{ik} > 0) \Big)\Big(\sum_{\substack{j = 1 \\ j \neq i}}^{r} \mu_{ji} \hat e_{ji} 1(\hat e_{ji} \geq 0) + Q_{ii} \sum_{j = 1}^{r} \mu_{ji} \hat f_{ji} 1(\hat f_{ji} \geq 0) \Big)1(\dot u_i > 0)\\
&-\sum_{i: \bar a_i =1}^{r}  1(\dot u_i = 0)\sum_{\substack{j = 1\\ j \neq i}}^{r}  \mu_{ji} \hat e_{ji} 1\big( \hat e_{ji} \geq 0\big)+ \sum_{i: \bar a_i =1}^{r} \sum_{\substack{j = 1\\ j \neq i}}^{r} Q_{ij}1\big( \hat e_{ij} > 0\big) \sum_{k = 1}^{r} \mu_{ki} \hat f_{ki} 1(\hat f_{ki} \leq 0) \\
&- \sum_{i: \bar a_i =1}^{r}\Big(1 - Q_{ii} 1(\dot u_i > 0)- \sum_{\substack{j = 1\\ j \neq i}}^{r} Q_{ij}1( \hat e_{ij} > 0)\Big) \sum_{j=1}^{r} \mu_{ji} \hat f_{ji} 1(\hat f_{ji} \geq  0).
\end{align*}
The form of $\frac{1}{2} \dot V(e,f)$ we obtain by plugging the above back into \eqref{eq:nbegin}--\eqref{eq:n3} can be compared with \eqref{eq:vdot_case1} to see that it matches.

Now suppose that $\{i :\ \bar a_{i} = 1\} = \emptyset$, or $\{i :\ \bar a_{i} = 1\} \neq \emptyset$ and $\sum_{i : \bar a_{i} = 1}^{r} e_{ii} > \bar m$. We argue that we have already done all the hard work to justify \eqref{eq:vdot_case2}. Indeed, the same logic used to derive \eqref{eq:l1}--\eqref{eq:l3} implies that
\begin{align*}
 \frac{1}{2}\dot V(e,f) =&\ \sum_{i = 1}^{r} \bigg(-\lambda_i \big( 1 - \dot u_i\big) + \sum_{\substack{j = 1 \\ j \neq i}}^{r} \mu_{ji} e_{ji} + Q_{ii} \sum_{j = 1}^{r} \mu_{ji} f_{ji}\bigg)  \\
&+ \sum_{i=1}^{r} \sum_{\substack{j = 1\\ j \neq i}}^{r} \bigg( - \mu_{ij} e_{ij} + Q_{ij} \sum_{k = 1}^{r} \mu_{ki} f_{ki} \bigg) 1\big( \hat e_{ij} > 0\big)\\
&+\sum_{i=1}^{r} \sum_{j=1}^{r} \Big(\lambda_{i} P_{ij}\big( 1 - \dot u_i\big) - \mu_{ij} f_{ij}\Big) 1\big(\hat f_{ij} >  0\big).
\end{align*}
The difference between the equation above and that of \eqref{eq:l1}--\eqref{eq:l3} is that the summation in the first line is over all $i = 1, \ldots, r$, as opposed to only those $i$ for which $\bar a_{i} < 1$. It can therefore be verified, by repeating the logic of this proof, that $\frac{1}{2} \dot V(e,f)$ equals \eqref{eq:nbegin}--\eqref{eq:nend} with summations over all $i = 1, \ldots, r$, instead of only $i$ such that $\bar a_i < 1$. This verifies \eqref{eq:vdot_case2} and concludes the proof of this lemma.
\endproof

\subsection{Proof of Lemma~\ref{lem:khallemma}}
\label{sec:khallemma}
\proof{Proof of Lemma~\ref{lem:khallemma}.}
Lemma~\ref{lem:khallemma} follows from the argument used in \cite[Lemma 4.1]{Khal2002} after the following observation. In \cite{Khal2002}, the author states Lemma 4.1 for functions $x : \R \to \R^n$ that satisfy
\begin{align}
\dot x(t) = g(x(t)), \label{eq:alt_def}
\end{align}
for some locally Lipschitz $g: \R^n \to \R^n$. Only three properties of $x(t)$ are really necessary for Lemma 4.1: that $x(t)$ is guaranteed to exist, be unique, and be continuous with respect to its initial condition. The latter property means that if $x(t)$ and $\tilde x(t)$ both satisfy \eqref{eq:alt_def} with $x(0) \neq \tilde x(0)$, then for every $t \geq 0$ and every $\epsilon > 0$, there exists a $\delta > 0$ such that
\begin{align}
\abs{x(0) - \tilde x(0)} < \delta \quad \Rightarrow \quad \abs{x(t) - \tilde x(t)} < \epsilon. \label{eq:ctny_init_cond}
\end{align}
 In our case, $\big(e(t),f(t)\big)$ are defined as the solution to \eqref{eq:ufluid_fij}--\eqref{eq:ufluid_regulator} and cannot be written in the form of \eqref{eq:alt_def}. In fact, $\big(\dot e(t), \dot f(t) \big)$ is not even guaranteed to exist for all $t \geq 0$. Nevertheless, we know that $\big(e(t),f(t)\big)$ is unique by Lemma~\ref{lem:continuity}, and that it satisfies the analogue of \eqref{eq:ctny_init_cond} by \eqref{eq:ups_lipschitz} of Lemma~\ref{lem:continuity}. Hence, the proof of Lemma 4.1 can be carried through to prove Lemma~\ref{lem:khallemma} as well.
\endproof

\section{Two Examples of Networks}
In this section we describe both the 9-region, and 5-region networks we used for the numerical examples in Section~\ref{sec:numerical}.
\subsection{9-Region Didi Network} \label{app:nineregion}
 The Di-Tech Challenge data set contains individual order information for trips taken between January 1, $2016$ until January 21, 2016, in an unspecified city in China. The city is partitioned into a number of distinct geographical regions. The data is divided into $10$-minute time slots, i.e.\ $144$ times slots per day. Each data entry represents a single order. An order is a passenger request for a car, and may or may not be fulfilled due to lack of cars in proximity. A single data entry contains information about the origin and intended destination of the order, a time-stamp of when the order was made, whether the order was fulfilled, and in the case when the order was fulfilled, the ID of the car fulfilling the order, as well as the total trip price. \bluetwo{Although the data set had more than nine regions, we focus on these because they are the `major' ones. That is, they have a much higher volume of request rates compared to the rest of the regions.}



\bluetwo{We restrict our attention to data between 5-6pm each day, because we identified this to be the time of the evening rush hour; see Figure~\ref{fig:arrivals} in Section~\ref{sec:lookahead}. Using the data set, we extract a nine-region network. We estimate $\lambda$, $P$, and $\mu$ as follows.}
To calculate $P_{ij}$, we consider all rides happening in the 5-6pm window each day. We tally the total number of orders from $i$ to $j$, and divide by the total number of orders originating at $i$. The result is the matrix $P$, which equals
\begin{align}
\left(
            \begin{array}{c|c|c|c|c|c|c|c|c|c}
            \hbox{Region} & 10 & 11 & 18 & 13 & 19 & 27 & 45 & 47 & 50\\
            \hline
		  10 & 0.230 & 0.297 & 0.372 & 0.004 & 0.026 & 0.029 & 0.009 & 0.018 & 0.015 \\
		  11 & 0.044 & 0.655 & 0.146 & 0.005 & 0.079 & 0.038 & 0.018 & 0.005 & 0.011 \\
		  18 & 0.165 & 0.291 & 0.288 & 0.007 & 0.054 & 0.126 & 0.017 & 0.025 & 0.027 \\
		  13 & 0.0013 & 0.010 & 0.006 & 0.139 & 0.031 & 0.185 & 0.101 & 0.117 & 0.409 \\
		  19 & 0.005 & 0.096 & 0.026 & 0.037 & 0.25 & 0.333 & 0.218 & 0.012 & 0.027 \\
		  27 & 0.004 & 0.031 & 0.032 & 0.088 & 0.121 & 0.426 & 0.148 & 0.059 & 0.092 \\
		  45 & 0.002 & 0.023 & 0.011 & 0.066 & 0.142 & 0.269 & 0.399 & 0.020 & 0.069 \\
		  47 & 0.004 & 0.008 & 0.023 & 0.067 & 0.011 & 0.095 & 0.019 & 0.400 & 0.374 \\
		  50 & 0.001 & 0.004 & 0.005 & 0.095 & 0.010 & 0.059 & 0.030 & 0.185 & 0.610 \\
            \end{array}
          \right). \label{eq:pdidi}
\end{align}
To calculate $\mu$, we used the average trip cost as a proxy for travel times, since travel times are not provided in the data set. Trip costs are a reasonable proxy because the price of a trip is typically a linear function of distance traveled, and time spent in car. To estimate $\mu_{ij}$, we first calculated the average trip cost between regions $i$ and $j$, and then set the average travel time to equal the average trip cost. Since our time unit is a time-slot, which is $10$-minute interval, we set $1/\mu_{ij}$ to equal the average trip cost divided by $10$. For example, the average trip cost between region 47 and region 50 is 14.1 CNY, so we assumed the average travel time is 14 minutes, and set $\mu_{47,50} = 1.41$. The resulting matrix $1/\mu$ is
\begin{align}
\left(
            \begin{array}{c|c|c|c|c|c|c|c|c|c}
            \hbox{Region} & 10 & 11 & 18 & 13 & 19 & 27 & 45 & 47 & 50\\
            \hline
		  10 & 0.83 & 1.87 & 1.07 & 3.89 & 3.25 & 2.79 & 4.25 & 2.94 & 4.37 \\
		  11 & 1.78 & 0.89 & 1.18 & 3.24 & 1.24 & 1.99 & 2.89 & 3.46 & 4.18 \\
		  18 & 1.02 & 1.31 & 0.78 & 2.82 & 1.45 & 1.36 & 3.26 & 2.17 & 3.04 \\
		  13 & 3.52 & 3.13 & 2.76 & 0.93 & 1.5 & 1.26 & 1.49 & 1.75 & 1.6
 \\
		  19 & 2.86 & 1.42 & 1.64 & 1.55 & 0.84 & 1.04 & 1.45 & 2.88 & 2.89 \\
		  27 & 2.61 & 2.17 & 1.54 & 1.31 & 1.15 & 0.81 & 1.86 & 1.78 & 2.2 \\
		  45 & 4.38 & 3.02 & 2.79 & 1.36 & 1.35 & 1.65 & 0.94 & 3.1 & 3 \\
		  47 & 2.93 & 3.06 & 2.26 & 1.75 & 2.69 & 1.62 & 3.23 & 0.9 & 1.48 \\
		  50 & 3.58 & 4.18 & 2.8 & 1.49 & 2.46 & 2.02 & 2.72 & 1.43 & 1.01 \\
            \end{array}
          \right). \label{eq:mudidi}
\end{align}
To determine the arrival rate to region $i$, we counted the average number of orders to region $i$ per time slot. This gave us an estimate $N \lambda_i$. Our data did not provide the exact number $N$ of cars in the network. To determine a reasonable choice for $N$, we summed up the number of fulfilled orders across all 9 regions in Figure~\ref{fig: traffic}. As a result, we chose $N = 2000$. Although not exact, this number is of the correct order of magnitude. \bluetwo{Furthermore, our numerical results in Section~\ref{sec:numerical} remained (qualitatively) consistent even for different choices of $N$ around $2000$.} Hence, we set our vector $\lambda$ to
\begin{align}
\left(
            \begin{array}{cccccccccc}
            \hbox{Region $i$} & 10 & 11 & 18 & 13 & 19 & 27 & 45 & 47 & 50 \\ \hline
            \lambda_i &  0.0131 & 0.0624 & 0.0381& 0.0652 & 0.0870 & 0.1178 & 0.0762 & 0.1438 & 0.2751 \\
            \end{array}
          \right) \label{eq:lambdadidi}
\end{align}

\subsection{5-Region Network} \label{app:fiveregion}
\begin{figure}[h] 
    \includegraphics[width=0.45\textwidth]{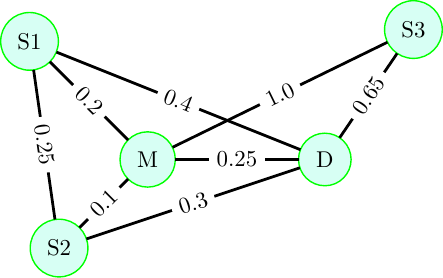}
    \centering
    \caption{The 5-region city where edge weights represent mean travel times (in hours). Travel times longer than $1$ hour are omitted to help with vizualisation.} \label{fig:toyTopology}
\end{figure}

Let us consider a simplified model of a city illustrated in Figure~\ref{fig:toyTopology} that consists of $5$ regions: a downtown area $D$, a midtown area $M$, and three suburban areas $S_1, S_2, S_3$.
The downtown area represents a central business district, where many people work but few people live. Midtown represents a region with restaurants and night-life, where people visit after work. The suburb regions are residential areas, and do not have as many entertainment options as midtown. For convenience, we enumerate $S_1, S_2, S_3, M,D$ as $1,2,3,4,5$, respectively.
 The parameters from 5-11pm are as follows:
\begin{enumerate}
\item From 5-7pm, the city experiences a rush hour as people go home from work. Most of the traffic originates in downtown, and flows into the suburbs as people go home after work. The parameters during this time-slot are
\begin{align*}
\lambda = 
\begin{pmatrix}
0.108 \\
0.108 \\
0.108 \\
0.108 \\
1.08
\end{pmatrix}, \quad 
P = 
\begin{pmatrix}
0.6 & 0.1 & 0 & 0.3 & 0 \\
0.1 & 0.6 & 0 & 0.3 & 0 \\
0 & 0 & 0.7 & 0.3 & 0 \\
0.2 & 0.2 & 0.2 & 0.2 & 0.2 \\
0.3 & 0.3 & 0.3 & 0.1 & 0
\end{pmatrix}, \quad 
\big(1/\mu_{ij}\big) = 
\begin{pmatrix}
0.15 & 0.25 & 1.25 & 0.2 & 0.4 \\
0.25 & 0.10 & 1.1 & 0.1 & 0.3 \\
1.25 & 1.1 & 0.1 & 1 & 0.65 \\
0.25 & 0.15 & 1 & 0.15 & 0.25 \\
0.5 & 0.4 & 0.75 & 0.25 & 0.2
\end{pmatrix}.
\end{align*}
Above, $\big(1/\mu_{ij}\big)$ is the matrix of mean travel times (in hours) from region $i$ to $j$, and $\lambda$ is a vector representing the number of passenger arrivals per hour, per car.
\item From 7-9pm, most of the traffic is headed into midtown as people go out in the evening to restaurants and for entertainment.  The parameters during this time-slot are
\begin{align*}
\lambda = 
\begin{pmatrix}
0.72 \\
0.48 \\
0.48 \\
0.48 \\
0.12
\end{pmatrix}, \quad 
P = 
\begin{pmatrix}
0.1 & 0 & 0 & 0.9 & 0 \\
0 & 0.1 & 0 & 0.9 & 0 \\
0 & 0 & 0.1 & 0.9 & 0 \\
0.05 & 0.05 & 0.05 & 0.8 & 0.05 \\
0 & 0 & 0 & 0.9 & 0.1
\end{pmatrix}, \quad 
\big(1/\mu_{ij}\big) = 
\begin{pmatrix}
0.15 & 0.25 & 1.25 & 0.2 & 0.4 \\
0.25 & 0.10 & 1.1 & 0.1 & 0.3 \\
1.25 & 1.1 & 0.1 & 1 & 0.65 \\
0.2 & 0.1 & 1 & 0.15 & 0.25 \\
0.4 & 0.3 & 0.65 & 0.25 & 0.2
\end{pmatrix}
\end{align*}
\item From 9-11pm, traffic flows mainly from midtown to the suburbs as people go home for the night.  The parameters during this time-slot are
\begin{align*}
\lambda = 
\begin{pmatrix}
 0.12 \\   
 0.12 \\
 0.12 \\
 1.32 \\
 0.12
\end{pmatrix}, \quad 
P = 
\begin{pmatrix}
0.9 & 0.05 & 0 & 0.05 & 0 \\
0.05 & 0.9 & 0 & 0.05 & 0 \\
0 & 0 & 0.9 & 0.1 & 0 \\
0.3 & 0.3 & 0.3 & 0.05 & 0.05 \\
0 & 0 & 0 & 0.1 & 0.9
\end{pmatrix}, \quad 
\big(1/\mu_{ij}\big) = 
\begin{pmatrix}
0.15 & 0.25 & 1.25 & 0.2 & 0.4 \\
0.25 & 0.10 & 1.1 & 0.1 & 0.3 \\
1.25 & 1.1 & 0.1 & 1 & 0.65 \\
0.2 & 0.1 & 1 & 0.15 & 0.25 \\
0.4 & 0.3 & 0.65 & 0.25 & 0.2
\end{pmatrix}
\end{align*}
\end{enumerate}

\section{\blue{Achieving 100\% Availability Globally}} 
\label{sec:perfectavailability}
It may be the case, e.g.\ during rush hours, that customer demand is so high that there are simply not enough cars in the system to fulfill all passenger requests, no matter what kind of routing policy is used. A natural question to ask then is how many cars does one need to achieve $100\%$ availability in every region? In this section we address this question using the fluid-based optimization problem. 

Fix $N > 0$ and some static routing policy $Q$. Since the CTMC $\big(E^{(N)}, F^{(N)}\big)$ is positive recurrent and has a finite state space, it must be the case that $A_{i}^{(N)}(\infty) < 1$. In other words, with a finite number of cars, the availability at any region can never be $100\%$. However, it is possible to have $A_{i}^{(N)}(\infty) \to 1$ as $N \to \infty$. We say that our system (asymptotically) achieves perfect availability (under routing policy $Q$) if 
\begin{align*}
\lim_{N \to \infty} A_{i}^{(N)}(\infty)  = 1, \quad 1 \leq i \leq r.
\end{align*}
We also say that perfect availability is feasible if there exists \emph{some} routing matrix $Q$ under which the system achieves perfect availability. 

The question of whether perfect availability is feasible is a question of whether there is enough supply of cars to meet passenger demand. If perfect availability is not feasible, we want to know the amount by which to increase our fleet size in order to make it feasible. Conversely, if perfect availability is feasible, we want to know how much slack, or excess capacity, our system has. Since the passenger arrival rate to each region, $N \lambda_i$, depends on the number of cars $N$, we need to clarify what we mean by increasing or decreasing the fleet size. 

Consider a system with $N$ cars, passenger arrival rate $N \lambda_i$, travel choices  $P_{ij}$, and travel time means $1/\mu_{ij}$, and recall that the associated fluid-based optimization problem is given by (\hyperref[eqlbm:objective]{4a})--(\hyperref[eqlbm:prob]{4g}). We can consider a related system, where all parameters are the same except that the number of cars is now $\kappa N$ for some $\kappa > 0$. In this system, $N$ is no longer the number of cars, but represents the size of the market.  The passenger arrival rate to region $i$ is $N \lambda_i = \kappa N (\lambda_i/\kappa)$, and therefore, the associated fluid-based optimization problem is (\hyperref[eqlbm:objective]{4a})--(\hyperref[eqlbm:prob]{4g}), but with $\lambda_i$ replaced by $\lambda_i/\kappa$ there. That is, multiplying the number of cars in the system by a factor of $\kappa$ is equivalent to dividing $\lambda_i$ by a factor of $\kappa$. This makes sense because $\lambda_i$ is the arrival rate of passengers \emph{per car} to region $i$. 

The following result states that determining the balance between supply and demand needed to achieve perfect availability can be done by solving a linear program. The proof can be found in Appendix~\ref{app:fleetsize}.


\begin{theorem} \label{thm:fleetsize}
The feasibility region of the linear program 
\begin{subequations}
\begin{align}
&\min_{q} \sum_{i=1}^{r}\sum_{j=1}^{r}\bar f_{ij} + \sum_{i=1}^{r}\sum_{j=1,j\neq i}^{r}\bar e_{ij} \label{feas:obj} \\
\text{ subject to } \quad & \lambda_iP_{ij} = \mu_{ij}\bar f_{ij}, \quad 1\leq i,j \leq r, \label{feas:one} \\
&\mu_{ij}\bar e_{ij} = q_{ij}\sum_{k=1}^{r}\mu_{ki}\bar f_{ki}, \quad 1\leq i,j \leq r, ~j\neq i, \label{feas:two} \\
&\lambda_i = \sum_{k=1,k\neq i}^{r}\mu_{ki}\bar e_{ki} + q_{ii}\sum_{k=1}^{r}\mu_{ki}\bar f_{ki}, \quad 1\leq i \leq r, \label{feas:three} \\
&q_{ij}\geq 0, \quad \sum_{j=1}^{r}q_{ij} = 1, \quad 1\leq i,j \leq r. \label{feas:four}
\end{align}
\end{subequations}
is non-empty. Let $\kappa > 0$ be the optimal objective value, and assume that 
\begin{align*}
\frac{1}{\mu_{ik}} \leq  \frac{1}{\mu_{ij}} + \frac{1}{\mu_{jk}}, \quad 1 \leq i \neq j \neq k \leq r,
\end{align*}
i.e.\ travel time means satisfy the triangle inequality. 
\begin{enumerate}
\item If $\kappa > 1$ then perfect availability is not feasible, but it becomes feasible if $\lambda$ is reduced to  $\lambda / \kappa$, i.e.\ an increase in fleet size by a factor of $\kappa$. 
\item If $\kappa \leq 1$ then perfect availability is feasible, and it remains feasible even if $\lambda$ is increased to $\lambda/\kappa$.
\end{enumerate}
\end{theorem}
The value $\kappa$ from Theorem~\ref{thm:fleetsize} is the minimal `fluid mass' needed in the system to achieve perfect availability, and can be interpreted a measure of imbalance of supply and demand in the system. If $\kappa > 1$, then passenger demand exceeds vehicle supply, and if $\kappa < 1$ then the reverse is true.


\subsection{Proof of Theorem~\ref{thm:fleetsize}}
\label{app:fleetsize}
\proof{Proof of Theorem~\ref{thm:fleetsize}}
To see that the feasibility region is non-empty (and hence an optimal solution exists), observe that 
\begin{align}
q_{ij} = \frac{\mu_{ji}f_{ji}}{\sum_{k=1}^{r} \mu_{ki} f_{ki} }, \quad 1 \leq i,j \leq r \label{eq:feas}
\end{align}
is a feasible solution. The intuition behind deriving \eqref{eq:feas} is that every car that completes a trip from $i$ to $j$ with a passenger must drive back empty from $j$ to $i$.

Recall that $\kappa > 0$ is the optimal value of the linear program, and let $q^\kappa $ be the solution that achieves this minimum. Also let $\bar e^\kappa $ and $\bar f^\kappa $ be the corresponding values of $\bar e$ and $\bar f$ under $q^\kappa $, i.e.\ the optimal value $\kappa  = \sum_{i=1}^{r}\sum_{j=1}^{r}\bar f_{ij}^{\kappa } + \sum_{i=1}^{r}\sum_{j=1,j\neq i}^{r}\bar e_{ij}^{\kappa }$; observe that \eqref{feas:one}--\eqref{feas:four} place no constraints on $\bar e^{\kappa }_{ii}$, so we will choose our $\bar e^{\kappa }$ with $\bar e^{\kappa }_{ii} = 0$. Lastly, let $\bar a^{\kappa } \in \R^{r}$ be a vector whose elements all equal to one. 

Suppose first that $\kappa  > 1$. To show that perfect availability is not feasible, we argue that if $(q,\bar e, \bar f, \bar a)$ is a point in the feasible region of the fluid-based optimization (\hyperref[eqlbm:objective]{4a})--(\hyperref[eqlbm:prob]{4g}), then it must be that $\bar a_i < 1$ for some $i$. Assume this is not the case, i.e.\ $(q,\bar e, \bar f, \bar a)$ satisfies the constraints of the fluid-based optimization and $\bar a_i = 1$ for all $1 \leq i \leq r$. Then (\hyperref[eqlbm: fullcsv]{4b})--(\hyperref[eqlbm: emptycsv2]{4d}) and (\hyperref[eqlbm:prob]{4g}) are identical to \eqref{feas:one}--\eqref{feas:three} and \eqref{feas:four}, respectively. Since the optimal value in \eqref{feas:obj} is greater than one, it means that constraint (\hyperref[eqlbm:carnum]{4f}) of the fluid-based optimization can never be satisfied under any $q$, which proves that perfect availability is not feasible.

Now consider a new system where $\lambda_i$ is replaced instead by $\lambda_i/\kappa $. It can be checked that \eqref{feas:one}--\eqref{feas:four} in this new system is satisfied by $q^{\kappa }$, and that the resulting objective value is 
\begin{align*}
\frac{1}{\kappa }\Big( \sum_{i=1}^{r}\sum_{j=1}^{r}\bar f_{ij}^{\kappa } + \sum_{i=1}^{r}\sum_{j=1,j\neq i}^{r}\bar e_{ij}^{\kappa }\Big) = 1.
\end{align*} 
One can then check that $(q^{\kappa }, \bar e^\kappa  / \kappa , \bar f^\kappa  / \kappa , \bar a^\kappa )$ is a feasible solution to the fluid-based optimization problem of this new system. Assuming for now that $q^{\kappa }_{ii} > 0$ for all $1 \leq i \leq r$, we can invoke Theorem~\ref{thm:asymptotic_utility} to conclude that the new system achieves perfect availability under $q^\kappa $. 

Now if $\kappa  \leq 1$, a similar argument can be used to see that the system achieves perfect availability under routing matrix $q^\kappa $, i.e.\ one confirms that $(q^\kappa , \bar e^{\kappa }, \bar f^{\kappa }, \bar a^{\kappa })$ is a feasible solution to the fluid-based optimization. Furthermore, if we consider a new system where $\lambda_i$ is replaced by $\lambda_i / \kappa $, then the new system still achieves perfect availability when $q^\kappa $ is taken to be the routing matrix; this can again be verified just like in the $\kappa  > 1$ case. 

To conclude the proof, we need to show that we can always find an optimal solution to \eqref{feas:obj}--\eqref{feas:four} such that $q_{ii}^{\kappa } > 0$ for all $1 \leq i \leq r$ (this is a minor technical condition needed to invoke Theorem~\ref{thm:asymptotic_utility}). Suppose $q^{\kappa }$ is an optimal solution and $q_{ii}^{\kappa } = 0$ for some $i$; we now construct another optimal solution $\hat q^{\kappa }$ such that $\hat q^{\kappa }_{ii} > 0$. We know that 
\begin{align*}
0 < \lambda_i = \sum_{k=1,k\neq i}^{r}\mu_{ki}\bar e_{ki}^{\kappa } + q_{ii}^{\kappa }\sum_{k=1}^{r}\mu_{ki}\bar f_{ki} = \sum_{k=1,k\neq i}^{r}\mu_{ki}\bar e_{ki}^{\kappa },
\end{align*}
where the first equality is from \eqref{feas:three}. The above implies that there must exist some $\ell$ such that $\bar e_{\ell i}^{\kappa } > 0$, and consequently (by \eqref{feas:two} and \eqref{feas:one}) $q^{\kappa }_{\ell i} > 0$. Furthermore, since $q^{\kappa }_{ii} = 0$, there must exist some $m$ such that $q^{\kappa }_{im} > 0$, and hence $\bar e_{im}^{\kappa } > 0$. 

The idea of the following argument is to redirect some empty cars going from $\ell$ to $i$ to instead go from $\ell$ to $m$, i.e.\ reduce $q^{\kappa }_{\ell i}$ to increase $q^{\kappa }_{\ell m}$. This will allow region $i$ to keep some of the cars it otherwise would have sent to $m$, i.e.\ reduce $q^{\kappa }_{im}$ to increase $q^{\kappa }_{ii}$. It remains to specify precisely the changes to $q^{\kappa }$ so that the objective value in \eqref{feas:obj} does not increase. To this end, fix $\varepsilon > 0$ and let 
\begin{align*}
\hat q^{\kappa }_{\ell i} =&\ q^{\kappa }_{\ell i} - \varepsilon\\
\hat q^{\kappa }_{\ell m} =&\ q^{\kappa }_{\ell m} + \varepsilon\\
\hat q^{\kappa }_{im} =&\ q^{\kappa }_{im} - \varepsilon \frac{\sum_{k=1}^{r} \lambda_k P_{k\ell}}{\sum_{k=1}^{r} \lambda_k P_{ki}} \\
\hat q^{\kappa }_{ii} =&\ q^{\kappa }_{ii} + \varepsilon \frac{\sum_{k=1}^{r} \lambda_k P_{k\ell}}{\sum_{k=1}^{r} \lambda_k P_{ki}},
\end{align*}
and let all other elements of  $\hat q^{\kappa }$ be the same as those of $q^{\kappa }$. Provided $\varepsilon$ is small enough,  $\hat q^{\kappa }$ satisfies \eqref{feas:four}, i.e.\ its rows are probability distributions. It is also not hard to check that $\hat q^{\kappa }$ satisfies \eqref{feas:one}--\eqref{feas:three}, and that the objective value under $\hat q^{\kappa }$ is 
\begin{align*}
\sum_{i=1}^{r}\sum_{j=1}^{r}\bar f_{ij}^{\kappa } + \sum_{i=1}^{r}\sum_{j=1,j\neq i}^{r}\bar e_{ij}^{\kappa } + \varepsilon \sum_{k=1}^{r} \lambda_k P_{k\ell} \Big(\frac{1}{\mu_{\ell m}} - \frac{1}{\mu_{\ell i}} - \frac{1}{\mu_{im}}\Big),
\end{align*}
which is not larger than the objective value under $q^{\kappa }$ by our assumption that the mean travel times satisfy the triangle inequality. Therefore, $\hat q^{\kappa }$ is an optimal solution, and this concludes the proof.
\endproof

\section{Miscellaneous Proofs}
\subsection{Proof of Lemma~\ref{lem:relax0}}\label{app:relax0}
\proof{Proof of Lemma \ref{lem:relax0}.}
Given $(\bar{e}, \bar{f}, \bar{a})$ and $q$ that satisfy (\hyperref[eqlbm: fullcsv]{4b})--(\hyperref[eqlbm:prob]{4g}), it can be easily verified that $(\bar{e}, \bar{f}, \bar{a})$ satisfies conditions \eqref{rrp: fullcsv}--\eqref{rrp: boundry} based on the fact that $0\leq q_{ij}\leq 1$.

Now given $(\bar{e}, \bar{f}, \bar{a})$ that satisfies \eqref{rrp: fullcsv}--\eqref{rrp: boundry}, we define
\begin{align*}
q_{ij}=&\ \frac{\mu_{ij}\bar e_{ij}}{\sum_{k=1}^{r}\mu_{ki}\bar f_{ki}}, \quad 1 \leq i\neq j \leq r, \quad q_{ii}= \frac{\lambda_i\bar a_i -\sum_{k=1,k\neq i}^{r}\mu_{ki}\bar e_{ki}}{\sum_{k=1}^{r}\mu_{ki}\bar f_{ki}} \quad 1 \leq i \leq r,
 \end{align*}
and let $q$ be the $r\times r$ matrix whose entries are $q_{ij}$.   Then conditions (\hyperref[eqlbm: emptycsv1]{4c}) and (\hyperref[eqlbm: emptycsv2]{4d}) hold according to the definition of $q$. Furthermore, $q_{ij}\geq 0$ because $\mu_{ij}\geq 0,$ $\bar e_{ij}\geq 0$, and $\lambda_i\bar a_i\geq \sum_{k=1,k\neq i}^{r}\mu_{ki}\bar e_{ki}$ by (\ref{rrp:prob2}). Finally,
\begin{eqnarray*}
\sum_{j=1}^r q_{ij}&=&\sum_{j=1, j\not=i}^r \frac{\mu_{ij}\bar e_{ij}}{\sum_{k=1}^{r}\mu_{ki}\bar f_{ki}}+\frac{\lambda_i\bar a_i -\sum_{k=1,k\neq i}^{r}\mu_{ki}\bar e_{ki}}{\sum_{k=1}^{r}\mu_{ki}\bar f_{ki}}\\
&=&\frac{\lambda_i\bar a_i+ \sum_{j=1, j\not=i}^r\mu_{ij}\bar e_{ij}- \sum_{k=1,k\neq i}^{r}\mu_{ki}\bar e_{ki}}{\sum_{k=1}^{r}\mu_{ki}\bar f_{ki}}\\
&\overset{(a)}{=}&\frac{\sum_{k=1}^{r}\mu_{ki}\bar f_{ki}}{\sum_{k=1}^{r}\mu_{ki}\bar f_{ki}}\\
&=&1,
\end{eqnarray*} where equality (a) is obtained based on (\ref{rrp: carcsv}).  Therefore, (\hyperref[eqlbm:prob]{4g}) holds and we can conclude that $(\bar{e}, \bar{f}, \bar{a})$ and our newly defined $q$ satisfy (\hyperref[eqlbm: fullcsv]{4b})--(\hyperref[eqlbm:prob]{4g}).
\endproof

\subsection{Proof of Lemma~\ref{lem:relax}}\label{app:relax}
\proof{Proof of Lemma~\ref{lem:relax}.}
 Suppose $\bar a_i^*<1$, for all $1 \leq i \leq r$. We argue by contradiction that this implies $\bar e_{ii}^*=0$, for all $i$. Suppose there exists region $i'$ such that  $\bar e^*_{i'i'} > 0$. We now construct another solution $({\tilde e}, {\tilde f}, {\tilde a})$ that is better than $(\bar{e}^*, \bar{f}^*, \bar{a}^*)$. Assume for now that there exists an $r\times r$ matrix with non-negative entries $\pi_{ij}$, such that $\sum_{i=1}^{r}\sum_{j=1}^{r} \pi_{ij} = 1$, and
\begin{align}
P_{ij}\sum_k \pi_{ki}\mu_{ki}= \mu_{ij}\pi_{ij}, \quad 1 \leq i,j \leq r. \label{eq:flow_bal}
\end{align}
Fix $\epsilon > 0$ to be specified later, and let
\begin{align*}
\tilde{e}_{i'i'}=&\ \bar{e}_{i'i'}^*-\epsilon \\
 \tilde{e}_{ii}=&\ \bar{e}_{ii}^*, \quad  i \neq i', \\
\tilde f_{ij}=&\ \bar f_{ij}^* + \epsilon \pi_{ij} , \quad 1\leq i,j \leq r,\\
\tilde e_{ij}=&\ \bar e_{ij}^*, \quad i \neq j,\\
\tilde a_i=&\ \bar a_i^*+\epsilon\frac{\sum_{k=1}^r \pi_{ki}\mu_{ki}}{\lambda_i}, \quad 1\leq i \leq r.
\end{align*}
Since $\tilde{a}_i\geq \bar{a}^*_i$ for all $i$, \delete{and $\tilde{f}_{ij}\geq \bar f^*_{ij}$ for all $i,j$, and $\tilde e_{ij}=\bar e_{ij}^*$ for $i\not=j,$} it follows that
\blue{
\begin{align}
\sum_{i=1}^{r} \sum_{j=1}^{r} \bar a_i^{*}\lambda_i P_{ij} c_{ij} \leq \sum_{i=1}^{r} \sum_{j=1}^{r} \tilde a_i\lambda_i P_{ij} c_{ij}.
\end{align}}
We now check that $(\tilde{e}, \tilde{f}, \tilde{a})$ satisfies \eqref{rrp: fullcsv}--\eqref{rrp: frac}, and is therefore a feasible solution. Since $(\bar e^*,\bar f^*,\bar a^*)$ is a feasible solution and satisfies \eqref{rrp: fullcsv}, it follows that
\begin{align*}
\lambda_iP_{ij}\tilde a_i=\lambda_iP_{ij}\left(\bar a_i^* +\epsilon\frac{\sum_{k=1}^r \pi_{ki}\mu_{ki}}{\lambda_i}\right)=&\  \mu_{ij}\bar f_{ij}^*+\epsilon P_{ij}{\sum_{k=1}^r \pi_{ki}\mu_{ki}} \\
=&\ \mu_{ij}\bar f_{ij}^*+\epsilon\mu_{ij}\pi_{ij} \\
=&\  \mu_{ij} \tilde f_{ij}, \quad 1 \leq i,j \leq r,
\end{align*}
meaning $(\tilde{e}, \tilde{f}, \tilde{a})$ satisfies \eqref{rrp: fullcsv}. Next, we see that
\begin{align*}
&\mu_{ij}\tilde e_{ij} =\mu_{ij}\bar{e}_{ij}^* \leq \sum_{k=1}^{r}\mu_{ki}\bar f^*_{ki} \leq \sum_{k=1}^{r}\mu_{ki}\tilde f_{ki}, \quad 1 \leq i \neq j \leq r,
\end{align*}
where the first inequality follows because $(\bar e^*,\bar f^*,\bar a^*)$ satisfies \eqref{rrp:prob1}. Thus, $(\tilde{e}, \tilde{f}, \tilde{a})$ satisfies \eqref{rrp:prob1}. Next we check constraint \eqref{rrp:prob2}. Since  $(\bar e^*,\bar f^*,\bar a^*)$ satisfies \eqref{rrp:prob2}, we know that
\begin{align*}
\sum_{k=1,k\neq i}^r\mu_{ki} \tilde{e}_{ki}   = \sum_{k=1,k\neq i}^r\mu_{ki} e^{*}_{ki}
 \leq  \lambda_i\bar a_i^*
  \leq \lambda_i\tilde a_i, \quad 1 \leq i \leq r,
\end{align*}
and
\begin{align*}
\lambda_i\tilde a_i =\lambda_i\left(\bar a_i^*+\epsilon\frac{\sum_{k=1}^r \pi_{ki}\mu_{ki}}{\lambda_i}\right) \leq&\  \sum_{k=1,k\neq i}^{r}\mu_{ki}\bar e^*_{ki} + \sum_{k=1}^{r}\mu_{ki}\bar f^*_{ki} + \epsilon\sum_{k=1}^r \pi_{ki}\mu_{ki} \\
=&\  \sum_{k=1,k\neq i}^{r}\mu_{ki}\bar e^*_{ki} + \sum_{k=1}^{r}\mu_{ki}(\bar f^*_{ki} + \epsilon \pi_{ki}) \\
=&\  \sum_{k=1,k\neq i}^r\mu_{ki} \tilde{e}_{ki} + \sum_{k=1}^{r}\mu_{ki} \tilde{f}_{ki}, \quad 1 \leq i \leq r.
\end{align*}
Therefore $(\tilde{e}, \tilde{f}, \tilde{a})$ satisfies \eqref{rrp:prob2}. To verify \eqref{rrp: carcsv}, observe that
\begin{align*}
&\lambda_i\tilde a_i + \sum_{j=1,j\neq i}^r\mu_{ij}\tilde e_{ij} - \sum_{k=1,k\neq i}^r\mu_{ki}\tilde e_{ki} - \sum_{k=1}^{r}\mu_{ki}\tilde f_{ki}\\
=&\lambda_i\bar a^*_i + \sum_{j=1,j\neq i}^r\mu_{ij}\bar e^*_{ij} - \sum_{k=1,k\neq i}^r\mu_{ki}\bar e^*_{ki} - \sum_{k=1}^{r}\mu_{ki}\bar f^*_{ki} + \epsilon {\sum_{k=1}^r \pi_{ki}\mu_{ki}}-\epsilon \sum_{k=1}^{r}\mu_{ki}\pi_{ki}\\
=&0,
\end{align*}
and therefore \eqref{rrp: carcsv}  holds under $(\tilde e, \tilde f, \tilde a).$ Lastly, $\epsilon$ can be chosen small enough so that both \eqref{rrp: carnum} and \eqref{rrp: frac} hold (the former also uses the fact that $\sum_{ij}\pi_{ij}=1$). We conclude that $(\tilde e, \tilde f, \tilde a)$ is a feasible solution to \eqref{opt-r:obj}--\eqref{opt-r:cdts} that is better than $(\bar e^*, \bar f^*, \bar a^*)$, which contradicts the fact that  $(\bar e^*, \bar f^*, \bar a^*)$ is an optimal solution.

It remains to verify that we can choose non-negative $\pi_{ij}$'s to satisfy $\sum_{ij} \pi_{ij} = 1$ and \eqref{eq:flow_bal}. Consider a CTMC defined on the space $\{1, \cdots, r\}^2$. For all $1 \leq i,j,k \leq r$, the transition rate from $(k,i)$ to $(i,k)$ is $\mu_{ki}P_{ik}$. No other transitions are possible. Since the CTMC is defined on a finite state space, it has a stationary distribution. Furthermore, any stationary distribution $\nu$ must satisfy the flow-balance equations
\begin{align*}
\sum_{k=1}^{r} \nu_{ki} \mu_{ki} P_{ij} =  \sum_{\ell = 1}^{r} \nu_{ij} \mu_{ij}  P_{j \ell}, \quad 1 \leq i,j \leq r,
\end{align*}
or
\begin{align*}
P_{ij}\sum_{k=1}^{r} \nu_{ki} \mu_{ki}  =  \nu_{ij} \mu_{ij}, \quad 1 \leq i,j \leq r,
\end{align*}
which are precisely the same as \eqref{eq:flow_bal}. Therefore we can take the $\pi_{ij}$'s to be any of the stationary distributions of such a CTMC.

In the case when there exists $i'$ such that $\bar{a}^*_{i'}=1$, the claim of the lemma is straightforward to verify.
\endproof

\subsection{Proof of Theorem~\ref{thm:asymptotic_utility}}
\label{app:main_thm_proof}
The proof of \eqref{eq:main1}--\eqref{eq:main4} hinges on combining Theorems~\ref{thm:main_conv} and \ref{thm:stability} with Theorem 5.1 of \cite{AnanBenc1993}. Below, we repeat the argument from \cite{AnanBenc1993} for completeness. We know that the sequence  $\{(\bar E^{(N)}(\infty), \bar F^{(N)}(\infty))\}_{N}$ is tight, because the support of $(\bar E^{(N)}(\infty), \bar F^{(N)}(\infty))$ is the compact set $\mathcal{T}$. It follows by Prohorov's Theorem \cite{Bill1999} that the sequence is also relatively compact. We will now show that any subsequence of $\{(\bar E^{(N)}(\infty), \bar F^{(N)}(\infty))\}_{N}$ has a further subsequence that converges weakly to a probability measure that assigns a mass of one to the equilibrium set $\mathcal{E}$, thereby proving \eqref{eq:main1}--\eqref{eq:main4}.

Fix $N > 0$ and initialize $(\bar E^{(N)}(0), \bar F^{(N)}(0))$ according to $(\bar E^{(N)}(\infty), \bar F^{(N)}(\infty))$. Prohorov's Theorem implies that for any subsequence
\begin{align*}
\{(\bar E^{(N')}(0), \bar F^{(N')}(0))\}_{N'} \subset \{(\bar E^{(N)}(0), \bar F^{(N)}(0))\}_{N},
\end{align*}
there exists a further subsequence
\begin{align*}
\{(\bar E^{(N'')}(0), \bar F^{(N'')}(0))\}_{N''} \subset \{(\bar E^{(N')}(0), \bar F^{(N')}(0))\}_{N'}
\end{align*}
that converges weakly to some probability measure $(e_0, f_0)$ with support in $\mathcal{T}$. Now for any $t \geq 0$,
\begin{align*}
(\bar E^{(N'')}(0), \bar F^{(N'')}(0)) \stackrel{d}{=} (\bar E^{(N'')}(t), \bar F^{(N'')}(t)) \Rightarrow (e(t), f(t)),
\end{align*}
as $N'' \to \infty$, where $(e(t), f(t))$ is the fluid model with intial condition $(e(0),f(0)) = (e_0, f_0)$, and the weak convergence follows from Theorem~\ref{thm:main_conv}. Since $(e(t), f(t))$ converges to the set $\mathcal{E}$ as $t \to \infty$, it must be the case that $(e_0, f_0) \in \mathcal{E}$ with probability one. This proves \eqref{eq:main1}--\eqref{eq:main4}.

 To prove \eqref{eq:busy_probs}, we need to use the generator of $(\bar E^{(N)}, \bar F^{(N)})$, which we call $G^{(N)}$. Since $(\bar E^{(N)}, \bar F^{(N)})$ takes values in a bounded set, Proposition 3 of \cite{GlynZeev2008} tells us that any function $g: \mathcal{T} \to \R$ satisfies
\begin{align}
\E \big[ G^{(N)} g\big(\bar E^{(N)}(\infty), \bar F^{(N)}(\infty)\big)\big] = 0. \label{bar}
\end{align}
In particular, fix $i,j$ between $1, \ldots, r$ and choose $g(e,f) = f_{ij}$. Then
\begin{align*}
G^{(N)} g(e,f) = N \lambda_i P_{ij}1(e_{ii} > 0)\big( (f_{ij}+1/N) - f_{ij}\big) + \mu_{ij}N f_{ij} \big( (f_{ij}-1/N) - f_{ij}\big),
\end{align*}
which implies
\begin{align}
\E \big[ G^{(N)} g\big(\bar E^{(N)}(\infty), \bar F^{(N)}(\infty)\big)\big] = \lambda_{i} P_{ij} \Prob(\bar E^{(N)}_{ii}(\infty) > 0) - \mu_{ij} \E \bar F^{(N)}_{ij}(\infty) = 0. \label{eq:bar_fij}
\end{align}
Hence,
\begin{align*}
\lim_{N\to\infty} \Prob(\bar E^{(N)}_{ii}(\infty) > 0) = \lim_{N\to\infty} \frac{\mu_{ij}}{\lambda_i P_{ij}} \E \bar F^{(N)}_{ij}(\infty) = \frac{\mu_{ij}}{\lambda_i P_{ij}} \bar f_{ij} = \bar a_{i},
\end{align*}
where in the second equality we used \eqref{eq:main1} and the fact that $\bar F^{(N)}_{ij}(\infty) \in [0,1]$ to conclude that the sequence of expected values  $\E [\bar F^{(N)}_{ij}(\infty)]$ converges to $\bar f_{ij}$, and in the last equality we used \eqref{eq:one}.

\subsection{Proof of Theorem~\ref{thm:lb}}\label{app:lb}
\proof{Proof of Theorem \ref{thm:lb}.}
We will show that the performance measures $\big(\E[\bar E^{(N)}(\infty)],  \E[\bar F^{(N)}(\infty)], A^{(N)}(\infty) \big)$ are a feasible solution to the optimization problem \eqref{opt-r:obj}--\eqref{opt-r:cdts}. Lemma~\ref{lem:relax} then implies part (a) of the theorem is satisfied but \blue{only with a non-strict inequality. To show the inequality is strict, note that the lemma also tells us that the optimal solution to \eqref{opt-r:obj}--\eqref{opt-r:cdts} can never be achieved by performance measures of a CTMC with finitely many cars under any routing policy. This is because the CTMC is positive recurrent and has finitely many states, and so its stationary distribution assigns positive mass to each state. It follows that $A_i^{(N)}(\infty) = \Prob(\bar E_{ii}^{(N)}(\infty) > 0 ) < 1$ and $\E [\bar E_{ii}^{(N)}(\infty)] > 0$ for all regions $i$. The first claim in Lemma~\ref{lem:relax} then prevents any performance measures from being the optimal solution, which proves the strict inequality.
}

 Part (b) of the theorem is an immediate consequence of Theorem~\ref{thm:asymptotic_utility} by setting $Q=q^*.$
Recall from \eqref{bar} that any function $g: \mathcal{T} \to \R$ satisfies
\begin{align}
\E \big[ G^{(N)} g\big(\bar E^{(N)}(\infty), \bar F^{(N)}(\infty)\big)\big] = 0. \label{bar-p}
\end{align} We now show that $\big(\E[\bar E^{(N)}(\infty)],  \E[\bar F^{(N)}(\infty)], A^{(N)}(\infty) \big)$ satisfies \eqref{rrp: fullcsv}--\eqref{rrp: frac}.

\begin{itemize}
\item Condition \eqref{rrp: fullcsv} was already verified in \eqref{eq:bar_fij}.

\item To check condition \eqref{rrp:prob1}, we fix $i\not=j,$ and use the test function $g(e,f)=e_{ij}$. Then
\begin{align*}
G^{(N)} g(e,f) =  Q_{ij}(e,f)\sum_{k=1}^{r} \mu_{ki} N f_{ki}\big((e_{ij} + 1/N) - e_{ij} \big) + \mu_{ij} N e_{ij} \big((e_{ij} - 1/N) - e_{ij} \big),
\end{align*}
where $Q_{ij}(e,f)$ is the probability that upon dropping a passenger off at region $i$, a car drives empty to region $j$ given the current state of the system is $(e,f)$. Using \eqref{bar-p} and the fact that $Q_{ij}(e,f) \in [0,1]$, we see that

\begin{align}
0&=\E \left[Q_{ij}\left(\bar E^{(N)}(\infty), \bar F^{(N)}(\infty)\right)\sum_{k=1}^r \mu_{ki}\bar{F}^{(N)}_{ki}(\infty)  -\mu_{ij}\bar{E}^{(N)}_{ij}(\infty)\right] \label{eq:partial_1} \\
&\leq \sum_{k=1}^r \mu_{ki}\E[\bar{F}^{(N)}_{ki}(\infty)] -\mu_{ij}\E[\bar{E}^{(N)}_{ij}(\infty)] \notag .
\end{align}

\item For condition \eqref{rrp:prob2}, we fix $i$ and use the test function $g(e,f) = e_{ii}$. Then
\begin{align*}
G^{(N)} g(e,f) =&\  N\lambda_{i} 1(e_{ii} > 0) \big((e_{ii}-1/N) - e_{ii} \big) \\
&+ \Big(Q_{ii}(e,f)\sum_{k=1}^{r} \mu_{ki} N f_{ki} + \sum_{k=1 , k \neq i}^{r} \mu_{ki} N e_{ki} \Big)\big((e_{ii}+1/N) - e_{ii} \big).
\end{align*}
Taking the expected value and using \eqref{bar-p}, we see that
\begin{align}
&\lambda_{i} \Prob(\bar{E}^N_{ii}(\infty) > 0) \notag \\
=&\  \E \bigg[Q_{ii}(\bar{E}^{(N)}(\infty),\bar{F}^{(N)}(\infty))\sum_{k=1}^{r} \mu_{ki}  \bar{F}_{ki}^{(N)}(\infty) + \sum_{k=1 , k \neq i}^{r} \mu_{ki}  \bar{E}_{ki}^{(N)}(\infty) \bigg] \label{eq:partial_2}
\end{align}
Using the fact that $Q_{ii}(e,f) \in [0,1]$, we conclude that
\begin{align*}
0&\leq \sum_{k=1}^r \mu_{ki}\E[\bar{F}^{(N)}_{ki}(\infty)] +\sum_{k=1,k\not=i}^r \mu_{ki}\E[\bar{E}^{(N)}_{ki}(\infty)]   -\lambda_i \Prob(\bar E_{ii}^{(N)}(\infty)>0),
\end{align*}
and
\begin{align*}
0&\geq \sum_{k=1,k\not=i}^r \mu_{ki}\E[\bar{E}^{(N)}_{ki}(\infty)] - \lambda_i \Prob(\bar E_{ii}^{(N)}(\infty)>0).
\end{align*}

\item Fix $i$. To check condition \eqref{rrp: carcsv}, we could use the test function $g(e,f) = \sum_{j=1}^{r} e_{ij}$. Alternatively, it is easier to just add up \eqref{eq:partial_1} for all $j \neq i$ together with \eqref{eq:partial_2} to arrive at
\begin{align*}
&\ \lambda_i \Prob(\bar E_{ii}^{(N)}(\infty)>0) +\sum_{j=1,j\not=i}^r \mu_{ij}\E[\bar{E}^{(N)}_{ij}(\infty)] \\
=&\ \sum_{k=1}^r \mu_{ki}\E[\bar{F}^{(N)}_{ki}(\infty)] +\sum_{k=1,k\not=i}^r \mu_{ki}\E[\bar{E}^{(N)}_{ki}(\infty)].
\end{align*}

\item Conditions \eqref{rrp: carnum} and \eqref{rrp: frac} hold trivially.
\end{itemize}
\endproof

\subsection{Proof of Lemma~\ref{lem:equil}}
\proof{Proof of Lemma~\ref{lem:equil}.}
Substituting \eqref{eq:one} and \eqref{eq:two} into \eqref{eq:three}, we obtain
\begin{align*}
\lambda_{i} \bar a_i = \sum_{\substack{\ell = 1 \\ \ell \neq i}}^{r} Q_{\ell i}\sum_{k = 1}^{r} \mu_{k \ell }\bar f_{ k \ell} + Q_{ii}\sum_{k = 1}^{r} \mu_{k i} \bar f_{k i} =&\ \sum_{ \ell = 1}^{r} Q_{\ell i} \sum_{k = 1}^{r} \lambda_{k} P_{k \ell} \bar a_k \\
=&\ \sum_{k = 1}^{r} \lambda_{k} \bar a_k \sum_{\ell = 1}^{r} P_{k \ell} Q_{\ell i}, \quad 1 \leq i \leq r.
\end{align*}
In matrix form, these equations can be written as
\begin{align}
\big( I - B \big) \Lambda \bar a = 0, \label{eq:traffic}
\end{align}
where $I$ is the $r \times r$ identity matrix, and  $B$ and $\Lambda$ are $r\times r$ matrices defined as
\begin{align}
B_{ij} = \sum_{\ell = 1}^{r} P_{j \ell} Q_{\ell i}, \quad \text{ and } \quad  \Lambda = \text{ diag}(\lambda).  \label{eq:QL}
\end{align}
Observe that $B$ is a column stochastic matrix, i.e.\ columns sum to one. We now argue that $B$ is irreducible because the CTMC $\big(E^{(N)}, F^{(N)}\big)$ is.  For any $1 \leq i,j \leq r$, the entry $B_{ij}$ is the probability that a car picks up a passenger at region $i$, drives him to some region $k$, and then drives empty to region $j$ to wait for a new passenger there (or stay and wait at region $j$ if $k = j$). Therefore, $B$ can be interpreted as the transition probability matrix of a discrete-time Markov chain (DTMC) that describes the motion of a single car in a network, i.e.\ how it serves passengers and makes routing decisions, as if travel times were zero. Irreducibility of $B$ then means that starting from any region the car in the DTMC can visit any other region, which is clearly satisfied when $P_{ij} > 0$ and $Q_{ii} > 0$ for all $i,j= 1, \ldots, r.$

Since $B$ is column stochastic and irreducible, \eqref{eq:traffic} has a unique solution up to a multiplicative constant. That is, any solution to \eqref{eq:traffic} must be of the form $c a^* \geq 0$, where $c > 0$, and $a^* \geq 0$ is a unique vector in $\R_+^{r}$. We now argue that \eqref{eq:ycomp} and \eqref{eq:balance} uniquely define $c$. First, we use $\eqref{eq:balance}$ and \eqref{eq:one}--\eqref{eq:two} to write 
\begin{align*}
1 =&\ \sum_{i=1}^{r} \sum_{j=1}^{r} \bar f_{ij} + \sum_{i=1}^{r} \sum_{\substack{j=1\\ j \neq i}}^{r} \bar e_{ij} + \sum_{i=1}^{r} \bar e_{ii} \\
=&\ \sum_{i=1}^{r} \sum_{j=1}^{r} \frac{\lambda_{i} P_{ij} }{\mu_{ij}}\bar a_i + \sum_{i=1}^{r} \sum_{\substack{j=1\\ j \neq i}}^{r} \frac{Q_{ij}}{\mu_{ij}}  \sum_{k = 1}^{r} \mu_{ki}\bar f_{ki} + \sum_{i=1}^{r} \bar e_{ii} \\
=&\ \sum_{i=1}^{r} \sum_{j=1}^{r} \frac{\lambda_{i} P_{ij} }{\mu_{ij}}\bar a_i + \sum_{i=1}^{r} \sum_{\substack{j=1\\ j \neq i}}^{r} \frac{Q_{ij}}{\mu_{ij}}  \sum_{k = 1}^{r} \lambda_{k} P_{ki}\bar a_k + \sum_{i=1}^{r} \bar e_{ii}.
\end{align*}
The equation above can be written as
\begin{align*}
1 =&\ \sum_{i=1}^{r} \tilde c_{i} \bar a_i + \sum_{i=1}^{r} \bar e_{ii} = c\sum_{i=1}^{r} \tilde c_{i} a_i^{*} + \sum_{i=1}^{r} \bar e_{ii}
\end{align*}
where $\tilde c_1, \ldots, \tilde c_{r} > 0$, and in the second equality we used $\bar a = c a^{*}$. Now if 
\begin{align}
\frac{1}{\sum_{i=1}^{r} \tilde c_{i} a_i^{*}} a^{*}_{i} \leq 1, \quad 1 \leq i \leq r, \label{eq:semi_final}
\end{align}
then $c = 1/{\sum_{i=1}^{r} \tilde c_{i} a_i^{*}}$ and $\bar e_{ii} = 0$ for all $i$ are the unique choices under which both \eqref{eq:semi_final} and \eqref{eq:ycomp} hold, and we are done. 

Now suppose \eqref{eq:semi_final} is violated. We cannot choose $c = 1/{\sum_{i=1}^{r} \tilde c_{i} a_i^{*}}$, because that would violate the requirement that $\bar a \in [0,1]^{r}$. Instead, choose 
\begin{align}
c = \sup \{s > 0 : s a^{*} \in [0,1]^{r}\}, \label{eq:csup}
\end{align}
and observe that such a choice necessarily satisfies $c < 1/{\sum_{i=1}^{r} \tilde c_{i} a_i^{*}}$. To satisfy \eqref{eq:ycomp}, we must set $\bar e_{ii} = 0$ for all $i$ such that $ \bar a_{i} < 1$. The only restriction on $\bar e_{ii}$ for $i$ such that $\bar a_{i} = 1$ is that 
\begin{align*}
\sum_{i: \bar a_{i}=1}^{r} \bar e_{ii} = 1 - c\sum_{i=1}^{r} \tilde c_{i} a_i^{*}.
\end{align*}
It remains to verify that the only viable choice of $c$ is given by \eqref{eq:csup}. Choosing $c < \sup \{s > 0 : s a^{*} \in [0,1]^{r}\}$ implies $\bar a \in (0,1)^{r}$, and consequently, \eqref{eq:ycomp} implies that $\bar e_{ii} = 0$ for all $i$. In such a case, 
\begin{align*}
c\sum_{i=1}^{r} \tilde c_{i} a_i^{*} + \sum_{i=1}^{r} \bar e_{ii} = c\sum_{i=1}^{r} \tilde c_{i} a_i^{*} < 1,
\end{align*}
because $\sup \{s > 0 : s a^{*} \in [0,1]^{r}\} < 1/{\sum_{i=1}^{r} \tilde c_{i} a_i^{*}}$, and a contradiction is reached. Choosing $c > \sup \{s > 0 : s a^{*} \in [0,1]^{r}\}$ forces some element of $\bar a$ to be greater than one, which violated the requirement that $\bar a \in [0,1]^{r}$. This concludes the proof.
\endproof
\ACKNOWLEDGMENT{The authors thank Siddhartha Banerjee and Daniel Freund for feedback and stimulating discussion on this work. They also thank Peter Frazier for arranging a visit to Uber headquarters, where they received invaluable feedback. This research is supported in part by NSF Grants CNS-1248117, CMMI-1335724, and CMMI-1537795.}

\bibliographystyle{informs2014} 
\bibliography{dai20170216} 



\end{document}